\documentclass[preprint, 3p, times]{elsarticle}

\usepackage{amsfonts,amsmath,amssymb,amsthm}
\usepackage{graphicx}

\usepackage{algorithm}
\usepackage{algorithmic}

\newcommand{\RR}{\mathbb{R}}

\newcommand{\ZZ}{\mathbb{Z}}

\newtheorem{definition}{Definition}[section]

\newtheorem{theorem}{Theorem}[section]

\newtheorem{problem}{Problem}[section]
\newtheorem{remark}{Remark}[section]

\journal{Mathematics and Computers in Simulation}

\begin{document}

\begin{frontmatter}



\title{Local interpolation schemes for landmark-based image registration:\\ a comparison}


\author{Giampietro Allasia} 
\ead{giampietro.allasia@unito.it}

\author{Roberto Cavoretto\corref{cor1}} 
\ead{roberto.cavoretto@unito.it}
\cortext[cor1]{Corresponding author.}

\author{Alessandra De Rossi}
\ead{alessandra.derossi@unito.it}

\address{Department of Mathematics \lq\lq G. Peano\rq\rq, University of Torino, via Carlo Alberto 10, I--10123 Torino, Italy}

\begin{abstract}
In this paper we focus, from a mathematical point of view, on properties and performances of some local interpolation schemes for landmark-based image registration. Precisely, we consider modified Shepard's interpolants, Wendland's functions, and Lobachevsky splines. They are quite unlike each other, but all of them are compactly supported and enjoy interesting theoretical and computational properties. In particular, we point out some unusual forms of the considered functions. Finally, detailed numerical comparisons are given, considering also Gaussians and thin plate splines, which are really globally supported but widely used in applications.
\end{abstract} 


\end{frontmatter}

\section{Introduction}
Image registration is an important challenging topic in image processing. It consists mainly in finding a suitable transformation between two images (or image data), called \textsl{source} and \textsl{target images}, taken either at different times or from different sensors or viewpoints. The scope is to determine a transformation such that the transformed version of the source image is similar to the target one. There is a large number of applications demanding image registration, including astronomy, biology, computer vision, genetics, physics, medicine, robotics, to name a few. For an overview, see e.g. 
\cite{fischer08, fitzpatrick08, goshtasby05, hajnal01, maintz98, Modersitzki04, Modersitzki09, pluim03, rohr01, scherzer06, zitova03} 
and references therein. In medicine, for example, registration is required for combining different modalities (X-ray, computer tomography (CT), magnetic resonance imaging (MRI) and positron emission tomography (PET) images, for instance), monitoring of diseases, treatment validation, comparison of the patient's data with anatomical atlases, and radiation therapy. In particular, the \textsl{landmark-based image registration} process is based on two finite sets of landmarks, i.e. sparse data points located on images, usually not uniformly distributed, where each landmark of the source image has to be mapped onto the corresponding landmark of the target image (see \cite{Modersitzki04, Modersitzki09, rohr01}). Now, in order to give a more formal idea, we consider the sets ${\cal S}_N=\{ {\bf x}_j \in \RR^m,j=1,2,\ldots,N\}$ and ${\cal T}_N=\{ {\bf t}_j\in\RR^m,j=1,2,\ldots,N\}$  each containing $N$ point-landmarks in the source and target images, respectively. Thus, the registration problem involves a transformation $\textbf{F}:\RR^m \rightarrow \RR^m$, such that
\begin{equation}
 \textbf{F}({\textbf{x}_j})= {\textbf{t}_j}, \hspace{1cm} j=1,2,\ldots,N,	\nonumber
\end{equation}
where each coordinate $F_k$ of the transformation function $\textbf{F}=(F_1,F_2,\ldots,F_m)^T$ is separately calculated, that is, the interpolation problem involving $F_k: \RR^m\rightarrow \RR$ is solved for $k=1,2,\ldots,m$, with the corresponding conditions $F_k({\bf x}_j)={t}_{jk}$, $j=1,2,\ldots,N$. This problem can be formulated in the context of multivariate scattered data interpolation, and solved by different techniques, among which radial basis functions (RBFs) play a preminent role (see, e.g., \cite{Buhmann03,Iske03,Wendland05}). The use of RBF transformations, in particular of the thin plate splines,  for point-based image registration was first proposed by Bookstein \cite{bookstein89}, and it is still common (see \cite{Quatember10a} and the software package MIPAV \cite{mipav}). A number of authors have investigated the most popular radial basis function transformations in the image registration context: thin plate spline 
\cite{arad94,little97b},
multiquadric \cite{little97a,ruprecht93}, inverse multiquadric \cite{ruprecht93}, and Gaussian transformations \cite{arad94}. A more specific application which involves registration and includes imaging techniques, such as computer tomography and magnetic resonance imaging, can be found in \cite{Quatember10a,Quatember12}.

Since using globally supported RBFs, as for example the Gaussians, a single landmark pair change may influence the whole registration result, in the last two decades several methods have been presented to circumvent this disadvantage, such as weighted least squares and weighted mean methods (WLSM and WMM, respectively) \cite{goshtasby88}, compactly supported  radial basis functions (CSRBFs), especially Wendland's and Gneiting's functions \cite{Cavoretto-DeRossi12,Cavoretto-DeRossi13,fornefett01}, and elastic body splines (EBSs) \cite{kohlrausch05}. 

A certain number of papers have been dedicated to recall and compare these methods for nonrigid image registration: main contributions, advantages and drawbacks of radial basis functions, compactly supported radial basis functions and elastic body splines are mentioned in \cite{zitova03}; thin plate splines, multiquadrics, piecewise linear and weighted mean transformations are explored and their performances are compared in \cite{zagorchev06}; finally, radial basis functions, Wendland's functions and elastic body splines are reviewed, together with B-splines and wavelets, in \cite{holden08}. 

Several authors have shown the superiority of local registration methods over the global ones in some situations, for instance, in medical imaging and in airborne imaging. In fact, a global mapping cannot properly handle images locally deformed. For this reason, more recently, local methods, already known in interpolation theory, have been proposed in landmark-based image registration: the modified Shepard's method (also known as the inverse distance weighted method (IDWM)) \cite{Cavoretto-DeRossi08a,Cavoretto-DeRossi08b}, and Lobachevsky spline method \cite{Allasia10,Allasia11b}. These interpolation techniques, giving rise to local mappings, handle well images locally deformed. Moreover, they are in general stable and the computational effort to determine transformations is low and, therefore, a large number of landmarks can be used.

In this paper we focus, from a mathematical point of view, on properties and performances of some local interpolation schemes for landmark-based image registration. Precisely, we consider modified Shepard's interpolants, Wendland's functions, and Lobachevsky splines. They are quite unlike each other, but all of them are compactly supported and enjoy interesting theoretical and computational properties (see \cite{Wendland05,Cavoretto10}). In particular, we point out some unusual forms of the considered functions. Moreover, referring to Wendland's functions, we consider for the first time in this context, as far as we know, compactly supported interpolants given by products of univariate Wendland's functions \cite{CDSV10}. All these methods are also compared with Gaussians and thin plate splines, which are globally supported but are still among the most widely used methods in applications. 

Numerical experiments point out differences in accuracy and smoothness of the considered methods. The comparison can be useful to users in the choice of the appropriate transformation for their scopes. Moreover, since some schemes need parameters, our numerical tests might be of interest in the choice of them.

The paper is organized as follows. Section 2 introduces some preliminaries: the landmark-based registration problem and the solvability of the associated interpolation problem. In Section 3 we briefly recall radial basis functions, like Gaussians, multiquadrics, inverse multiquadrics and thin plate splines to construct globally supported transformations. Section 4 is devoted to describe local transformations, given by the modified Shepard's formula which uses RBFs as local approximants. In Section 5 Wendland's functions are presented to define compactly supported transformations, whereas in Section 6 we focus on Lobachevsky splines which define again compactly supported transformations. Finally, Section 7 contains several numerical results obtained in some test and real-life examples: special emphasis is devoted to comparing accuracy of local interpolation schemes and to determining optimal values of parameters. 


\section{Preliminaries}
Let ${\cal S}_N=\{ {\bf x}_j\in \RR^m,j=1,2,\ldots,N\}$ be a given set of landmarks in the source image $S$ and let ${\cal T}_N=\{ {\bf t}_j\in \RR^m,j=1,2,\ldots,N\}$ be the given set of corresponding landmarks in the target image $T$. The registration problem reads as follows.

\begin{problem} \label{lit}
Let the landmark sets ${\cal S}_N$ and ${\cal T}_N$ be given. Find a transformation ${\bf F}:\RR^m \rightarrow \RR^m$ within a suitable space ${\cal F}$ of admissible functions, such that
\begin{eqnarray}
 \label{interp}
 {\bf F}({\bf x}_j)= {\bf t}_j, \hspace{1cm} j=1,2,\ldots,N.
\end{eqnarray}
\end{problem}

Each coordinate $F_k$ of the transformation function is calculated separately, i.e. the interpolation problem $F_k: \RR^m\rightarrow \RR$ is solved for each $k=1,2,\ldots,m$, with the corresponding conditions 
\begin{equation}
F_k({\bf x}_j)={t}_{jk}, \hspace {1.cm} j=1,2,\ldots,N. \nonumber
\end{equation}

In order to have a class of basis functions that generate non-singular interpolation matrices for any set of distinct points, we introduce the concept of strictly positive definite functions \cite{Fasshauer07}. We suppose that the interpolant $F_k:\RR^m\rightarrow\RR$ has the form
\begin{equation}\label{intform}
F_k({\bf x})=\sum_{j=1}^N c_j \Psi({\bf x}-{\bf x}_j), 
\end{equation}
$c_j$ being the coefficients to be found. A necessary condition to have unique solvability of the interpolation problem 
\begin{eqnarray}\label{IC}
F_k({\bf x}_j)=t_{jk}, \hspace{0.5cm} j=1,2,\ldots,N,
\end{eqnarray}
is given by the following result \cite{Fasshauer07}.

\begin{theorem} \label{nec_cond}
The interpolation problem (\ref{IC}), where $F_k$ is of the form (\ref{intform}), has a unique solution if the function $\Psi$ is strictly positive definite on $\RR^m$, that is 
\begin{eqnarray}
	\label{real_positive_definite_function}
	\sum_{i=1}^{N}\sum_{j=1}^{N} c_i c_j \Psi(\textbf{x}_i - \textbf{x}_j) > 0,
\end{eqnarray}
for any $N$ pairwise different points $\textbf{x}_1,\textbf{x}_2,\ldots,\textbf{x}_N \in \RR^m$, and $\textbf{c}=\left[c_1,c_2,\ldots,c_N\right]^T\in \RR^N$, $\textbf{c}\neq \textbf{0}$. 
\end{theorem}

Moreover, we remark that Theorem \ref{nec_cond} is also satisfied for a strictly conditionally positive definite function $\Psi$ of order $v$ if the quadratic form (\ref{real_positive_definite_function}) holds and
\begin{equation}
	\sum_{i=1}^{N} c_i p({\bf x}_i) = 0,	\nonumber
\end{equation}
for any polynomial $p$ of degree at most $v-1$. 

Finally, we give the following theorem which allows us to construct multivariate strictly positive definite functions from univariate ones (see, e.g., \cite{Wendland05}).

\begin{theorem} \label{spdf_dD}
Suppose that $\psi_1,\psi_2, \ldots, \psi_m$ are strictly positive definite and integrable functions on $\RR$, then  
\begin{eqnarray}
 \Psi({\bf x}) = 	\psi(x_1)\psi(x_2)\cdots \psi(x_m), \hspace{1cm} \textbf{x}=(x_1,x_2,\ldots,x_m) \in \RR^m, \nonumber
\end{eqnarray}
is a strictly positive definite function on $\RR^m$.
\end{theorem}

\section{Radial basis functions} \label{g_tps}

In this section we consider globally supported radial basis functions, which are well-known in the field of approximation theory and scattered data interpolation (see \cite{Fasshauer07,Wendland05}), and widely used in landmark-based image registration as well (see, e.g., \cite{Modersitzki09} and references therein).

Using radial basis functions, the general coordinate $F_k({\bf x})$, $k=1,2,\ldots m$ of the transformation function is assumed to have the form 
\begin{equation}
\label{s}
	F_k({\bf x})= \sum_{j=1}^N a_j \Phi(\vert\vert{\textbf{x}}-{\textbf{x}_j}\vert\vert_2)+\sum_{k=1}^U b_k \pi_k({\textbf{x}}) ,
\end{equation} 
where $\Phi(\vert\vert \textbf{x}-{\textbf{x}_j}\vert\vert_2)$ is a radial basis function depending only on the Euclidean distance $r=\vert\vert \textbf{x} -{\textbf{x}}_j \vert\vert_2$, and $a_j$ and $b_k$ are coefficients to be determined. The space ${\cal P}_{v-1}^m$ $\equiv $ ${\cal P}_{v-1}({\RR}^m)$ $=$ ${\rm span}$$\{ \pi_k\}_{k=1}^U$, where the $\pi_k$ are a basis of polynomials up to degree $v-1$, has dimension $U=(m+v-1)!/(m!(v-1)!)$, which must be lower than $N$. Therefore, in order to compute the coefficients ${\bf a}=(a_1,a_2,\ldots,a_N)^T$ and ${\bf b}=(b_1,b_2,\ldots,b_U)^T$ in (\ref{s}), it is required to solve the following system of linear equations
\begin{equation}
\label{sist}
\left\{
\begin{array}{rcl}
{\textbf{M}\textbf{a}} + {\textbf{Q}\textbf{b}} &=& {\textbf{t}}, \\
{\textbf{Q}}^T {\textbf{a}} &=& {\textbf{0}} ,
\end{array}
\right.
\end{equation}
where ${\textbf{M}}=\{\Phi(||{\textbf{x}_i}-{\textbf{x}_j}||)\}$ is a $N\times N$ matrix, ${Q}=\{ \pi_k({\textbf{x}_j}) \}$ is a $N\times U$ matrix, and ${\bf t}$ denotes the column vector of the $k$-th coordinate of the target point-landmarks $\textbf{t}_j$ corresponding to the image point-landmarks $\textbf{x}_j$. Equations (\ref{sist}) are obtained by requiring that $\textbf{F}$ satisfies the interpolation conditions (\ref{interp}) and the side conditions $\sum_{j=1}^N a_j\pi_k({\textbf{x}}_j)=0$, for $k=1,2,\ldots,U$, i.e. ${\bf Q}^T {\bf a} = {\bf 0}$.

The most popular choices for $\Phi$ in landmark-based registration are
\begin{equation}
\left.
\begin{array}{rclll}
\Phi(r) & = & r^{2}\log r,      &                       & \hspace{0.5cm} \mbox{{\rm (thin plate spline)}} \nonumber \\
\Phi(r) & = & {\rm e}^{-\alpha^2 r^2}, &                      & \hspace{0.5cm} {\rm (Gaussian)} \\
\Phi(r) & = & (r^2+\gamma^2)^{\mu/2}, &                      & \hspace{0.5cm} \mbox{{\rm (generalized multiquadric)}} \\
\end{array}
\right.
\end{equation}
where $\alpha, \gamma \in \RR^+$ and $\mu \in \ZZ$. The Gaussian (G) and the inverse multiquadric (IMQ), which occurs for $\mu<0$ in the generalized multiquadric function, are strictly positive definite functions, whereas the thin plate spline (TPS) and the multiquadric (MQ), i.e. for $\mu>0$ in the generalized multiquadric function, are strictly conditionally positive definite functions of order $\mu$. The addition of a polynomial term of a certain order along with side conditions, in order to guarantee existence and uniqueness of the solution in the linear equation system (\ref{sist}), is required only for strictly conditionally positive definite functions. It allows us to have a nonsingular interpolation matrix. For thin plate spline the order of the polynomial is 1 and for multiquadric depends on the exponent $\mu$, the minimal degree being $v=\mu-1$. Polynomials have global support and therefore the polynomial part influences globally the registration result.

An important feature of some radial basis functions is the presence of a shape parameter, which allows us to control their influence on the registration result \cite{fornefett01}. On the other hand this property could be seen as a drawback, because the user needs to give the parameter value. Also for this reason the thin plate spline is usually considered more suitable than other RBFs for image registration, since this process of image registration must sometimes be automatic. 

The thin plate spline seems to be preferable to other radial basis functions for image registration also for other reasons (see, in particular, \cite{arad94} and references already cited). The thin plate spline minimizes a functional which represents the bending energy of a thin plate separately for each component $F_k$, $k=1,2,\ldots,m$, of the transformation $\textbf{F}$. Thus, the functional $J(\textbf{F})$ can be separated into a sum of similar functionals that only depend on one component $F_k$ of \textbf{F}, and the  problem of finding \textbf{F} can be decomposed into $m$ problems \cite{wahba90}.

Note that the interpolation matrices generated by radial basis functions are dense, since they are globally supported, and ill-conditioned \cite{Fasshauer07}, especially those generated by Gaussian. Ill-conditioning could happen even if it is required to interpolate a relatively small number of landmarks, since the landmarks may be very close to each other. 

Another possible disadvantage of the use of RBFs is given by the number of floating-point operations, which can be very time consuming, especially if the number of landmarks is high or the registration of 3D images is needed. 


\section{Modified Shepard's formula} 

Shepard's method or inverse distance weighted method is commonly used for multivariate interpolation and approximation of scattered data. In the literature, many versions of this method can be found (see \cite{Cavoretto10} for an overview). Here, we consider a modified version of Shepard's method, which is local and uses radial basis functions as local approximants (see \cite{Allasia11}). This interpolant has been already successfully used in landmark-based registration context \cite{Cavoretto-DeRossi08a,Cavoretto-DeRossi08b}.

Referring to the context of landmark-based image registration, we define the $k$-th coordinate of a modified Shepard's transformation $\textbf{F}:\RR^m \rightarrow \RR^m$ in the form
\begin{eqnarray}
	\label{sh}
	F_k(\textbf{x})=\sum_{j=1}^{N} L_j(\textbf{x}) \bar{W}_j(\textbf{x}), \hspace{0.5cm} k=1,2,\ldots,m.
\end{eqnarray}
The {\sl nodal functions} $L_j$, $j=1,2,\ldots,N$, are local approximants at $\textbf{x}_j$, constructed on the $N_L$ source landmarks closest to $\textbf{x}_j$ and satisfying the interpolation conditions $L_j(\textbf{x}_j)=\textbf{t}_j$. The weight functions $\bar{W}_j$, $j=1,2,\ldots,N$, are given by
\begin{equation}\label{wei1}
	\bar{W}_j(\textbf{x}) = \frac{W_j(\textbf{x})}{\sum_{k=1}^{N} W_k(\textbf{x})}, \hspace{0.5cm} j=1,2,\ldots,N,\nonumber
\end{equation}
where
\begin{equation} \label{wei2}
	W_j(\textbf{x}) = \tau(\textbf{x},\textbf{x}_j) / \alpha(\textbf{x},\textbf{x}_j), \nonumber
\end{equation}
$\alpha(\textbf{x},\textbf{x}_j) = \Vert \textbf{x}-\textbf{x}_j \Vert_2^2$, and $\tau(\textbf{x},\textbf{x}_j)$ is a non-negative localizing function defined by
\begin{equation}
\tau(\textbf{x},\textbf{x}_j)= \left\{
\begin{array}{ll}
1, & \mbox{if $\textbf{x} \in {\cal C}(\textbf{x}_j;\rho)$}, \\ \nonumber
0, & \mbox{otherwise}, 
\end{array}
\right.
\end{equation}
${\cal C}(\textbf{x}_j;\rho)$ being the hypercube of centre at $\textbf{x}_j$ and side $\rho$. Note that the choice of the hypercube side $\rho$ may not be unique, being partially empirical.

In this way, each component $F_k$, $k=1,2,\ldots,m$, is evaluated at $\textbf{x}$ considering only a certain number $N_W$ of landmarks closest to $\textbf{x}$. 

\begin{remark}
The weights $\bar{W}_j (\textbf{x})$ given in (\ref{sh}) satisfy the cardinality relations
\begin{equation}
\bar{W}_j (\textbf{x}_i) = \delta_{ij}, \ \ \ \ \ i,j=1,2,\ldots,N, \nonumber
\end{equation}
where $\delta_{ij}$ is the Kronecker delta, and give a partition of unity, namely
\begin{equation}
\sum_{j=1}^{N} \bar{W}_j (\textbf{x}) = 1. \nonumber
\end{equation}
\end{remark}

In (\ref{sh}) we take the local approximants $L_j$, $j=1,2,\ldots,N$, as local RBF interpolants, of the form
\begin{eqnarray}
\label{rbf}
  L_j(\textbf{x}) = \sum_{i=1}^{N_L}  a_i \Phi(\left\|\textbf{x}-\textbf{x}_i\right\|_2)+\sum_{k=1}^U  b_k \pi_k(\textbf{x}) ,
\end{eqnarray}
where the radial basis functions $\Phi(\left\|{\bf x}-\textbf{x}_i\right\|_2)$ depend only on the $N_L$ source landmark points of the considered neighborhood of $\textbf{x}_j$, the index $i$ refers to and renumbers the landmark points of the neighborhood of $x_j$, and the space ${\cal P}_{v-1}^m$ spanned by the ($v-1$)-degree polynomials $\pi_k (x)$ has a dimension $U=(m+v-1)!/(m!(v-1)!)$ which must be lower than $N_L$. It is required that $L_j$ satisfies the interpolation conditions
\begin{equation}
 L_j(\textbf{x}_i)=\textbf{t}_i, \hspace{0.5cm} i=1,2,\ldots,N_L ,	\nonumber
\end{equation}
and the side conditions 
\begin{equation}
\sum_{i=1}^{N_L} a_i\pi_k(\textbf{x}_i)=0, \hspace{.5cm} \hbox{for } k=1,2,\ldots,U. \nonumber
\end{equation}
Hence, to compute the coefficients $\textbf{a}=(a_1,a_2,\ldots,a_{N_L})^T$ and $\textbf{b}=(b_1,b_2,\ldots,b_U)^T$ in (\ref{rbf}), it is required to solve uniquely the system of linear equations
\begin{equation}
\left\{
\begin{array}{rcl}
{\tilde{\textbf M}{\textbf a}} + {\tilde{\textbf Q}{\textbf b}} & = & \textbf{t}, \\ \nonumber
{\tilde{\textbf Q}}^T \textbf{a} & = & \textbf{0}, 
\end{array}
\right.
\end{equation}
where $\tilde{{\textbf M}}=\{ \Phi(||\textbf{x}_j-\textbf{x}_i||_2) \}$ is a $N_L\times N_L$ matrix, $\tilde{{\textbf Q}}=\{ \pi_k(\textbf{x}_j) \}$ is a $N_L\times U$ matrix, and $\textbf{t}$ denotes the column vector of the target landmark points $\textbf{t}_j$ corresponding to the $\textbf{x}_j$.

Registration results obtained using local radial basis functions (in particular thin plate spline and Gaussian) point out that the local transformation technique compares well with the global one, and in some cases it is superior. In particular, the former is preferable when the number of landmarks is large, since its main features are locality and stability \cite{Cavoretto-DeRossi08a,Cavoretto-DeRossi08b}. The drawback of the technique is the need of a empirical determination of the local parameters $N_L$ and $N_W$.


\section{Wendland's functions}
In this section we focus on probably the most known family of compactly supported radial basis functions, namely Wendland's functions (see \cite{Fasshauer07}). They have been introduced in image registration context with the motivation that their influence around a landmark is limited, in 2D and 3D images on a circle or a sphere, respectively \cite{fornefett99,fornefett01}. This property allows us the registration of medical images where changes occur only locally. Here we propose also the use of compactly supported radial basis functions defined by products of univariate Wendland's functions. Some numerical tests using this type of functions, have pointed out their better accuracy but worse stability in comparison with multivariate Wendland's functions \cite{CDSV10}, the support being a square or a cube, in 2D or 3D setting.

Wendland's functions are obtained by using the truncated power function
$\varphi_s(r) = (1-r)_+^s$ (where $\left(x\right)_+$ is defined as $x$ for
$x > 0$ and $0$ for $x \leq 0$), which is strictly positive definite and
radial on $\RR^m$ for $s\geq \left\lfloor m/2\right\rfloor +1$, and
repeatedly applying the operator ${\cal I}$ given by $\left({\cal
I}\varphi\right) (r)=\int_{r}^{\infty} t\varphi(t) dt$, $r \geq 0$ (see,
e.g., \cite{Wendland05}). Then, Wendland's functions are given by
\begin{equation}
        \varphi_{m,h}={\cal I}^h \varphi_{\left\lfloor m/2\right\rfloor+h+1},
\nonumber
\end{equation}

Specifically, they are all supported on the interval $[0,1]$ and have a polynomial representation there. In addition to this, the following theorem states that any other compactly supported $C^{2h}$ polynomial function that is strictly positive and radial on $\RR^s$ will not have a smaller polynomial degree. Finally, they have minimal degree with respect to a given space dimension $m$ and smoothness $2h$ \cite{Wendland05}. 

\begin{theorem}
The functions $\varphi_{m,h}$ are strictly positive definite and radial on $\RR^m$ and are of the form
\begin{equation}
 \varphi_{m,h}(r) = \left\{
\begin{array}{ll}
p_{m,h}(r), & \mbox{ } \  r \in [0,1],\\ \nonumber
0, & \mbox{ } \  r > 1,  \nonumber
\end{array}
\right.
\end{equation}
with a univariate polynomial $p_{m,h}$ of degree $\left\lfloor m/2\right\rfloor+3h+1$. Moreover, $\varphi \in C^{2h}(\RR)$ are unique up to a constant factor, and the polynomial degree is minimal for given space dimension $m$ and smoothness $2h$.
\end{theorem}

Though there exist recursive formulas to compute the functions
$\varphi_{m,h}$ for all $m$ and $h$, here for convenience we only give
explicit forms of $\varphi_{m,h}$, for $h=0,1,2,3$.

\begin{theorem}
 The functions $\varphi_{m,h}$, $h=0,1,2,3$, have the form
\begin{equation}
 \left.
\begin{array}{l}
\varphi_{m,0}(r) \doteq \left( 1-r\right)_+^s ,\\ \nonumber
\varphi_{m,1}(r) \doteq \left( 1-r\right)_+^{s+1} \left[\left(s+1\right)r+1\right],\\ \nonumber
\varphi_{m,2}(r) \doteq \left( 1-r\right)_+^{s+2}\left[\left(s^2+4s+3\right)r^2+\left(3s+6\right)r+3\right],\\ \nonumber
\varphi_{m,3}(r) \doteq \left( 1-r\right)_+^{s+3} \left[\left(s^3+9s^2+23s+15\right)r^3+\left(6s^2+36s+45\right)r^2+\left(15s+45\right)r+15\right], \nonumber
\end{array}
\right.
	\end{equation}
where $s=\left\lfloor m/2\right\rfloor+h+1$, and the symbol $\doteq$ denotes equality up to a positive constant factor.
\end{theorem}

Since Wendland's functions are compactly supported, the interpolation matrices can be made sparse by appropriately scaling the support of the basic function. In the following we can consider only Wendland's functions depending on a shape parameter $c \in \RR^+$. Using for simplicity the symbol $\varphi_{m,h}$ as in Theorem 5.2, we list some of the most commonly used functions for $m=1$ along with their degree of smoothness $2h$, that is
\begin{equation}
\left.
\begin{array}{rclll}
\varphi_{1,0}(r) & \doteq & \left( 1-cr\right)_+,                                  &  & \hspace{0.5cm} {\rm C^0} \nonumber \\
\varphi_{1,1}(r) & \doteq & \left( 1-cr\right)_+^{3} \left(3cr+1\right),              &  & \hspace{0.5cm} {\rm C^2} \\ 
\varphi_{1,2}(r) & \doteq & \left( 1-cr\right)_+^{5}\left(8(cr)^2+5cr+1\right),        &  & \hspace{0.5cm} {\rm C^4} \\ 
\varphi_{1,3}(r) & \doteq & \left( 1-cr\right)_+^{7} \left(21(cr)^3+19(cr)^2 + 7cr+1\right),&  & \hspace{0.5cm} {\rm C^6} \\ 
\end{array}
\right.
\end{equation}
and for $m=2$, that is 
\begin{equation}
\left.
\begin{array}{rclll}
\varphi_{2,0}(r) & \doteq & \left( 1-cr\right)_+^2,                                  &  & \hspace{0.5cm} {\rm C^0} \nonumber \\
\varphi_{2,1}(r) & \doteq & \left( 1-cr\right)_+^{4} \left(4cr+1\right),              &  & \hspace{0.5cm} {\rm C^2} \\ 
\varphi_{2,2}(r) & \doteq & \left( 1-cr\right)_+^{6}\left(35(cr)^2+18cr+3\right),        &  & \hspace{0.5cm} {\rm C^4} \\ 
\varphi_{2,3}(r) & \doteq & \left( 1-cr\right)_+^{8} \left(32(cr)^3+25(cr)^2 + 8cr+1\right).&  & \hspace{0.5cm} {\rm C^6} \\ 
\end{array}
\right.
\end{equation}
We remark that the functions $\varphi_{2,k}$, $k=0,1,2,3$, are strictly positive definite and radial not only on $\RR^2$ but also on $\RR^m$, for $m \leq 3$ (see \cite{Wendland05}).

Referring to the image registration context we can define Wendland's transformations in two ways:
\begin{itemize}
	\item by multivariate radial Wendland's functions;
	\item by products of univariate radial Wendland's functions.
\end{itemize}

\begin{definition} \label{wetr}
Given a set of source landmark points ${\cal S}_N=\{ {\bf x}_j\in \RR^m,j=1,2,\ldots,N\}$, with associated the corresponding set of target landmark points ${\cal T}_N=\{ {\bf t}_j\in \RR^m,j=1,2,\ldots,N\}$, a {\sl Wendland's transformation} ${\bf F}:$ $ {\mathbb R}^m\rightarrow {\mathbb R}^m$ is such that each its component 
\begin{equation}
F_k:\RR^m\rightarrow \RR, \hskip0.5cm k=1,2,\ldots ,m, \nonumber
\end{equation}
assumes one of the following forms
\begin{eqnarray}
	\label{lsttr}
i) \ \	 F_k({\bf x}) &=& F_k(x_1,x_2,\ldots ,x_m)= \sum_{j=1}^N c_{j} 
\varphi_{m,h}(r),\\
\label{lsttr2}
ii) \ \	 F_k({\bf x}) &=& F_k(x_1,x_2,\ldots ,x_m)= \sum_{j=1}^N c_{j} \varphi_{1,h}(x_1-x_{j1}) \varphi_{1,h}(x_2-x_{j2})\ldots \varphi_{1,h}(x_m-x_{jm}),
\end{eqnarray}
and $(x_1,x_2,\ldots ,x_m)$, $(x_{j1},x_{j2},\ldots,x_{jm}) \in \RR^m$.
\end{definition} 
From Definition \ref{wetr} it follows that the transformation function $F_k:$ ${\mathbb R}^m\rightarrow {\mathbb R}$ is calculated for each $k=1,2,\ldots,m$, and the coefficients $c_{j}$ in (\ref{lsttr}) and (\ref{lsttr2}) are to be obtained by solving $m$ systems of linear equations. Since we are mainly interested in the bivariate transformation ${\bf F}:\RR^2 \rightarrow \RR^2$, we take $m=2$, and so from (\ref{lsttr}) and (\ref{lsttr2}) we have to consider respectively
\begin{equation}
F_k({\bf x})=F_k(x_1, x_2)=\sum_{j=1}^N c_{j}  \varphi_{2,h}(\vert\vert{\textbf{x}}-{\textbf{x}_j}\vert\vert_2), \nonumber
\end{equation}
\begin{equation}
F_k({\bf x})=F_k(x_1, x_2)=\sum_{j=1}^N c_{j}  \varphi_{1,h}(x_1-x_{j1}) \varphi_{1,h}(x_2-x_{j2}), \nonumber
\end{equation}
which, imposing the conditions $F_1({\bf x}_j)={t}_{j1}$ and $F_2({\bf x}_j)={t}_{j2}$, for $j=1,2,\ldots, N$, give two associated linear systems for each of them.

In general, when radial basis functions with compact support are used for landmark-based image registration, a significant reduction of the global influence  is observed. Moreover, even significant local deformations can be well registered by this approach. Therefore, their use is advantageous when it allows us to deal with local changes in images, as in medical images \cite{fornefett01}.


\section{Lobachevsky splines} \label{sf_pt}

In this section we consider a class of spline functions, called \textsl{Lobachevsky splines}, arising in probability theory.
They were introduced by Lobachevsky in 1842 (see, e.g., \cite{Renyi66}) and are identical (up to a scaling factor) to zero centered uniform B-splines (see below). The use of Lobachevsky splines for interpolation and integration of multivariate scattered data was considered in \cite{Allasia13a, Allasia13b, Allasia14} and for image registration in \cite{Allasia10,Allasia11b}.

The sequence of Lobachevsky splines $\{ f_n^*(x), n=1,2,\ldots\}$ converges for $n\to\infty$ to the normal density function with expectation 0 and standard deviation 1 (see \cite{Allasia13a}), i.e.
\begin{eqnarray}
\label{eq:a2}
\lim_{n\to\infty} f^*_n(x)=\frac{1}{\sqrt{2\pi}}\exp{\Big(\frac{-x^2}{2}\Big)},
\end{eqnarray}
and moreover the convergence is uniform for all $x\in \RR$.

To get an explicit expression of $f^*_n$ it is convenient referring to the function $f_n$, related to $f_n^*$ by the relation
\begin{equation}
f^*_n(x)= a\sqrt{\frac{n}{3}}~f_n\bigg(a\sqrt{\frac{n}{3}}~x\bigg). \nonumber
\end{equation}
Both $f_n^*$ and $f_n$ can be interpreted as probability density functions having supports on $[-\sqrt{3n},\sqrt{3n}]$ and $[-na,na]$, $a \in \RR^+$, respectively. Now, $f_n$ is given by the convolution product
\begin{equation}
	f_n(x)=\int_{-\infty}^{+\infty}f_1(u)f_{n-1}(x-u)du=\frac{1}{2a}\int_{-a}^{+a}f_{n-1}(x-u)du, \nonumber
\end{equation}
because by definition $f_1(x)=1/(2a)$ for $-a\le x\le +a$, and $f_1(x)=0$ elsewhere. Setting $x-u=t$, we get the recursive formula
\begin{equation}
\label{eq:a4}
f_n(x)=\frac{1}{2a}\int_{x-a}^{x+a}f_{n-1}(t)dt,\qquad n=2, 3, \ldots \nonumber
\end{equation} 

Thus, reasoning by induction (see \cite{Feller71}), we may express $f_n$, for $n=1,2,\ldots$, by the formula

\begin{eqnarray}
\label{eq:a5b}
f_n(x)={\displaystyle \frac{1}{(2a)^n(n-1)!}}\sum_{k=0}^{n}(-1)^k\Big(\displaystyle{n \atop k}\Big)[x+(n-2k)a]_+^{n-1},
\end{eqnarray}
where $(\cdot)_+$ is the truncated power function. Sometimes it may be convenient to consider a different form of (\ref{eq:a5b}), namely
\begin{eqnarray}
\label{eq:a5a}
f_n(x)={\displaystyle \frac{1}{(2a)^n(n-1)!}}\sum_{k=0}^{\big\lfloor \frac{na+x}{2a}\big\rfloor}(-1)^k
\Big(\displaystyle{n \atop k}\Big)[x+(n-2k)a]^{n-1}, 
\end{eqnarray}
for $-na\le x\le na$, and $f_n(x)=0$ otherwise, where $\lfloor\cdot\rfloor$ means the greatest integer less than or equal to the argument.

The piecewise function $f_n$ is graphically represented by arcs of parabolas of degree $n-1$; the first $n-2$ derivatives of different arcs of parabolas are equal at the knots (the points where the parabolic pieces are joined), i.e. $f_n\in C^{n-2}[-na, na]$. The knots are uniformly spaced with size $2a$, and $f_n$ is an even function with support $[-na, na]$.  In the case $a=1/2$ the knots are unit spaced and this class of splines enjoys very interesting applications in signal and image processing, wavelets theory, etc. (see, e.g., \cite{Brinks08, Cheney99, Unser91}). From a computational viewpoint it is convenient to evaluate $f_n$ starting from the pieces defined on $[-na, 0]$ and then obtain the pieces on $[0, na]$ by symmetry. Moreover, each piece on $[-na, 0]$ can be obtained by the preceding one by simply adding a term, as clearly appears from (\ref{eq:a5a}).

In order to show the connection between the Lobachevsky splines and the B-splines (see \cite{Allasia13a}), we may consider the following expression
$$
u_n(t)=2af_n\Big[2a\Big(t-\frac{n}{2}\Big)\Big]=\frac{1}{(n-1)!}\sum_{k=0}^n(-1)^k\Big(\displaystyle{n \atop k}\Big)(t-k)_+^{n-1},
$$
which is a well-known form of the classic B-splines on the support $[0, n]$ (see, e.g. \cite{schumaker93}). 
Lobachevsky splines satisfy a three-term recurrence relation for $n=2,3,\ldots$, namely,
\begin{equation}
f_n(x)=\frac{1}{n-1}\Big[\frac{na+x}{2a}f_{n-1}(x+a)+\frac{na-x}{2a}f_{n-1}(x-a)\Big],  \label{eq:a6c}
\end{equation} 
which may be interesting from a computational point of view. In fact, this relation is very stable, whereas (\ref{eq:a5b}) could suffer from loss-of-significance errors owing to subtraction of nearly equal quantities. Note that, for computational reasons, $f_1$ in (\ref{eq:a6c}) must be taken in the interval $[-a,a[$.

Considering the connection between $f_n(x)$ and $f^*_n(x)$, the limit (\ref{eq:a2}) becomes
\begin{equation}
\label{eq:a6}
\lim_{n\to \infty} f^*_n(x)=\lim_{n\to\infty} a\sqrt{\frac{n}{3}}~f_n\bigg(a\sqrt{\frac{n}{3}}~x\bigg)=\frac{1}{\sqrt{2\pi}}\exp{\Big(\frac{-x^2}{2}\Big)}. \nonumber
\end{equation}

Noteworthy convergence properties are satisfied also by integrals and derivatives of Lobachevsky splines. From the central limit theorem for the convergence in distribution (see, e.g., \cite{Renyi66}) we have 
\begin{equation}
\lim_{n\to\infty}\int_{-\infty}^x f_n^*(t)dt=\lim_{n\to\infty}\int_{-\infty}^x a\sqrt{\frac{n}{3}}~
f_n\bigg(a\sqrt{\frac{n}{3}}~t\bigg)dt=\int_{-\infty}^x \frac{1}{\sqrt{2\pi}}\exp{\Big(\frac{-t^2}{2}\Big)}dt. \nonumber
\end{equation}
This result is also a direct consequence of (\ref{eq:a2}) and Lebesgue's Dominate Convergence Theorem. 
 
The asymptotic behaviour of  derivatives of Lobachevsky splines is described by the following result \cite{Brinks08}
\begin{equation}
\lim_{n\to \infty}D^q f^*_n(x)=\lim_{n\to \infty}D^q \bigg[a\sqrt{\frac{n}{3}}~
f_n\bigg(a\sqrt{\frac{n}{3}}~x\bigg)\bigg]=D^q \bigg[\frac{1}{\sqrt{2\pi}}\exp{\Big(\frac{-x^2}{2}\Big)}\bigg],\nonumber
\label{eq:a6b}
\end{equation}
that is, $q$ being a fixed integer, the sequence $D^q f_n^*(x),~n=q+2, q+3, \ldots$, of the $q$-th derivatives of $f_n^*(x)$ converges to the $q$-th derivative of the standardized normal density function.

\begin{remark} It is well-known that using Gaussians it is convenient to introduce a shape parameter $\alpha$. 
This trick can be conveniently applied  in our case as well considering $f_n^*(\alpha x)$ instead of $f_n^*(x)$. 
Using the shape parameter as a factor has the effect that a decrease of the shape parameter produces  flat basis functions, while increasing $\alpha$ leads to more peaked (or localized) basis functions. The same result can be achieved considering $f_n(x)$ and acting on the parameter $a$. Numerical computations show that the approximation performance obtained by $f^*_n(x)$ or $f_n(x)$ are totally equivalent, if the values of the parameters $\alpha$ and $a$ are   conveniently chosen \cite{Allasia13a,Cavoretto10}.
\end{remark}

In order to formulate the Lobachevsky spline transformation in the context of image registration, we may consider the following definition (see \cite{Allasia10,Allasia11b}).

\begin{definition} \label{ddd}
Given a set of source landmark points 
${\cal S}_N=\{ {\bf x}_j\in \RR^m,j=1,2,\ldots,N\}$, with associated the corresponding set of target landmark points ${\cal T}_N=\{ {\bf t}_j\in \RR^m,j=1,2,\ldots,N\}$, a {\sl Lobachevsky  spline transformation} ${\bf L}_n:$ $ {\mathbb R}^m\rightarrow {\mathbb R}^m$ is such that each its component 
\begin{equation}
(L_n)_k:\RR^m\rightarrow \RR, \hskip0.5cm k=1,2,\ldots ,m, \nonumber
\end{equation}
takes the form
\begin{eqnarray}
	\label{lst}
	(L_n)_k({\bf x})= (L_n)_k(x_1,x_2,\ldots ,x_m)&=&\sum_{j=1}^N c_{j}  f_n^* (x_1-x_{j1})f_n^* (x_2-x_{j2})\ldots f_n^* (x_m-x_{jm})\\
	&=&\sum_{j=1}^N d_{j}  f_n (x_1-x_{j1})f_n (x_2-x_{j2})\cdots f_n (x_m-x_{jm}), \nonumber
\end{eqnarray}
where 
\begin{equation}
f^*_n(x_j-x_{ji})=a\sqrt{\frac{n}{3}}~f_n\bigg(a\sqrt{\frac{n}{3}}~(x_i-x_{ji})\bigg), \hspace{0.5cm} i=1,2,\ldots,m, \nonumber
\end{equation} 
$f_n$ is given in (\ref{eq:a5b}), $n$ is even, $(x_1,x_2,\ldots,x_m)$ is any point in $\RR^m$, and $(x_{j1},x_{j2},\ldots,x_{jm})$ is a data site.
\end{definition}

Note that $(L_n)_k(\textbf{x})$ in (\ref{lst}) is a linear combination of the product of univariate shifted strictly  positive define functions $f^*_n (x-x_j)$, which are spline functions with compact supports $[x_j-\sqrt{3n}, x_j+\sqrt{3n}]$ and continuous up to $n-2$ derivatives. Furthermore, we remark that $(L_n)_k(\textbf{x})$ might also be expressed as a linear combination (with different coefficients $d_{j}$) of the product of the shifted strictly positive define functions $f_n (x-x_j)$, whose compact supports are $[x_j-na, x_j+na]$. 

Referring to Definition \ref{ddd}, we observe that the transformation function $(L_n)_k:$ ${\mathbb R}^m\rightarrow {\mathbb R}$ has to be calculated for each $k=1,2,\ldots,m$, and the parameters $c_{j}$ (or $d_{j}$) in (\ref{lst}) are to be obtained by solving $m$ systems of linear equations. Since we are mainly interested in the bivariate transformation ${\bf L}_n:\RR^2 \rightarrow \RR^2$, we take $m=2$ in (\ref{lst}), requiring that ${\bf L}_n$ solves Problem \ref{lit}. Therefore, we have to consider
\begin{equation}
(L_n)_k({\bf x})=(L_n)_k(x_1, x_2)=\sum_{j=1}^N c_{j}  f_n^* (x_1-x_{j1})f_n^* (x_2-x_{j2}) = \sum_{j=1}^N d_{j}  f_n (x_1-x_{j1})f_n (x_2-x_{j2}), \nonumber
\end{equation}
which, imposing the conditions $(L_n)_1({\bf x}_j)={t}_{j1}$ and $(L_n)_2({\bf x}_j)={t}_{j2}$, for $j=1,2,\ldots, N$, gives two associated linear systems.

Numerical results in \cite{Allasia13a} have shown that Lobachevsky splines are accurate and stable in the image registration context. Moreover, being compactly supported, they give rise to sparse interpolation matrices (for suitable choices of the shape parameters).

\section{Numerical comparison}

In this section, we compare the performances of the above recalled methods when applied to give local image transformations. In particular, we consider the modified Shepard's formula with local approximants given by Gaussian and thin plate spline (Shep-G and Shep-TPS, respectively), Wendland's transformations defined by means of bivariate Wendland's functions $C^2$ and $C^4$ (W2-2D and W4-2D, respectively) and by means of products of univariate Wendland's functions (W2-1D$\times$1D and W4-1D$\times$1D), and Lobachevsky spline transformations with $n = 4, 6$ (denoted by L4 and L6, respectively). Moreover, to get a more complete picture, we show also the registration results obtained by Gaussian (G) and thin plate spline (TPS), which are widely used in this context.

In order to test the different local interpolation schemes, we obtain several numerical results on some test cases. Here, for brevity, we refer to two test examples given in \cite{fornefett01, kohlrausch05}. They simulate typical medical cases where image portions scale or shift. These image portions represent rigid objects embedded in elastic material changing their position or form. The approach we propose can cope with local differences between corresponding images. In general these differences may be caused by the physical deformation of human tissue due to surgeries or pathological processes such as tumor growth or tumor resection. However, the aim is here to determine a transformation function, which connects the points of the source and target images, so that the target image is affected by the slightest possible deformation. 

Finally, we briefly present some experimental results obtained by applying Gaussian, thin plate spline, Lobachev\-sky splines, bivariate Wendland's functions and products of univariate Wendland's functions to real image data. More precisely, we consider two X-ray images of the cervical of an anonymous patient taken at different times. Note that this example has simply an illustrative purpose.


\subsection{Test example 1: square shift and scaling}

\begin{figure}[ht!]
\begin{center}
\includegraphics[width=8.cm]{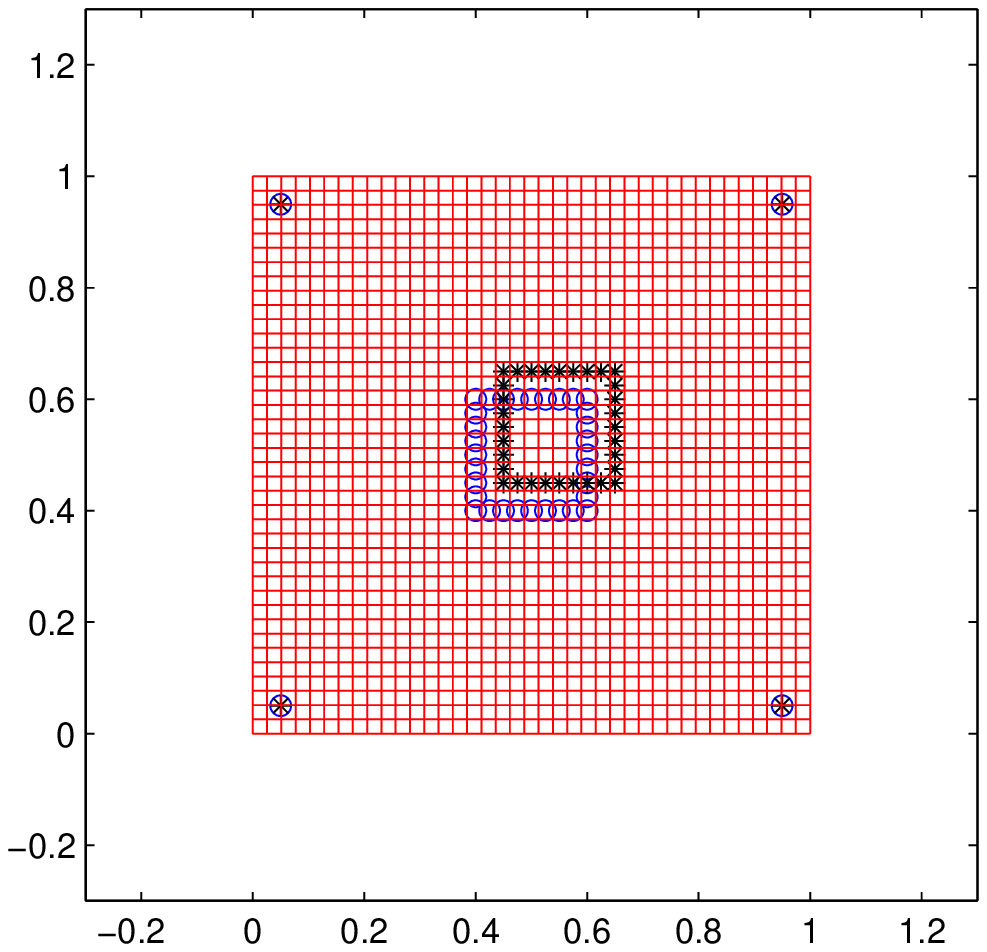}
\includegraphics[width=8.cm]{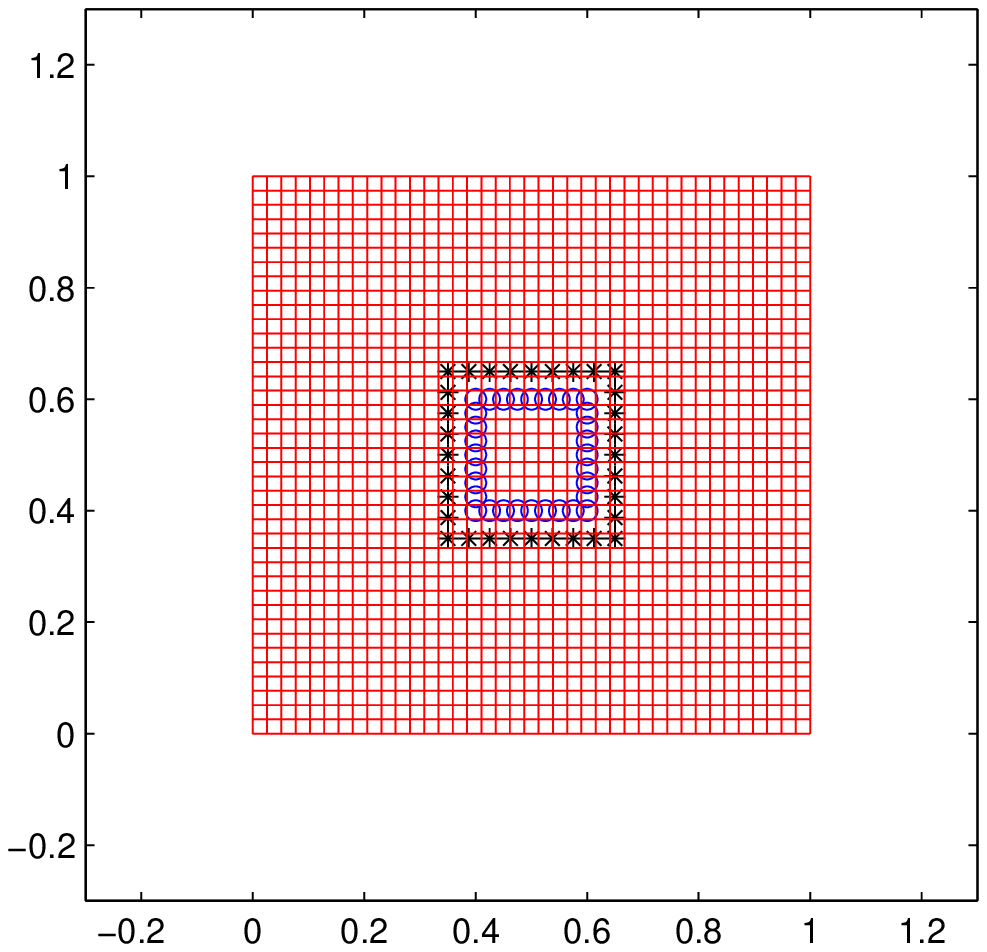}
\end{center}  
  \caption{Case 1, square shift (left); Case 2, square scaling (right); source and target landmarks.}
\label{MI_1}
\end{figure}

\begin{figure}[ht!]
\begin{center}
\includegraphics[width=8.cm]{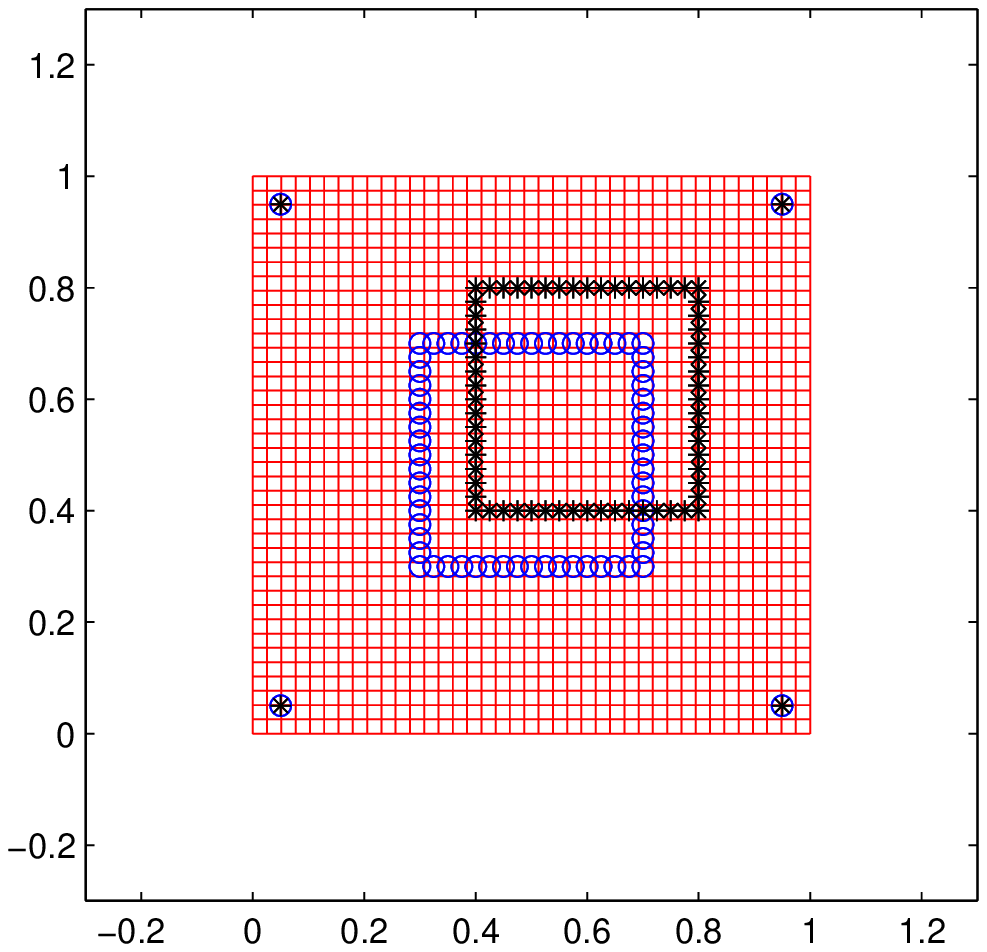}
\includegraphics[width=8.cm]{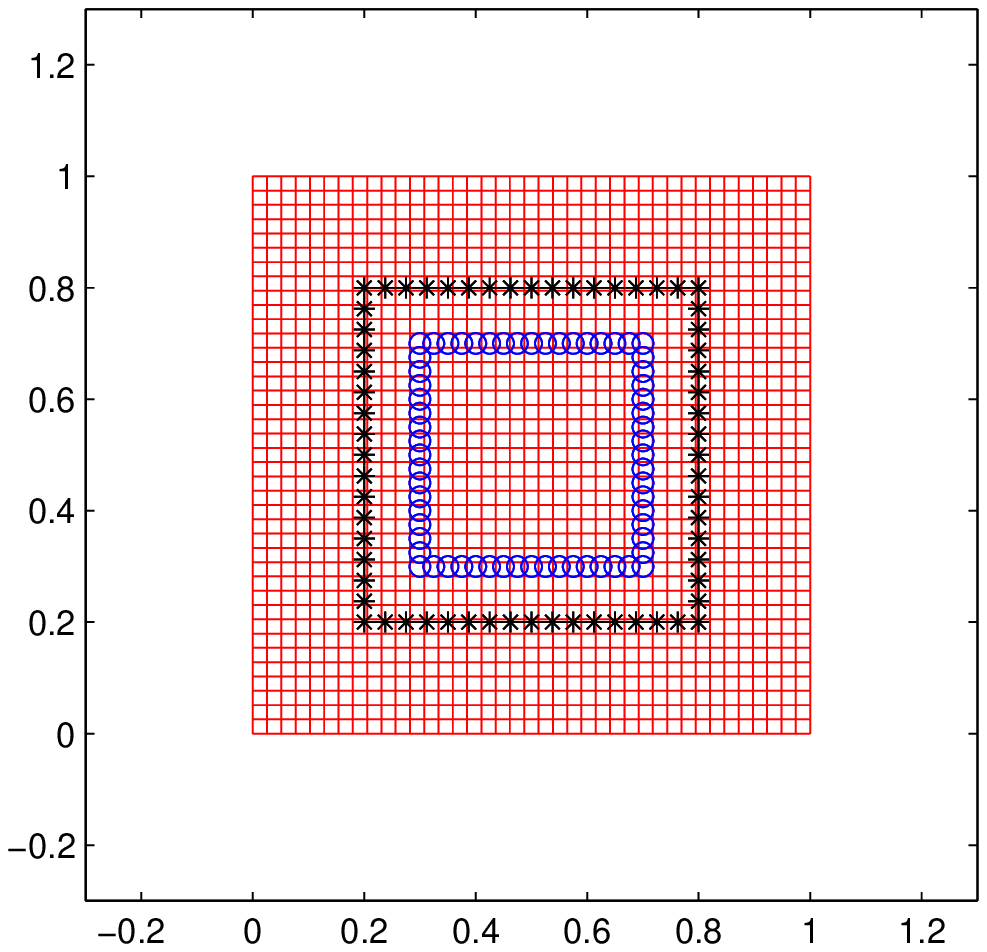}
\end{center}  
  \caption{Case 3, square shift (left); Case 4, square scaling (right); source and target landmarks.}
\label{new_MI_STI_large}
\end{figure}

We denote by ${\cal X}$ the set of $40 \times 40$ corner points of a regular grid superimposed on the source images (see Figures \ref{MI_1} and \ref{new_MI_STI_large}). The grid is transformed using 32 (Cases 1 and 2: square shift and scaling, respectively) and 64 (Cases 3 and 4: square shift and scaling, respectively) landmarks and, in the case of square shift, also 4 quasi--landmarks, i.e. landmarks which have the same positions in both source and target images, to prevent an overall shift. 
The choice to consider two different pair landmark sets is justified by the necessity of testing the interpolation schemes on different situations: the former with a deformation concentrated in a small region of the image, the latter with a larger deformation, which probably influences the entire image transformation. The source and target image landmarks, shown in the images in Figure \ref{MI_1} and in Figure \ref{new_MI_STI_large} for both cases, are marked by a circle ($\circ$) and a star ($\star$), respectively. 

At first, we analyze the behaviour of the root mean squares error (RMSE) by varying the shape parameter in each interpolation scheme in order to determine the \lq\lq optimal\rq\rq~value of the parameter, namely the value which gives us the smallest error.  Root mean squares error is found computing the distances between the displacements of grid points ${\bf x}\in {\cal X}$ and the values $\textbf{F}(\textbf{x})$ obtained by the transformations. It assumes the following form
 
\begin{equation}
   \text{RMSE} = \sqrt{\frac{\sum_{{\bf x}\in {\cal X}} \left\|{\bf x} - {\bf F}({\bf x}) \right\|^2_2}{\sum_{{\bf x}\in {\cal X}} 1}},  \nonumber
\end{equation}
where $\left\| \cdot \right\|_2$ is the Euclidean norm.

It is well-known that the choice of parameter values for a specific method is a very difficult task, because the \lq\lq optimal\rq\rq\ values depend, at least to a certain extent, on the actual application. Hence, in order to obtain some general information on the parameter variation, the most suitable way is given by considering geometric models (squares and circles) of underlying objects. Moreover, for a fixed model the parameter values depend on the set of landmarks, i.e. their number and collocation.

As for the use of the root mean squares error, we remark that RMSE is currently considered in the literature as a suitable tool for measuring the performance of an image transformation method (see, e.g., \cite{Allasia11b,kohlrausch05}). Other performance indices are theoretically possible, but the RMSE appears the most practical and effective one.

Quantitative graphics for the accuracy of the registration results are shown in Figures 3--6. Thus, exploiting RMSE graphs, we can pick optimal values for the shape parameters and compare the registration results, obtained by applying the different interpolation methods to the test Cases 1-4. Their comparison points out goodness and effectiveness of the various approaches from the point of view of errors. In particular, such graphs are plotted for equispaced values of $\alpha \in [0.2, 2.0]$ for G, Shep-G, L4 and L6 (Figures 3--6, left), and equispaced values of $c\in [0.1,1]$ for W2-2D, W4-2D, W2-1D$\times$1D and W4-1D$\times$1D (Figures 3--6, right). Here the choice of the localization parameters for Shep-G and Shep-TPS was $N_L=25$ and $N_W=25$.

\begin{figure}[ht!]
\begin{center}
  \includegraphics[width=8.cm]{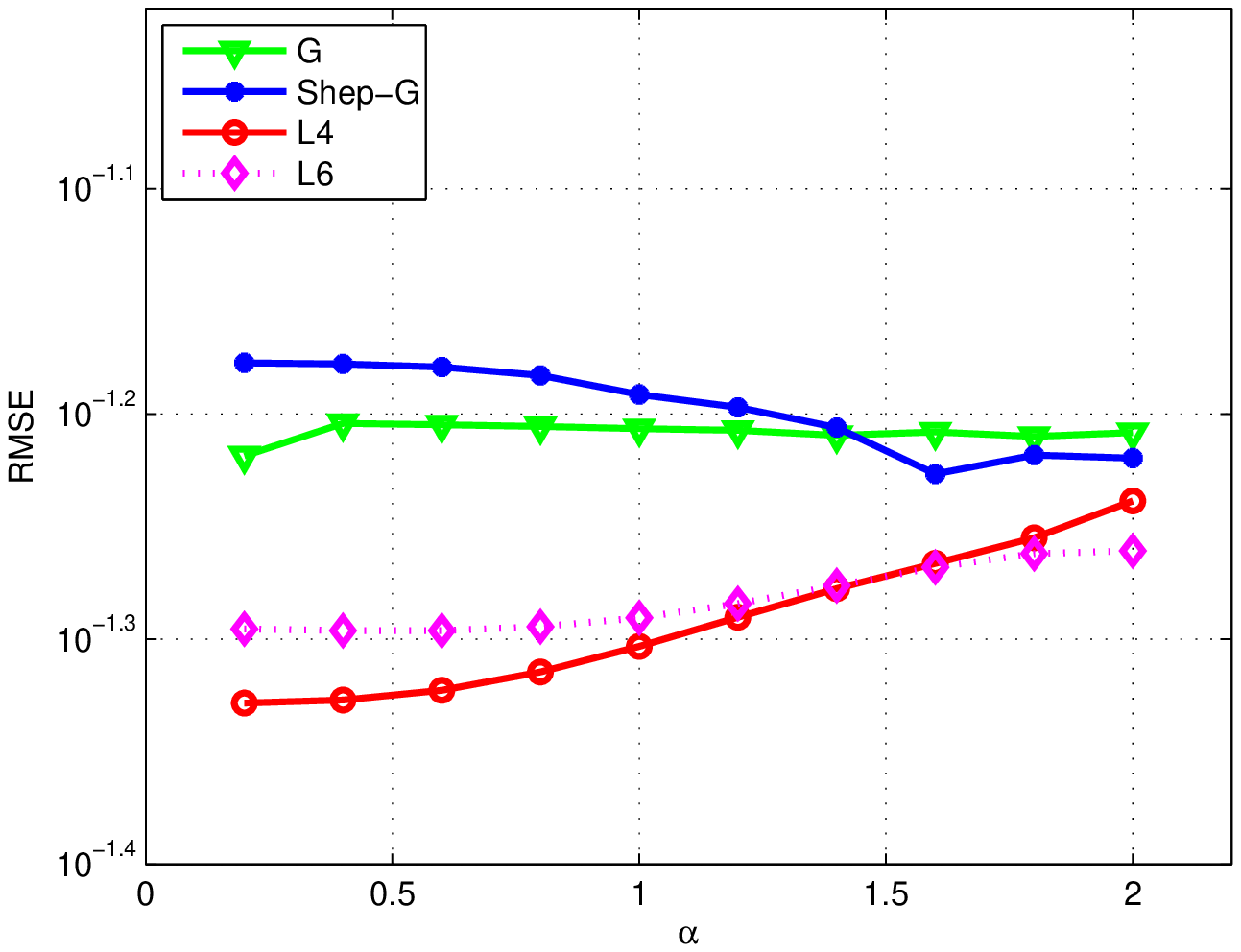}
  \includegraphics[width=8.cm]{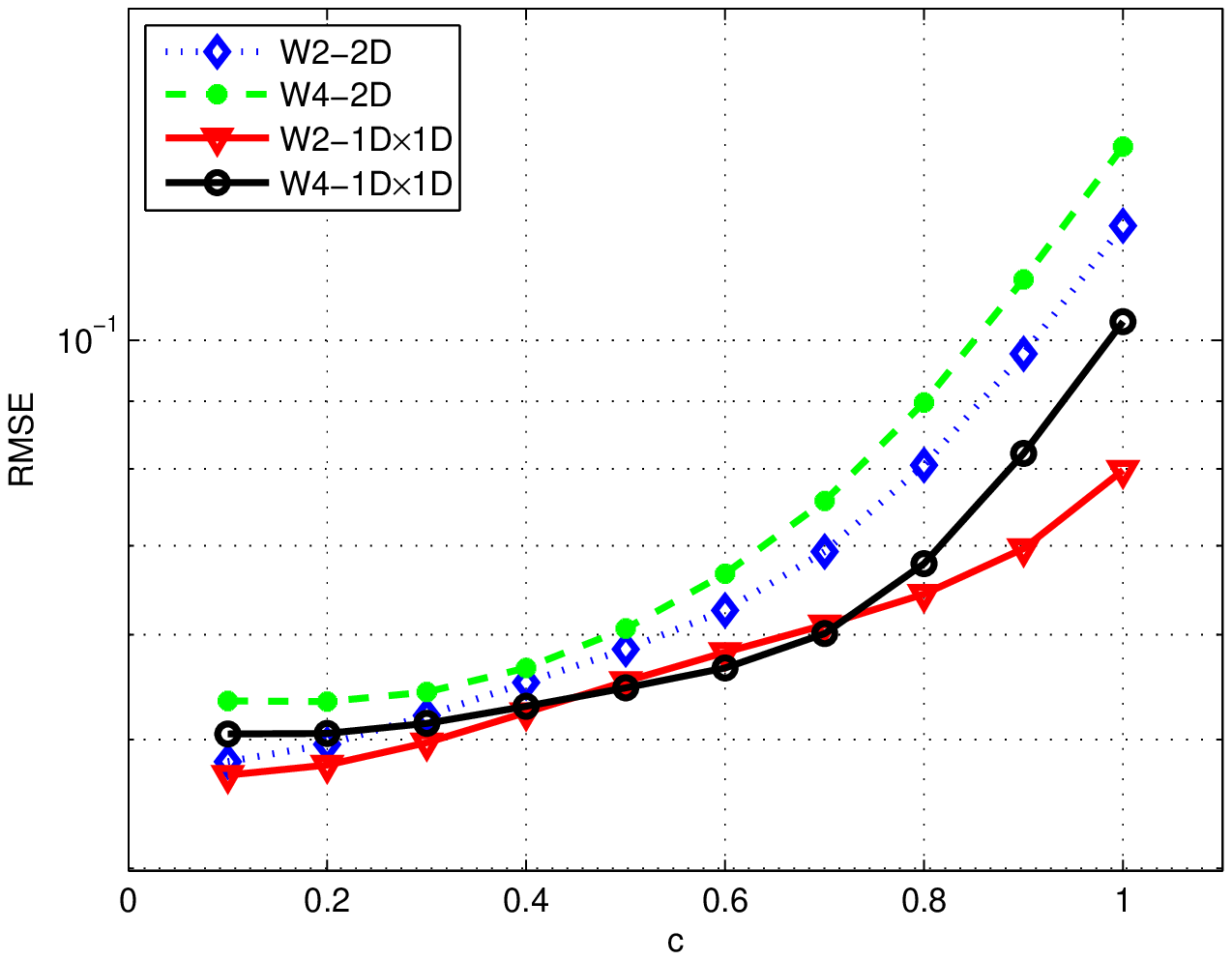}
\end{center}  
  \caption{Case 1, square shift, RMSEs by varying the shape parameter $\alpha \in [0.2,2.0]$ (left) and $c\in[0.1,1.0]$ (right).}
\label{RS1}
\end{figure}

\begin{figure}[ht!]
\begin{center}
  \includegraphics[width=8.cm]{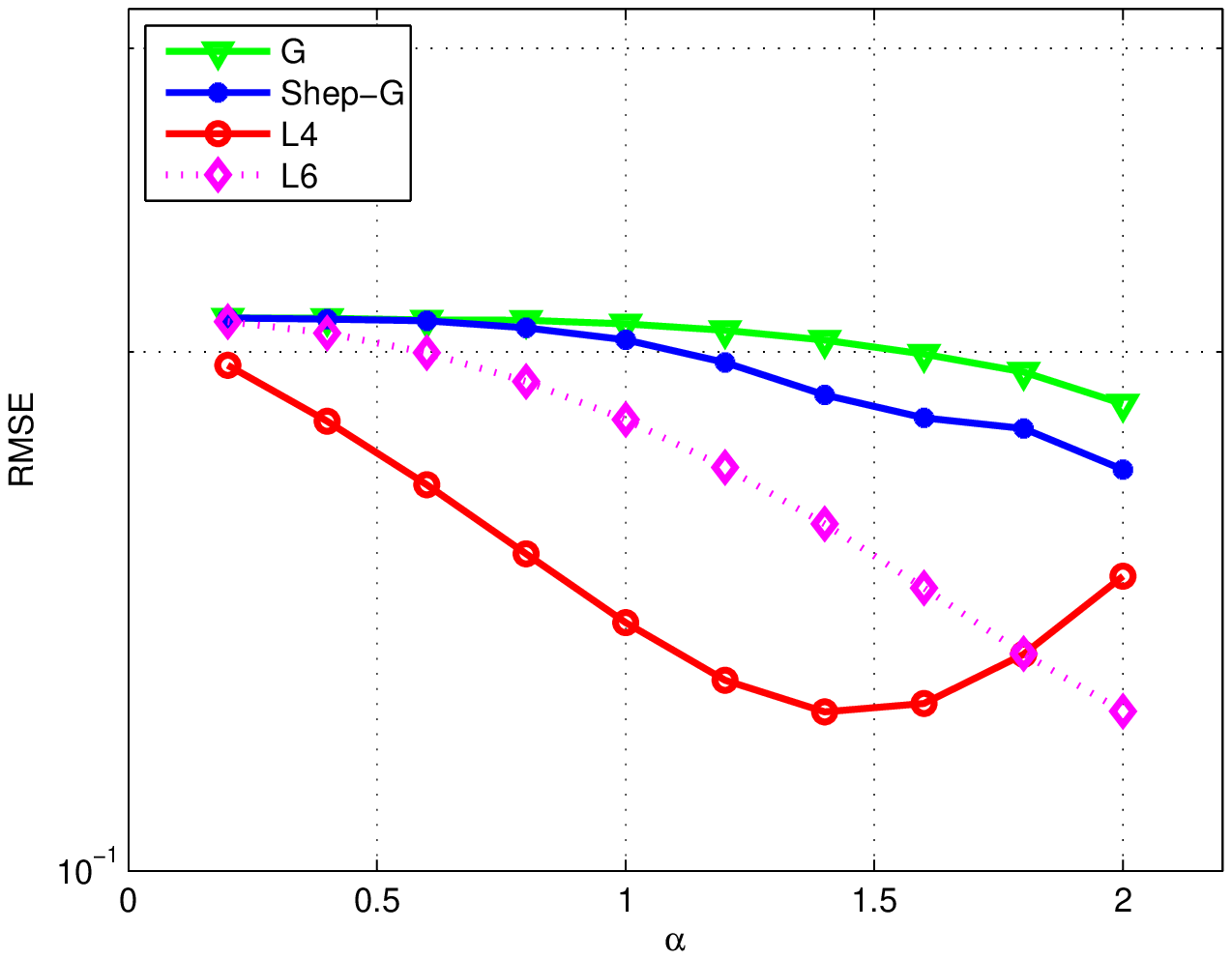}
  \includegraphics[width=8.cm]{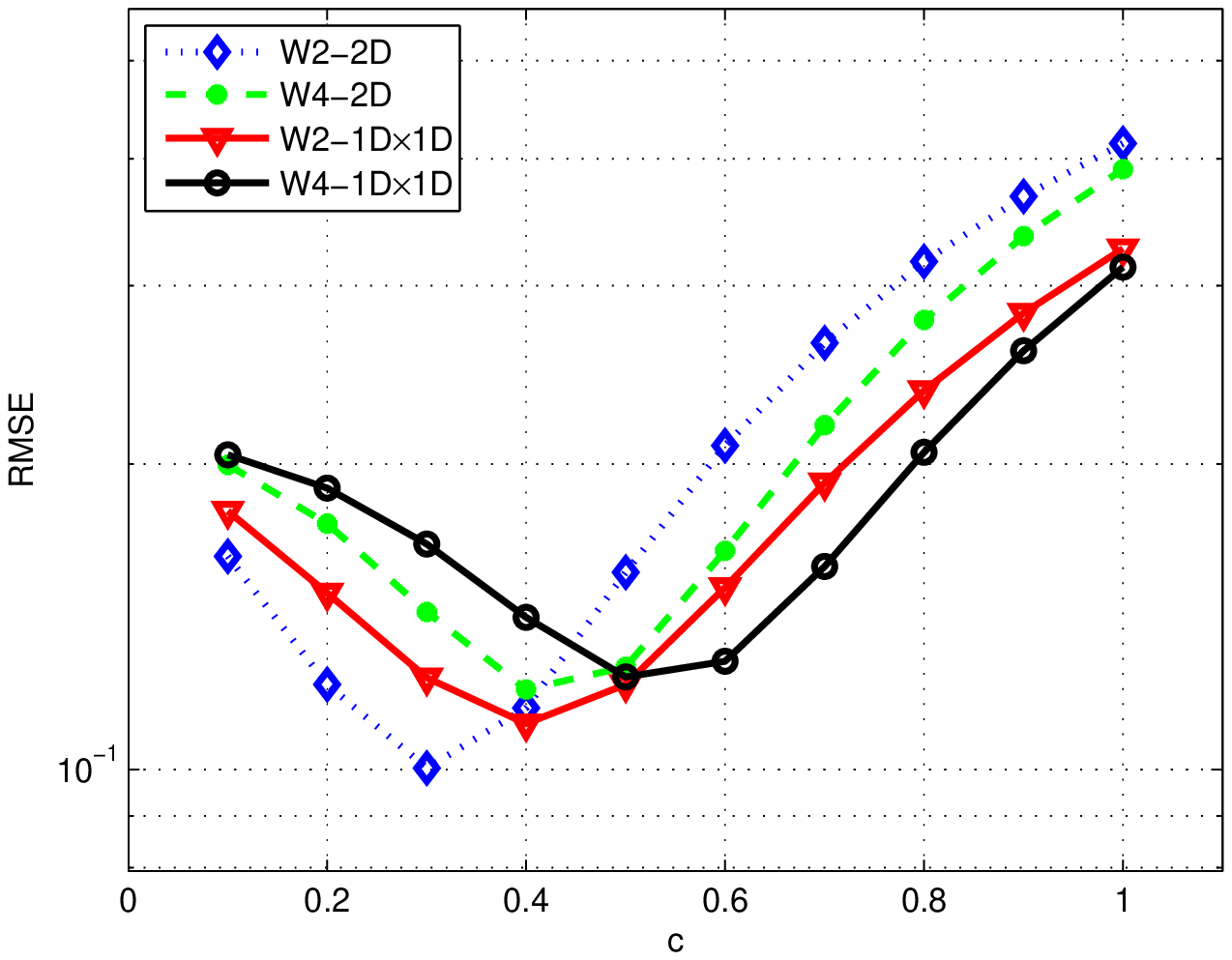}
\end{center}  
  \caption{Case 2, square scaling, RMSEs by varying the shape parameter $\alpha \in [0.2,2.0]$ (left) and $c\in[0.1,1.0]$ (right).}
\label{RS2}
\end{figure}

\begin{figure}[ht!]
\begin{center}
  \includegraphics[width=8.cm]{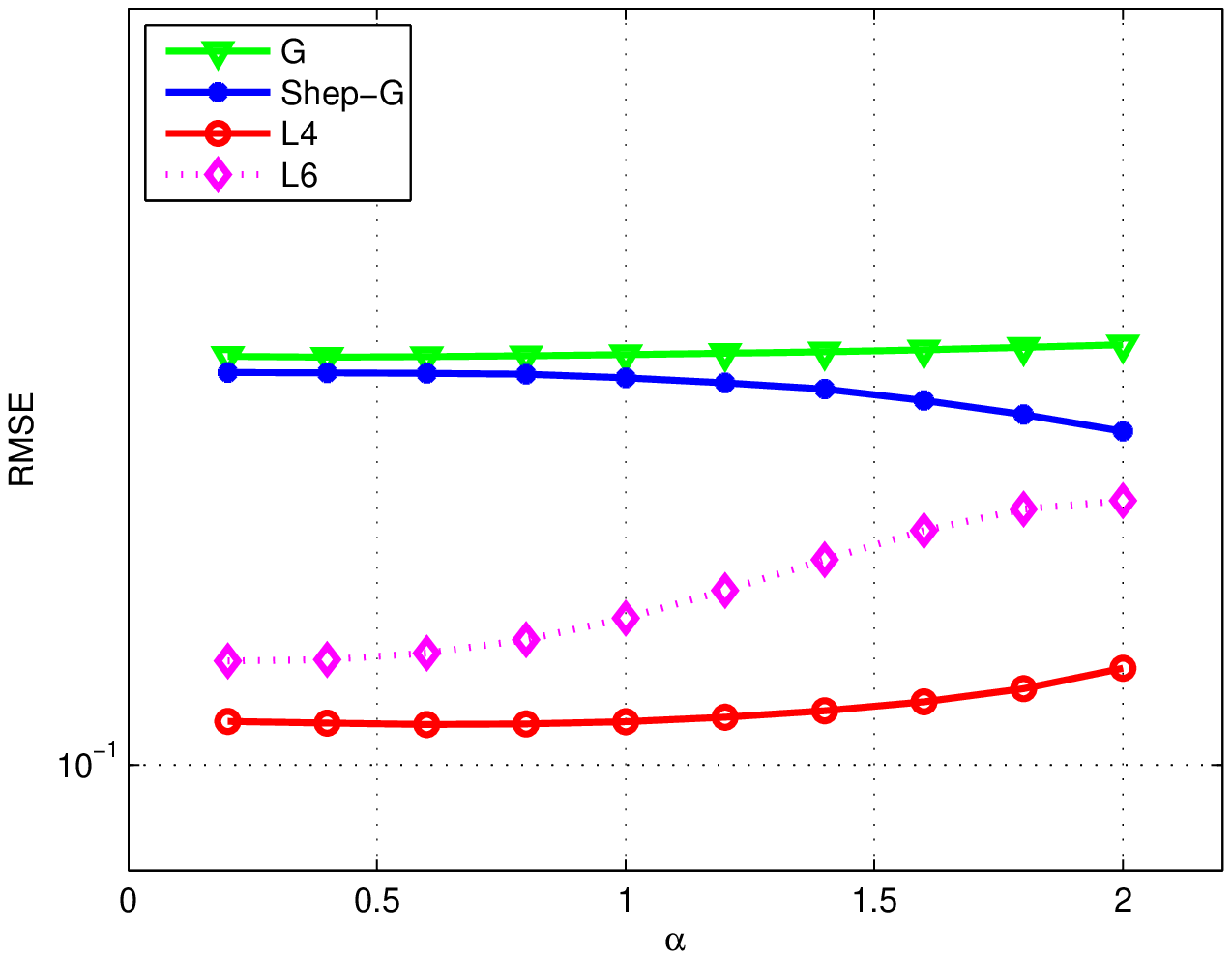}
  \includegraphics[width=8.cm]{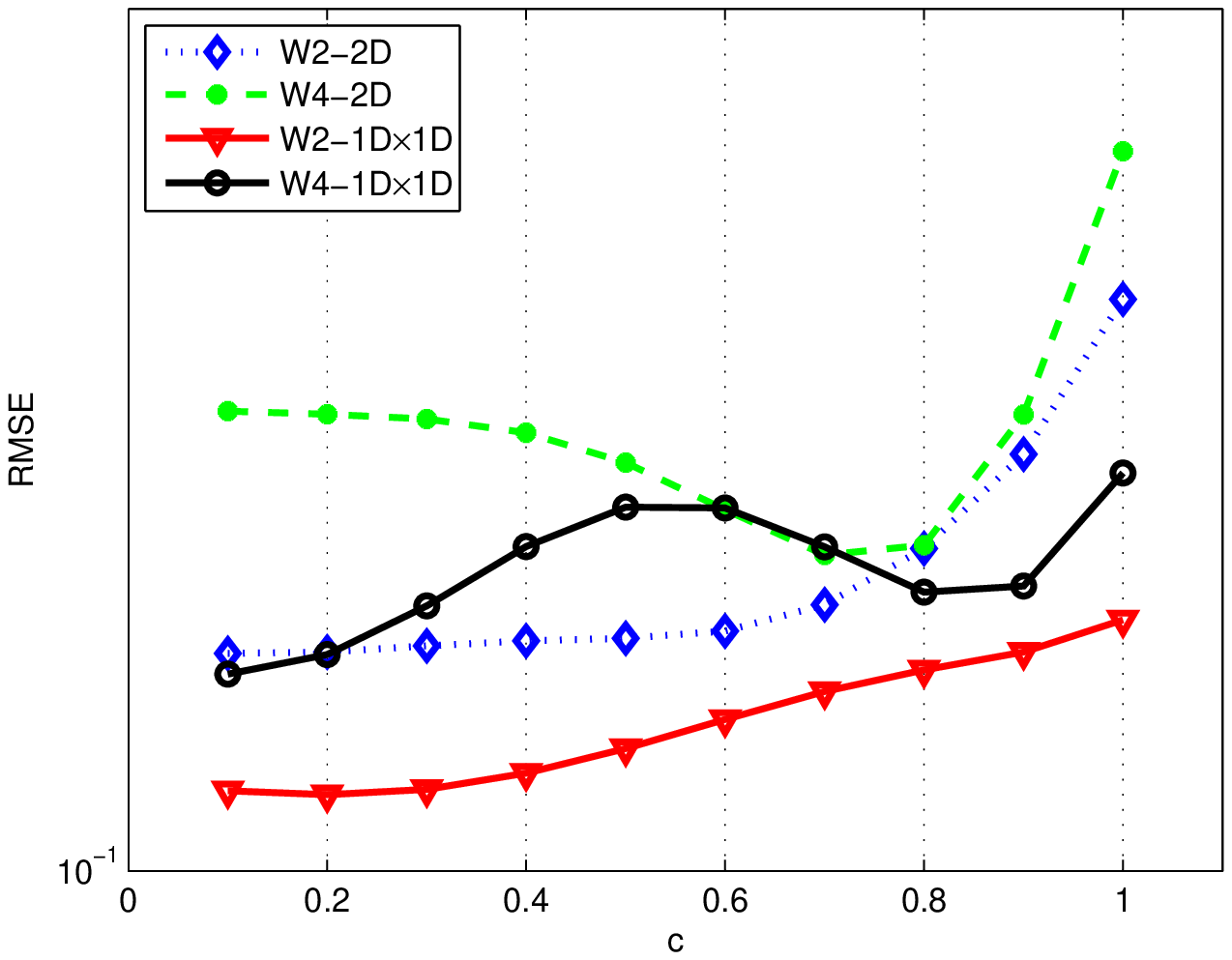}
\end{center}  
  \caption{Case 3, square shift, RMSEs by varying the shape parameter $\alpha \in [0.2,2.0]$ (left) and $c\in[0.1,1.0]$ (right).}
\label{RS3}
\end{figure}

\begin{figure}[ht!]
\begin{center}
  \includegraphics[width=8.cm]{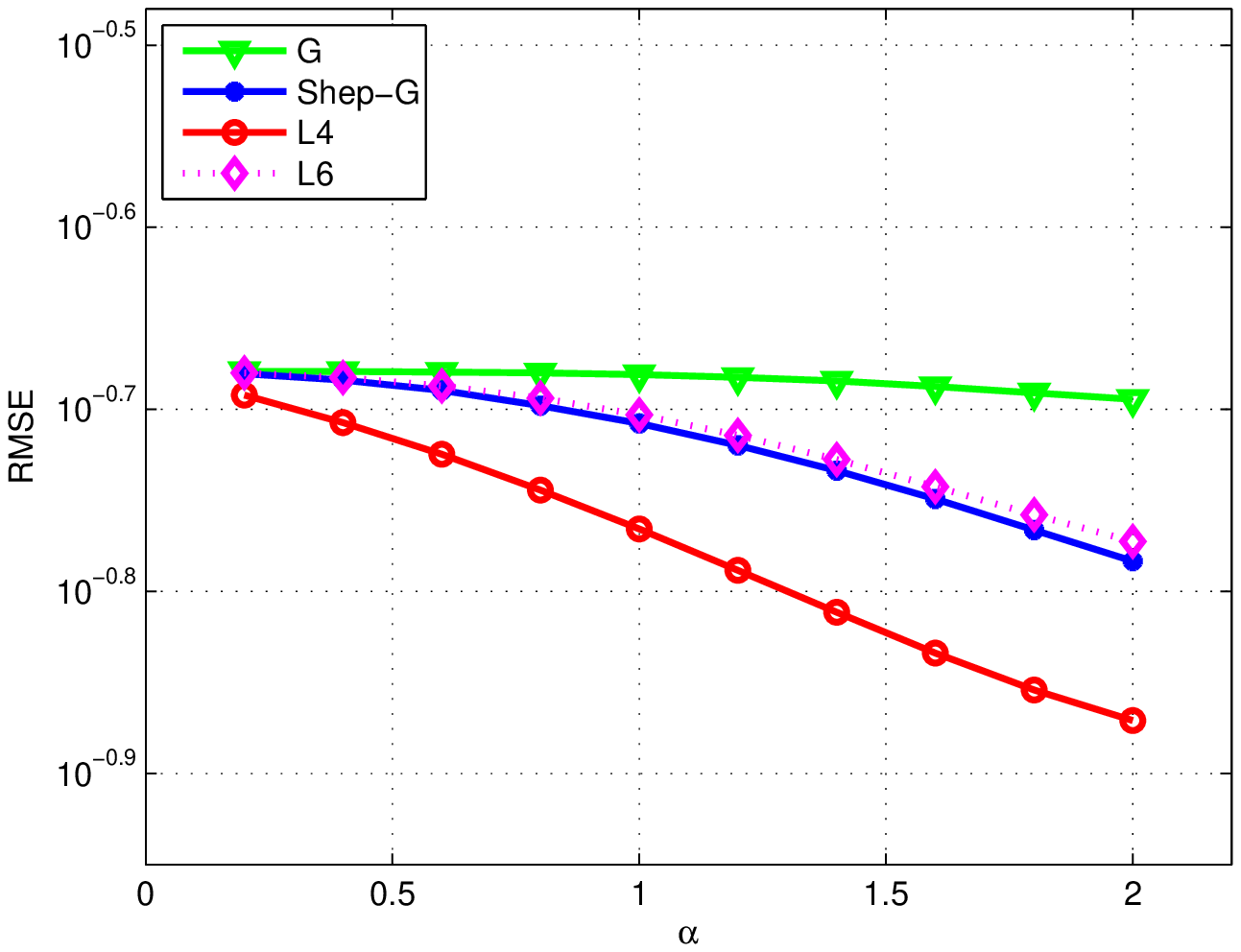}
  \includegraphics[width=8.cm]{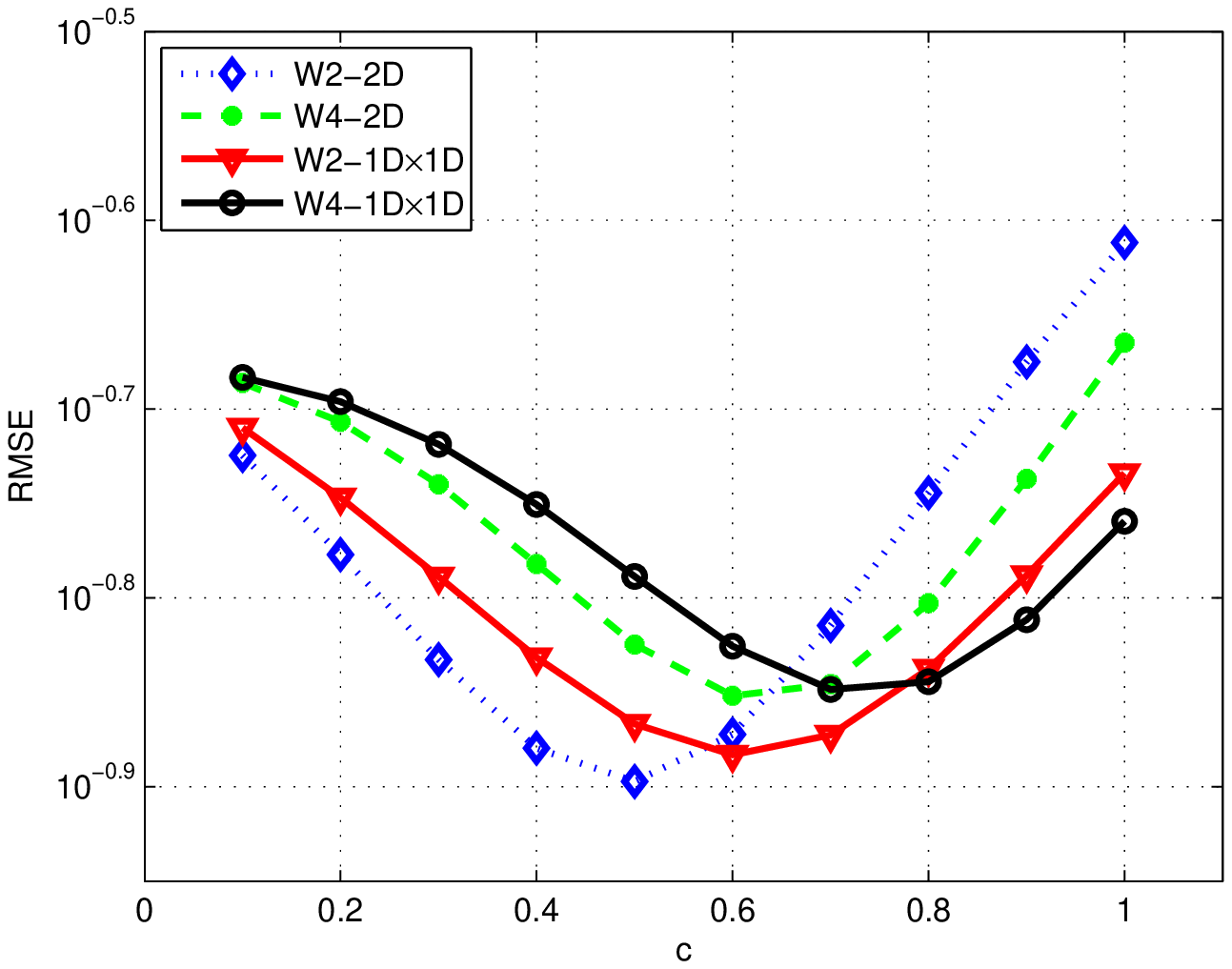}
\end{center}  
  \caption{Case 4, square scaling, RMSEs by varying the shape parameter $\alpha \in [0.2,2.0]$ (left) and $c\in[0.1,1.0]$ (right).}
\label{RS4}
\end{figure}

In Table \ref{C1_opt} we report the root mean squares errors (RMSEs) obtained by using the optimal values determined above. An analysis of errors is difficult, since all the considered interpolation schemes are known to be accurate and, in fact, all errors are of the same orders of magnitude. However, a comparison of RMSEs points out that in Cases 1 and 3 (when the square shifts using 32 and 64 landmarks, respectively) the results are better when thin plate spline, Wendland's function W2 (2D and 1D$\times$1D) and Lobachevsky spline L4 are used. In Case 2 (square scaling with 32 landmarks), errors are less when compactly supported Wendland's functions W2 and W4-2D are used, while in Case 4 (square scaling with 64 landmarks) better accuracy is obtained using W2 and Lobachevsky spline L4.  

Moreover, we note that, in the comparison of global and local methods, thin plate spline gives better results when the object shifts, whereas the local Shepard's method with Gaussian is preferable when the object scales, as already shown in \cite{Cavoretto-DeRossi08a,Cavoretto-DeRossi08b}.

\begin{Large}
\begin{table}[ht!]
\begin{center}
\begin{tabular}{|c|c|c|c|c|}
\hline
\rule[0mm]{0mm}{3ex}
  &  Case 1  & Case 2  &  Case 3 & Case 4   \\
\hline
\hline
\rule[0mm]{0mm}{3ex}
G               & $6.0461\text{E}-2$ & $1.8654\text{E}-1$ & $1.3639\text{E}-1$ & $2.0206\text{E}-1$  \\
$\alpha$       & $0.2$ & $2.0$ & $0.4$ & $2.0$  \\
\hline 
\rule[0mm]{0mm}{3ex}
TPS             & $4.3460\text{E}-2$ & $2.0929\text{E}-1$ & $1.0310\text{E}-1$ & $2.0929\text{E}-1$  \\
\hline
\rule[0mm]{0mm}{3ex}
Shep-G          & $5.9351\text{E}-2$ & $1.7087\text{E}-1$ & $1.2891\text{E}-1$ & $1.6464\text{E}-1$  \\
$\alpha$       & $1.6$ & $2.0$ & $2.0$ & $2.0$  \\
\hline 
\rule[0mm]{0mm}{3ex}
Shep-TPS        & $6.4435\text{E}-2$ & $2.0929\text{E}-1$ & $1.3275\text{E}-1$ & $2.0929\text{E}-1$  \\
\hline
\rule[0mm]{0mm}{3ex}
W2-2D           & $4.8120\text{E}-2$ & $1.0033\text{E}-1$ & $1.0911\text{E}-1$ & $1.2671\text{E}-1$  \\
$c$            & $0.1$ & $0.3$ & $0.1$ & $0.5$  \\
\hline 
\rule[0mm]{0mm}{3ex}
W4-2D           & $5.3417\text{E}-2$ & $1.1990\text{E}-1$ & $1.1349\text{E}-1$ & $1.4067\text{E}-1$  \\
$c$            & $0.2$ & $0.4$ & $0.7$ & $0.6$  \\
\hline 
\rule[0mm]{0mm}{3ex}
W2-1D$\times$1D & $4.7013\text{E}-2$ & $1.1098\text{E}-1$ & $1.0310\text{E}-1$ & $1.3089\text{E}-1$  \\
$c$            & $0.1$ & $0.4$ & $0.2$ & $0.6$  \\
\hline 
\rule[0mm]{0mm}{3ex}
W4-1D$\times$1D & $5.0482\text{E}-2$ & $1.2341\text{E}-1$ & $1.0820\text{E}-1$ & $1.4178\text{E}-1$  \\
$c$            & $0.1$ & $0.5$ & $0.1$ & $0.7$  \\
\hline 
\rule[0mm]{0mm}{3ex}
L4              & $4.6950\text{E}-2$ & $1.2368\text{E}-1$ & $1.0314\text{E}-1$ & $1.3462\text{E}-1$  \\
$\alpha$       & $0.2$ & $1.4$ & $0.6$ & $2.0$  \\
\hline 
\rule[0mm]{0mm}{3ex}
L6              & $5.0566\text{E}-2$ & $1.2374\text{E}-1$ & $1.0825\text{E}-1$ & $1.6880\text{E}-1$  \\
$\alpha$       & $0.4$ & $2.0$ & $0.2$ & $2.0$  \\
\hline
\end{tabular}
\end{center}
\caption{RMSEs for optimal values of $\alpha$ and $c$.}
\label{C1_opt}
\end{table}
\end{Large}

It is well known that also shape and smoothness of the transformed grids may be relevant to compare different transformations. To this aim in the following figures we present registration results obtained by the interpolation methods, using the previously determined optimal values for the parameters.  Here we observe significant differences. Registration results involving interpolation schemes W4-2D and W4-1D$\times$1D are not represented, since the results are worse than those obtained with the other schemes.

In Figure 7 we show registration results obtained by applying the interpolation schemes in the Case 1 of square shift when the number of landmarks is 32. Irregular grids are obtained using G and Shep-G, while smoother results are given by TPS, Shep-TPS, Lobachevsky splines and Wendland's functions. Transformed grids visibly smoother are obtained when TPS, Shep-TPS and W2-2D transformations are employed for the shift of a square using 64 landmarks (see Figure 9). Figure 8 and 10 concern the square scaling with 32 and 64 landmarks, respectively. In both cases, TPS and Shep-TPS transformations maintain the square form of the image, but the image is visibly enlarged, while the deformation is limited when Lobachevsky splines and Wendland's functions are used. Moreover, we point out that the condition numbers of interpolation matrices generated by Gaussian are very large in both cases, namely $10^{17}\div 10^{19}$.

\begin{figure}[ht!]
\begin{center}
\begin{minipage}{60mm}
\includegraphics[width=6.cm]{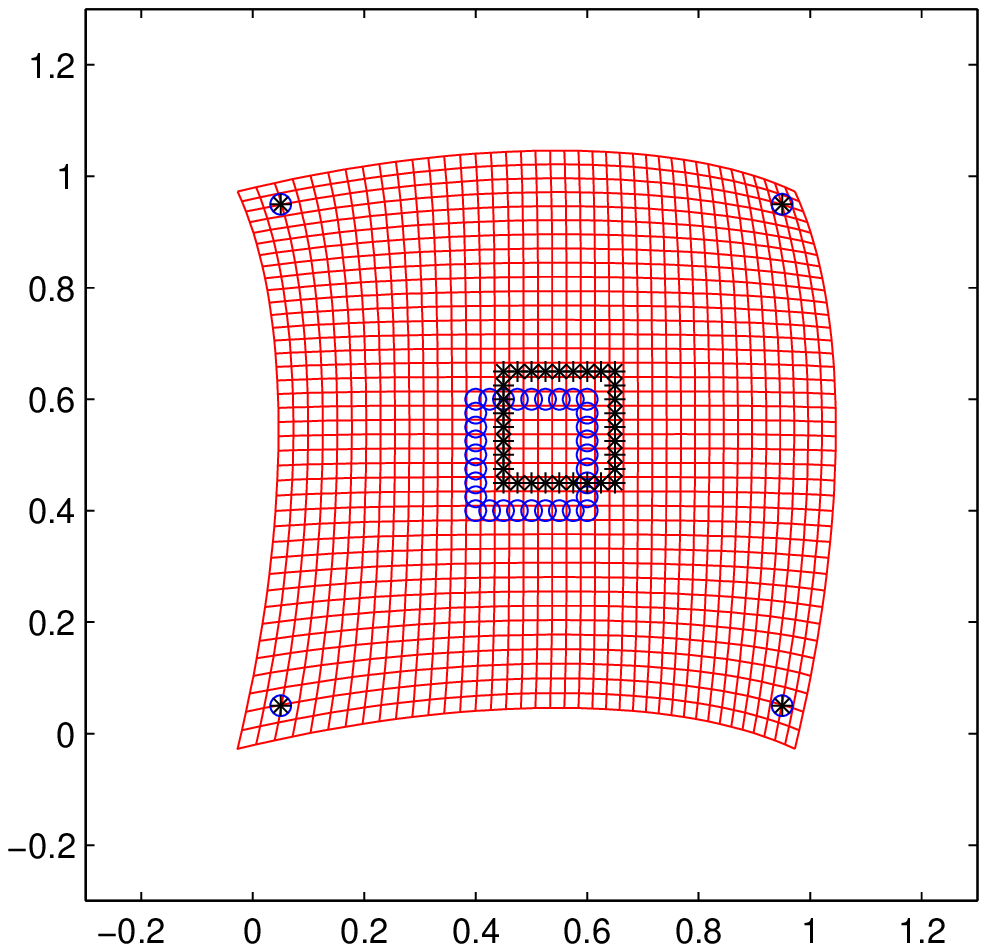}
\centerline{(a) G, $\alpha=0.2$}
\end{minipage}
\begin{minipage}{60mm}
\includegraphics[width=6.cm]{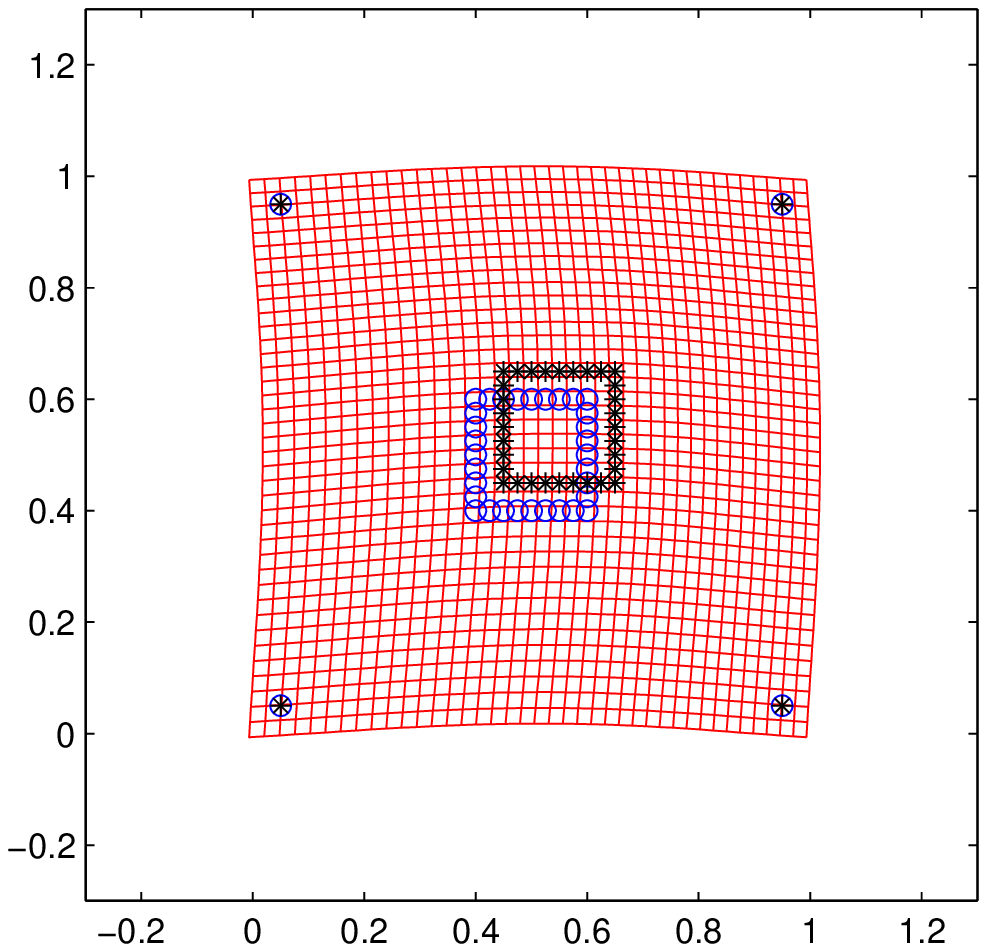}
\centerline{(b) TPS}
\end{minipage}\\
\begin{minipage}{60mm}
\includegraphics[width=6.cm]{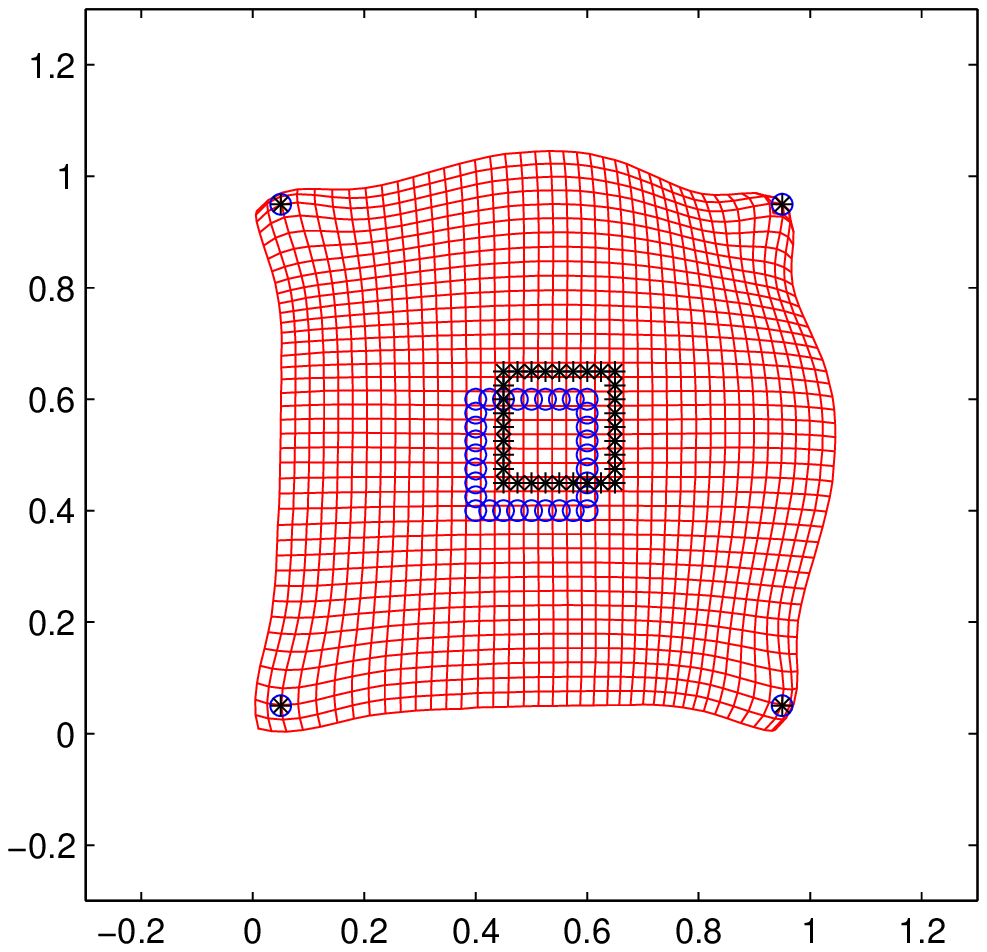}
\centerline{(c) Shep-G, $\alpha=1.6$}
\end{minipage}
\begin{minipage}{60mm}
\includegraphics[width=6.cm]{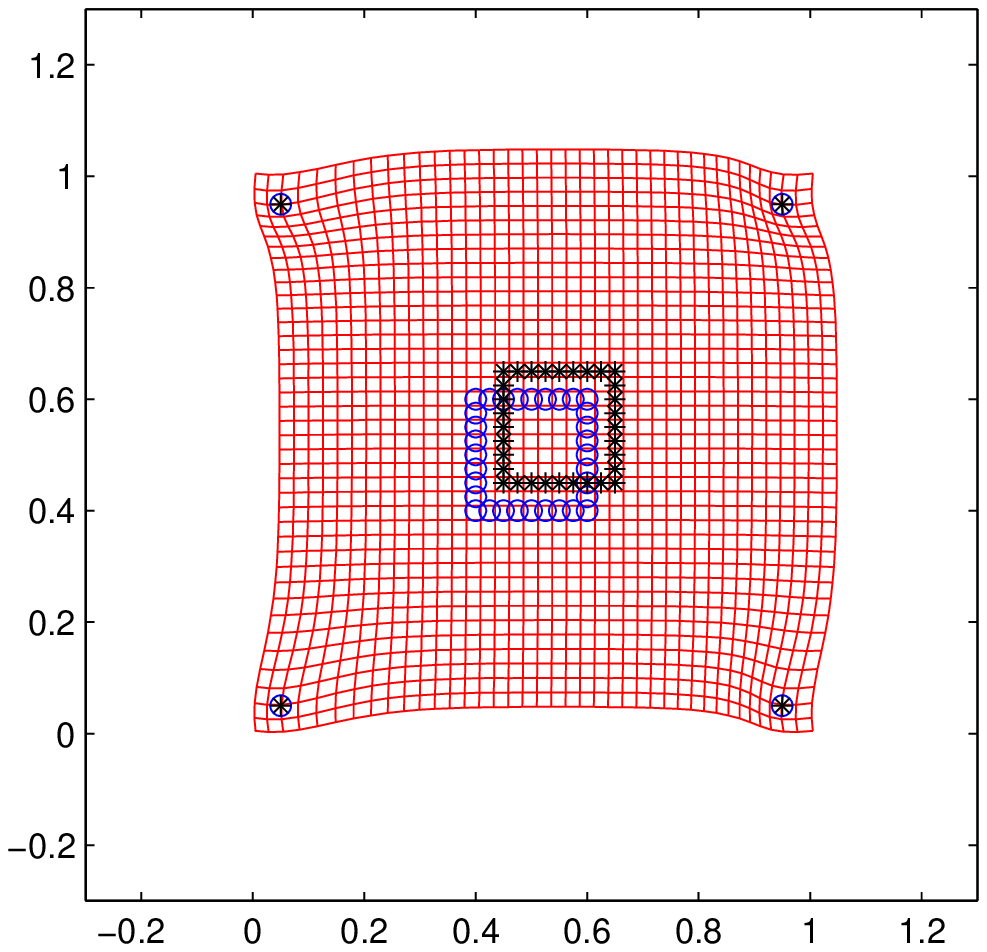}
\centerline{(d) Shep-TPS}
\end{minipage}\\
\begin{minipage}{60mm}
\includegraphics[width=6.cm]{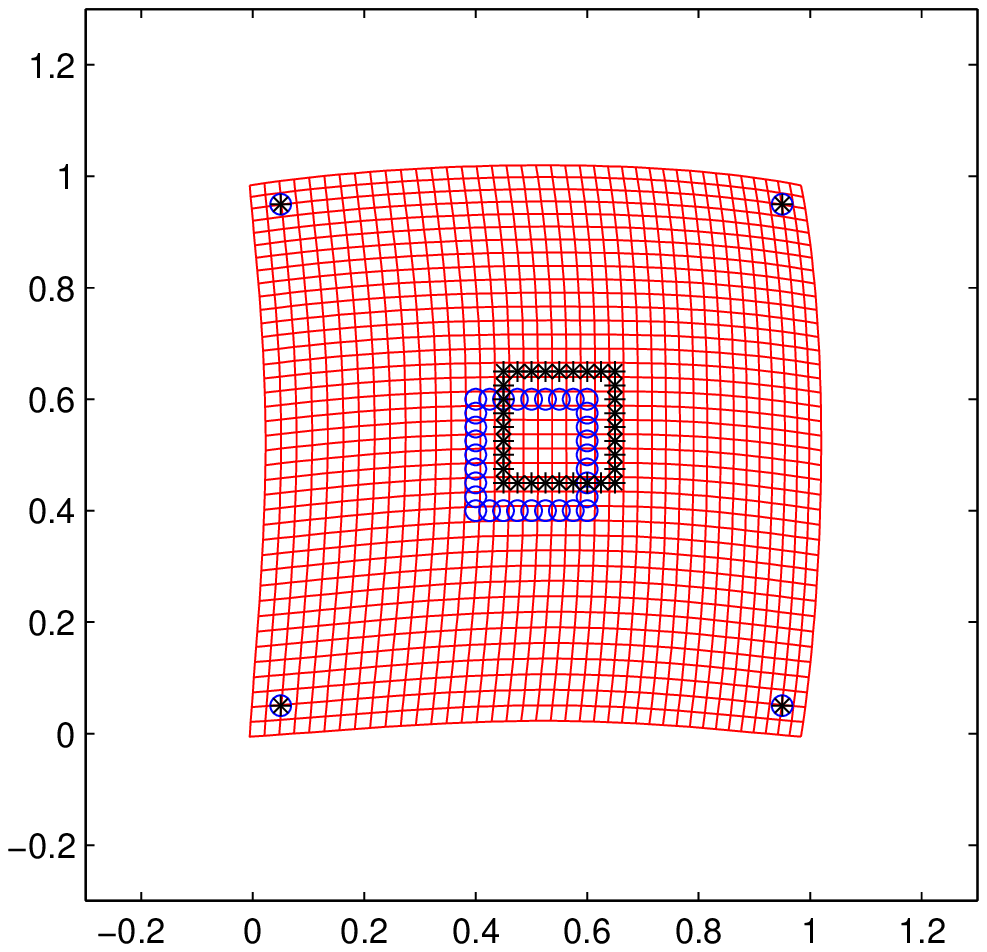}
\centerline{(e) W2-2D, $c=0.1$}
\end{minipage}
\begin{minipage}{60mm}
\includegraphics[width=6.cm]{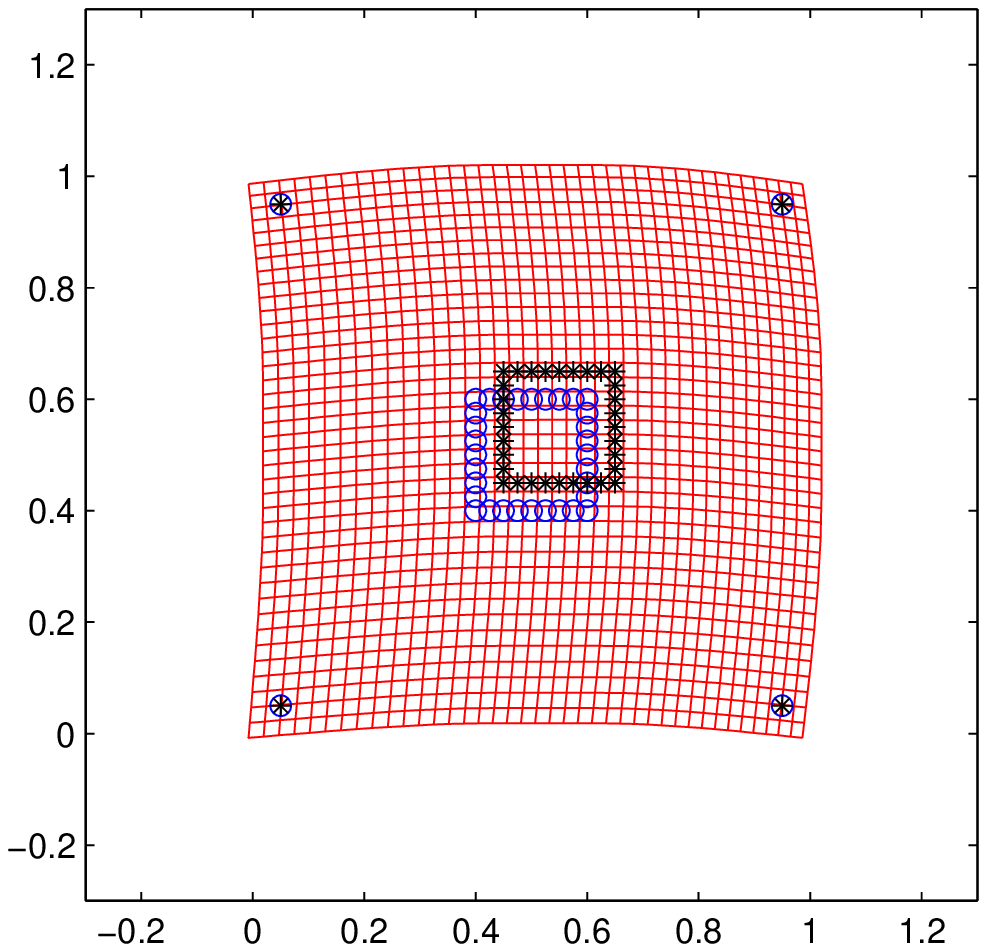}
\centerline{(f) W2-1D$\times$1D, $c=0.1$}
\end{minipage}
\begin{minipage}{60mm}
\includegraphics[width=6.cm]{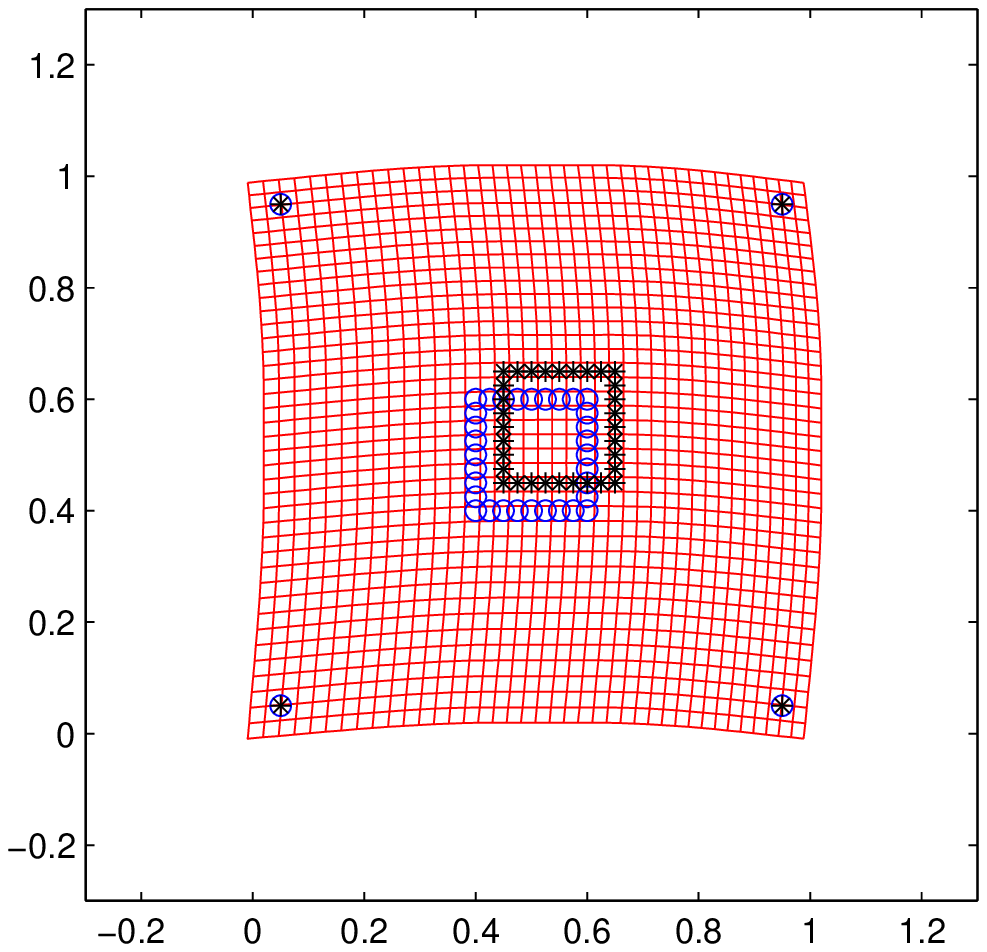}
\centerline{(g) L4, $\alpha=0.2$}
\end{minipage}
\begin{minipage}{60mm}
\includegraphics[width=6.cm]{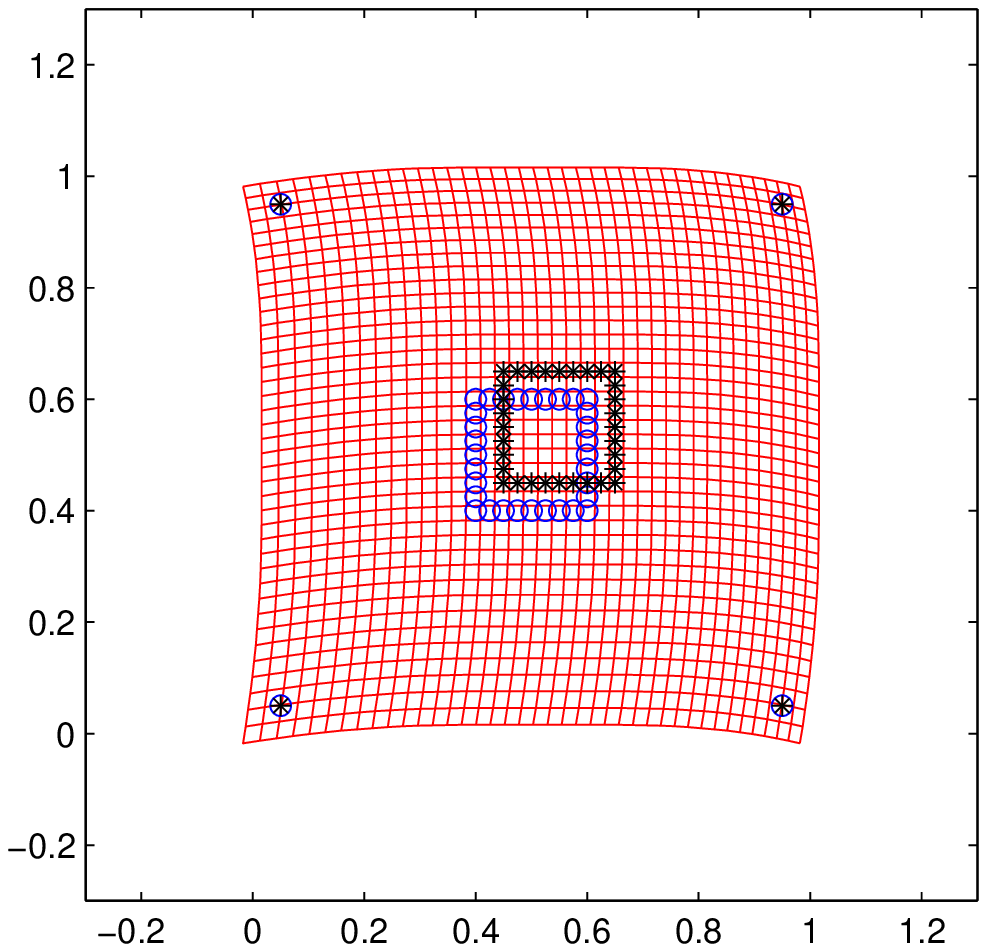}
\centerline{(h) L6, $\alpha=0.4$}
\end{minipage}\\
\end{center}
\caption{Case 1: registration results for the shift of a square using optimal values of $\alpha$ and $c$.}
\label{C1_TEST_opt}
\end{figure}

\begin{figure}[ht!]
\begin{center}
\begin{minipage}{60mm}
\includegraphics[width=6.cm]{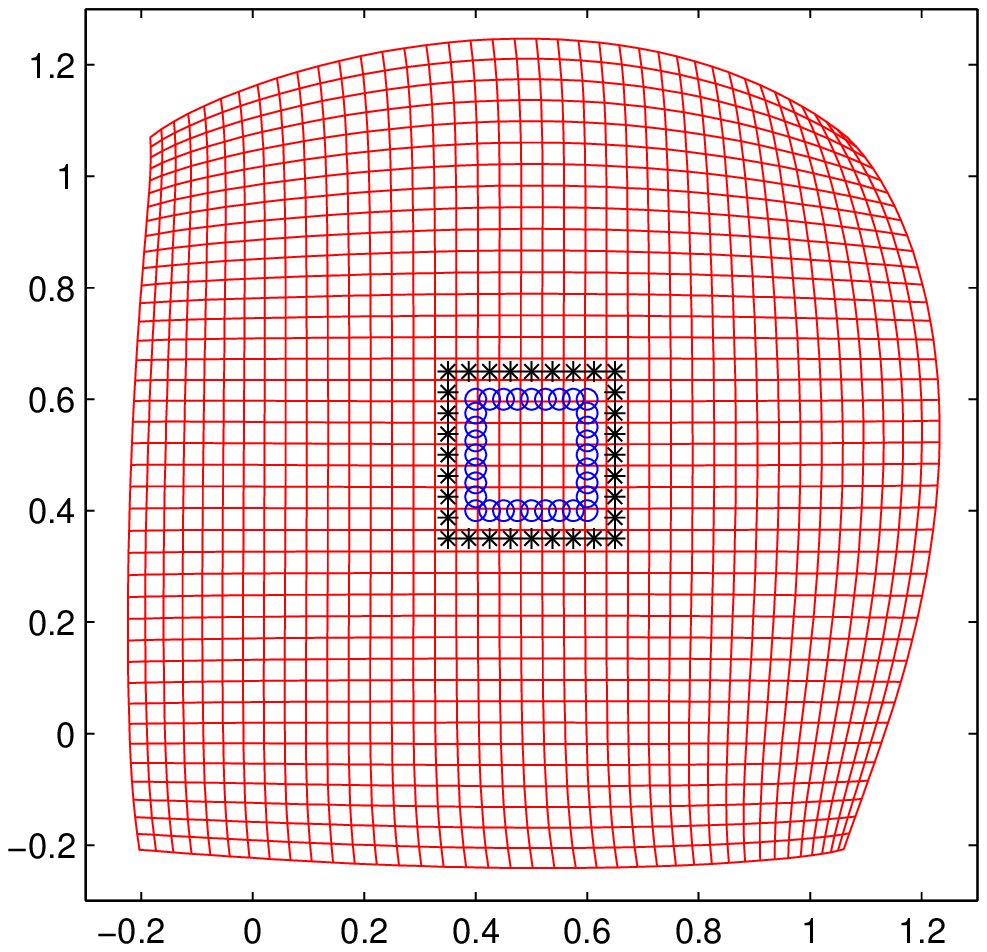}
\centerline{(a) G, $\alpha=2.0$}
\end{minipage}
\begin{minipage}{60mm}
\includegraphics[width=6.cm]{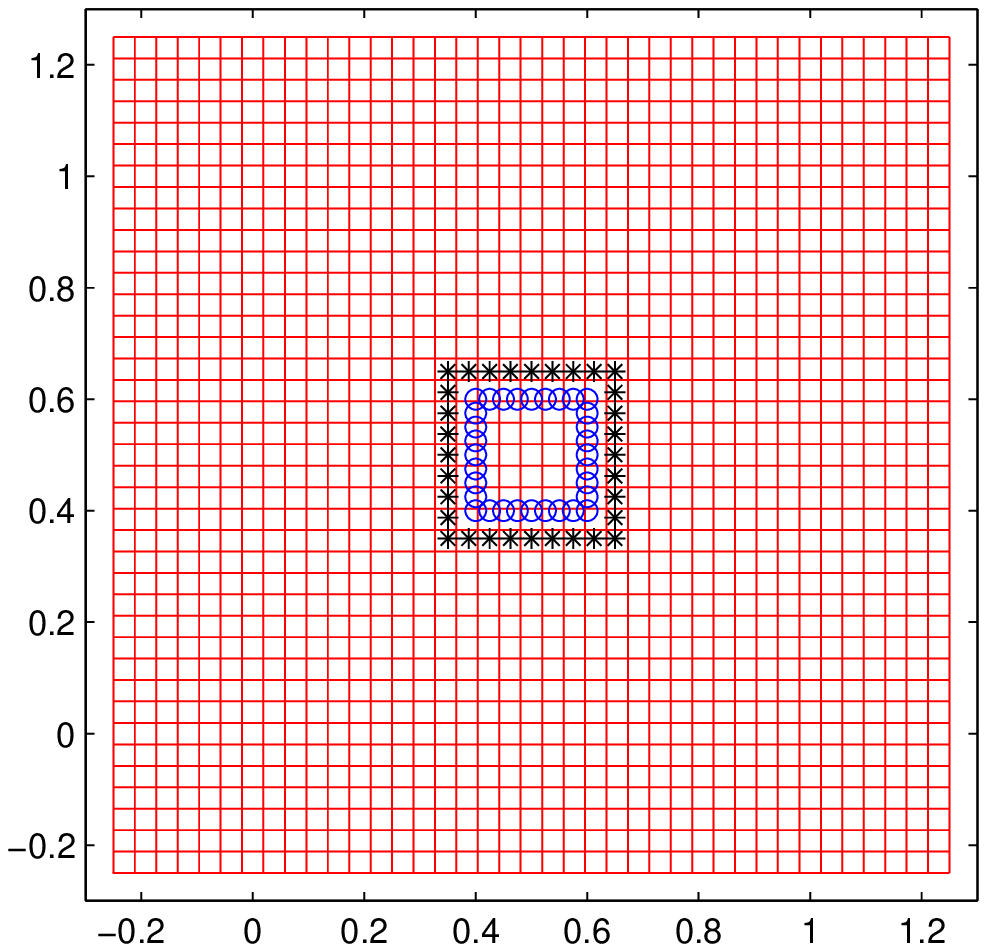}
\centerline{(b) TPS}
\end{minipage}\\
\begin{minipage}{60mm}
\includegraphics[width=6.cm]{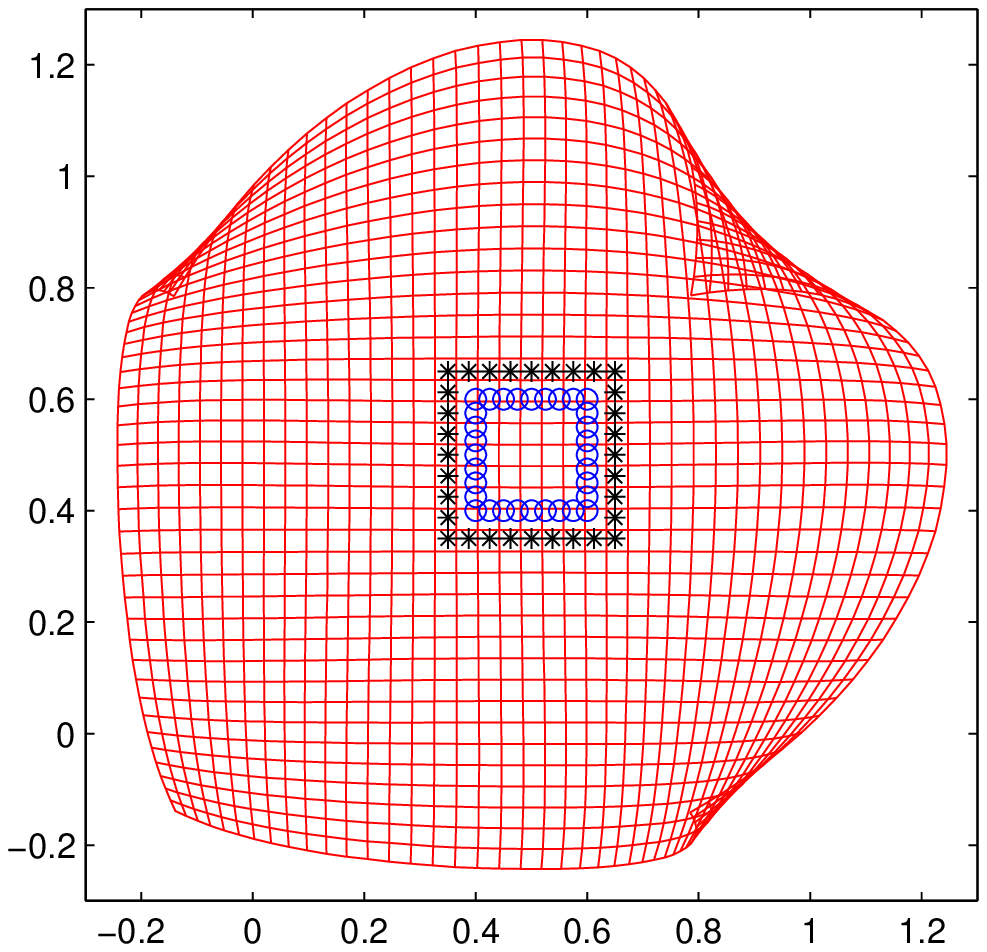}
\centerline{(c) Shep-G, $\alpha=2.0$}
\end{minipage}
\begin{minipage}{60mm}
\includegraphics[width=6.cm]{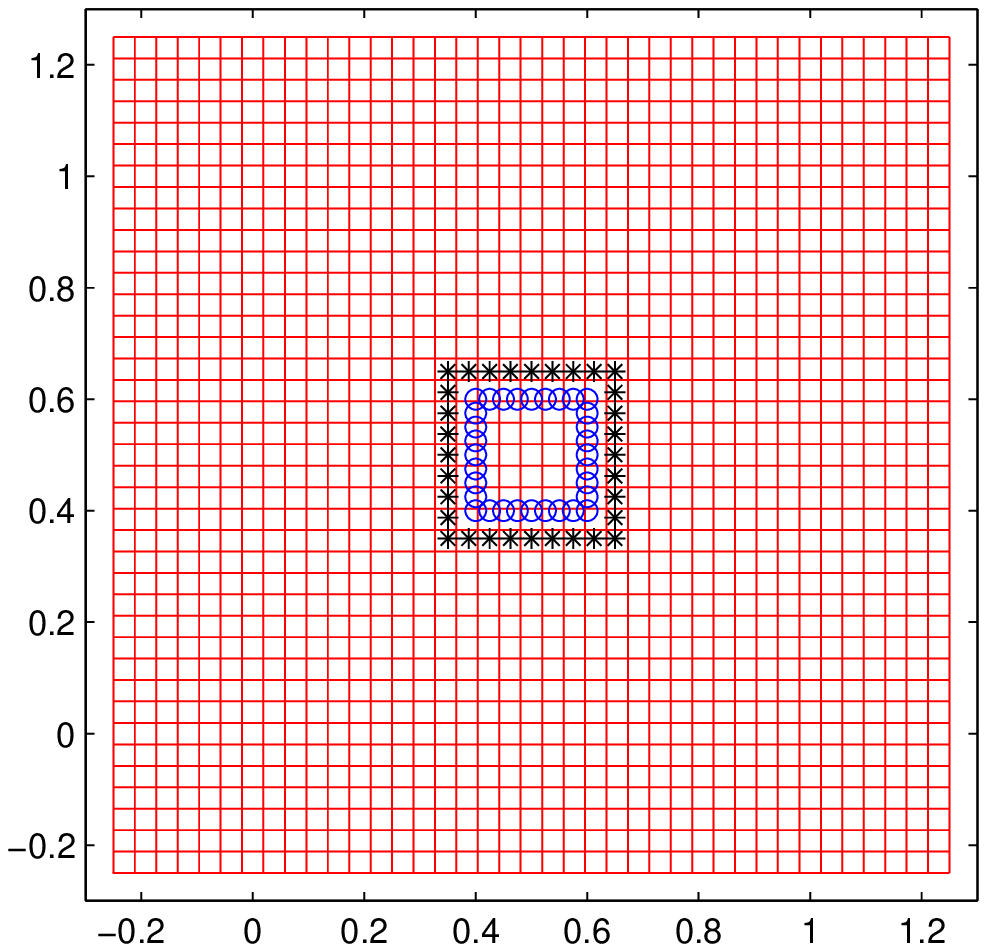}
\centerline{(d) Shep-TPS}
\end{minipage}\\
\begin{minipage}{60mm}
\includegraphics[width=6.cm]{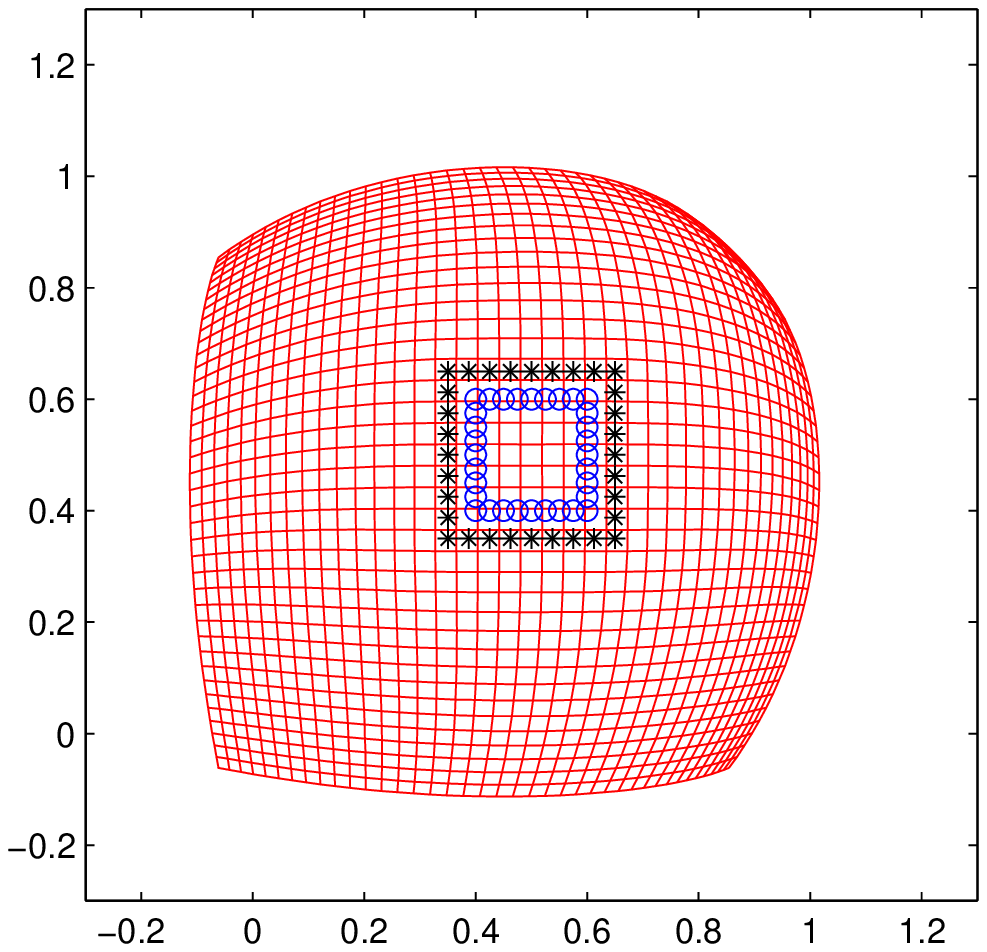}
\centerline{(e) W2-2D, $c=0.3$}
\end{minipage}
\begin{minipage}{60mm}
\includegraphics[width=6.cm]{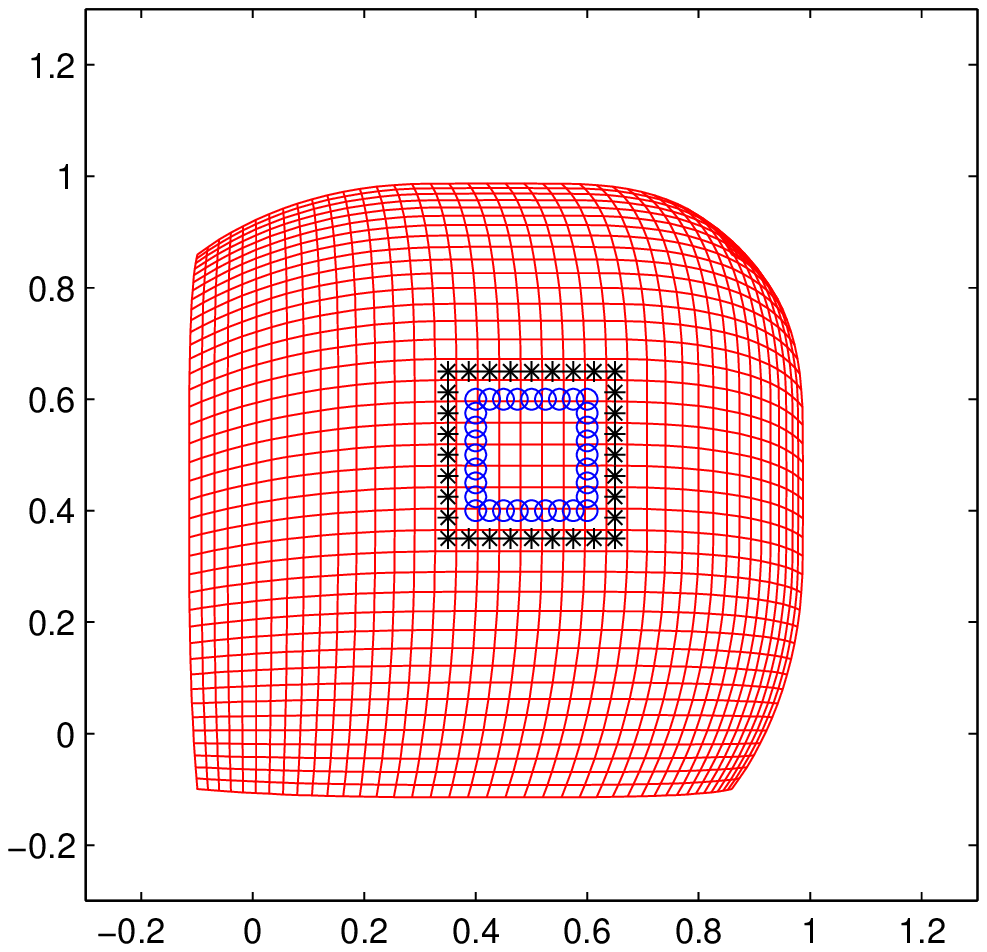}
\centerline{(f) W2-1D$\times$1D, $c=0.4$}
\end{minipage}
\begin{minipage}{60mm}
\includegraphics[width=6.cm]{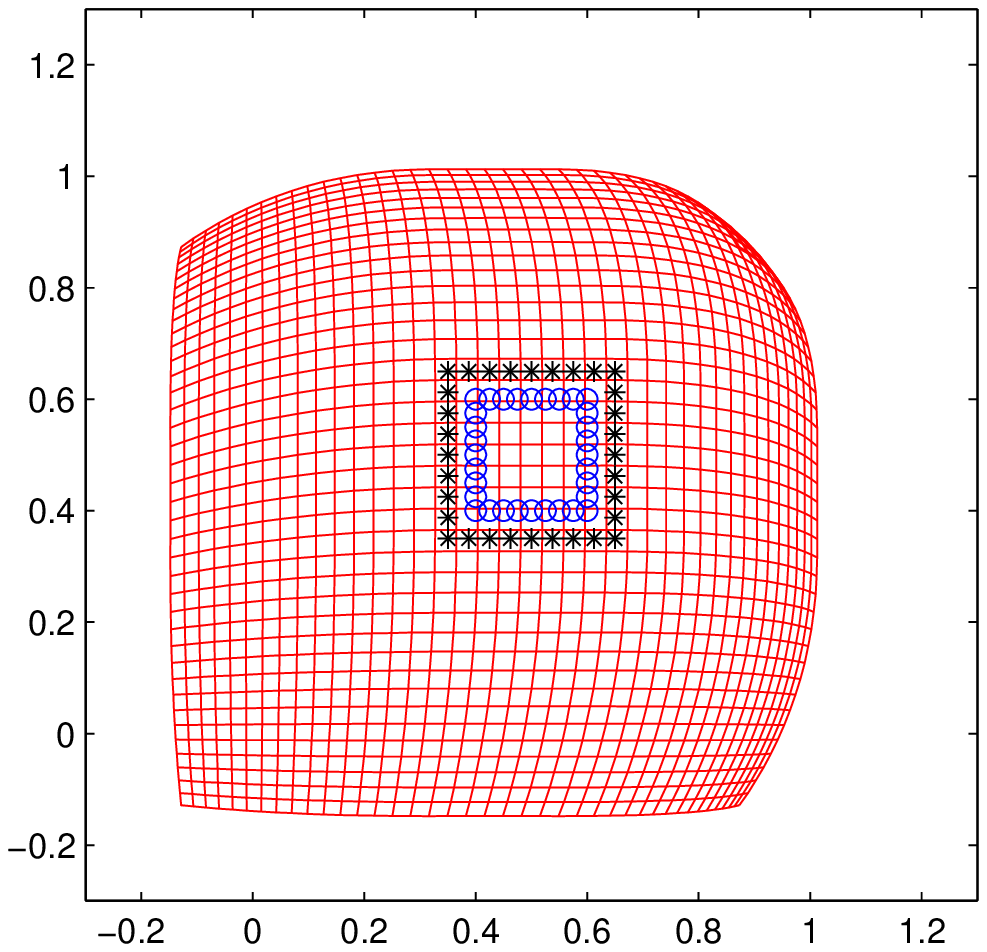}
\centerline{(g) L4, $\alpha=1.4$}
\end{minipage}
\begin{minipage}{60mm}
\includegraphics[width=6.cm]{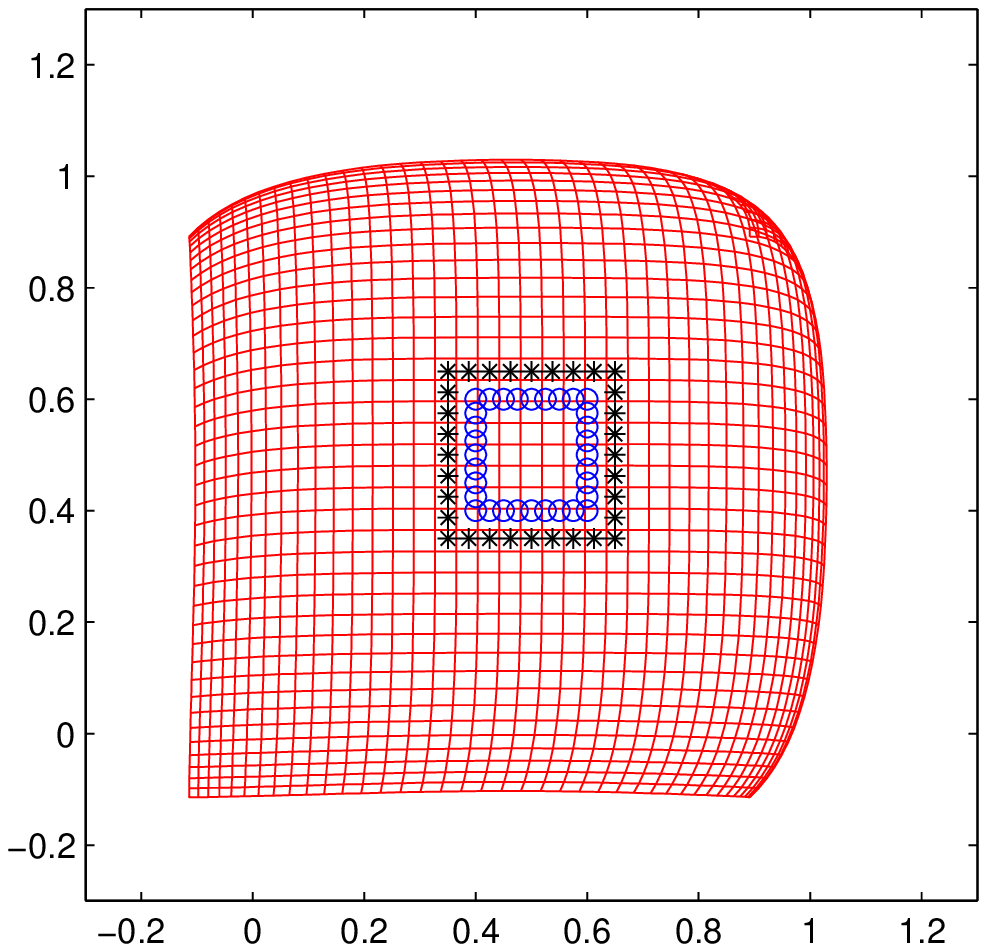}
\centerline{(h) L6, $\alpha=2.0$}
\end{minipage}\\
\end{center}
\caption{Case 2: registration results for the scaling of a square using optimal values of $\alpha$ and $c$.}
\label{C2_TEST_opt}
\end{figure}

\begin{figure}[ht!]
\begin{center}
\begin{minipage}{60mm}
\includegraphics[width=6.cm]{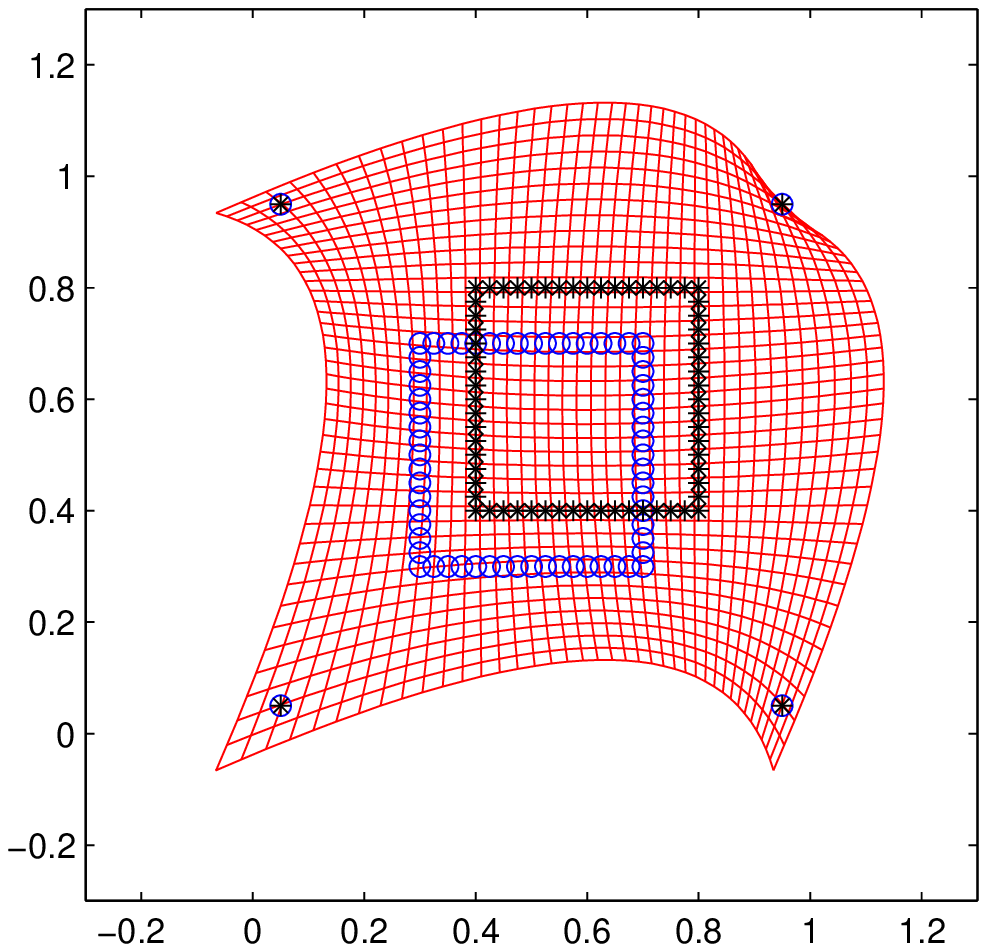}
\centerline{(a) G, $\alpha=0.4$}
\end{minipage}
\begin{minipage}{60mm}
\includegraphics[width=6.cm]{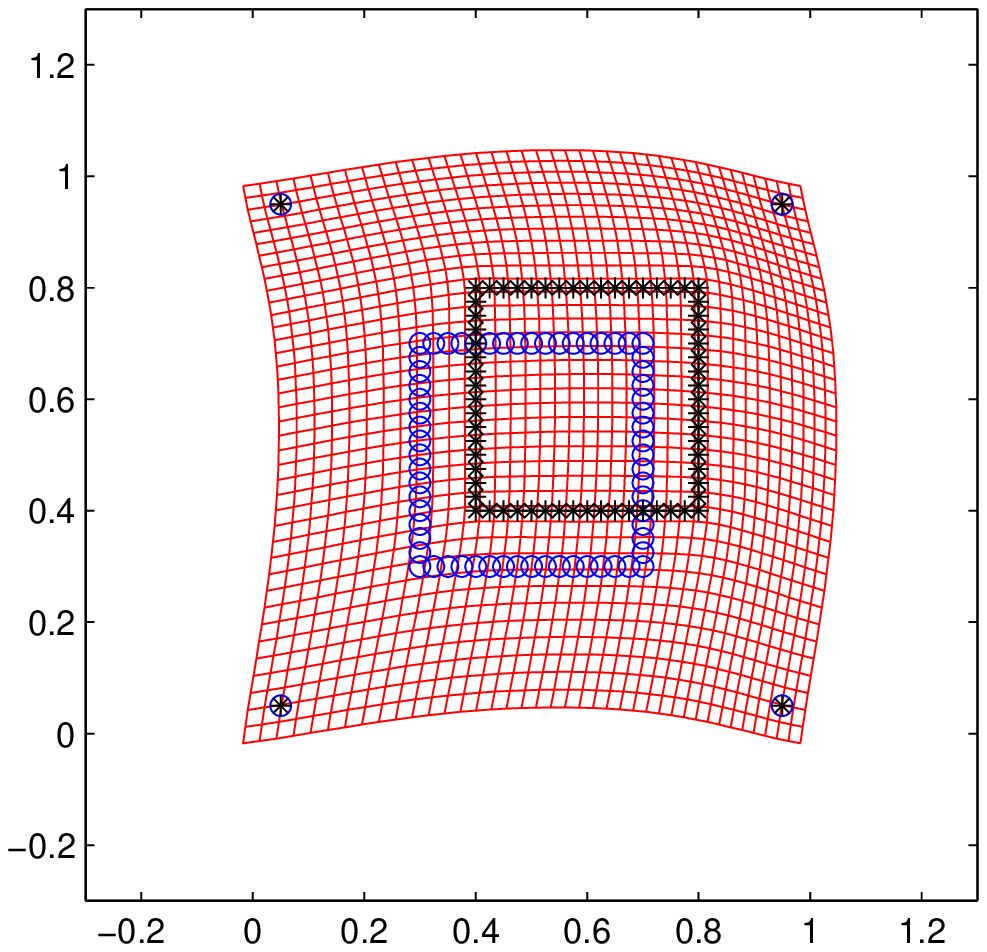}
\centerline{(b) TPS}
\end{minipage}\\
\begin{minipage}{60mm}
\includegraphics[width=6.cm]{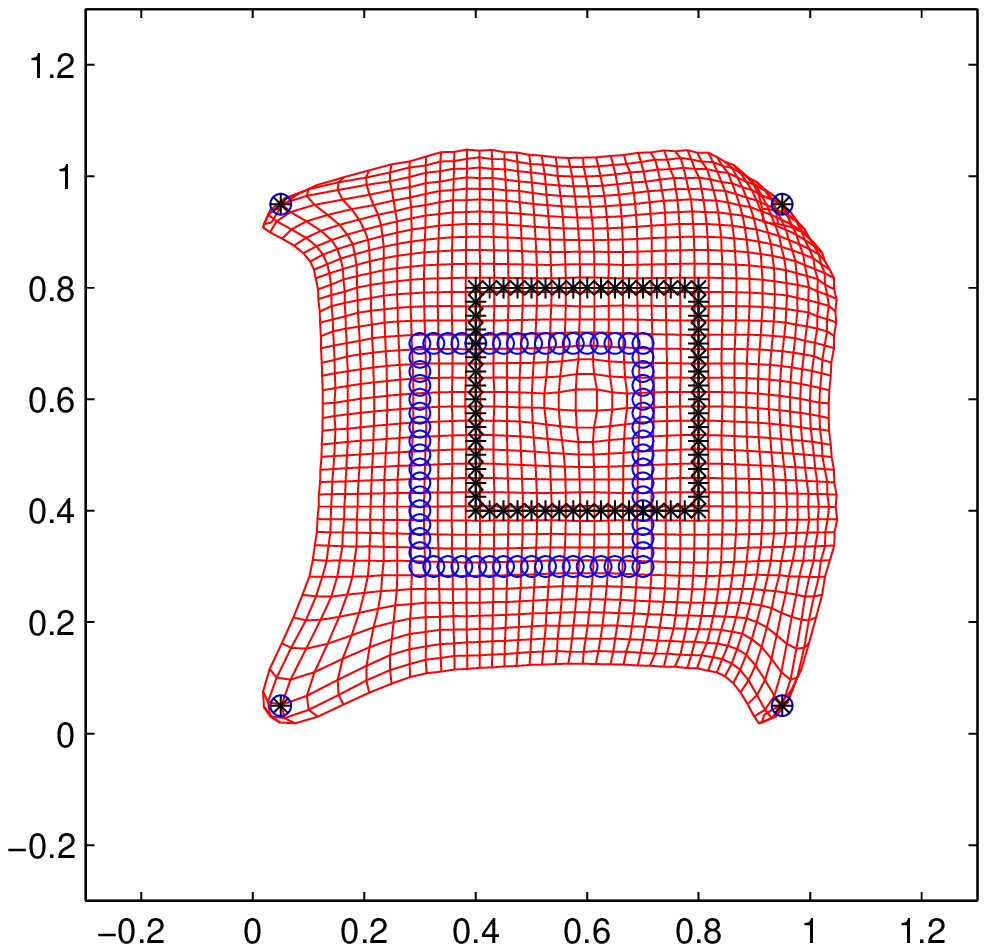}
\centerline{(c) Shep-G, $\alpha=2.0$}
\end{minipage} 
\begin{minipage}{60mm}
\includegraphics[width=6.cm]{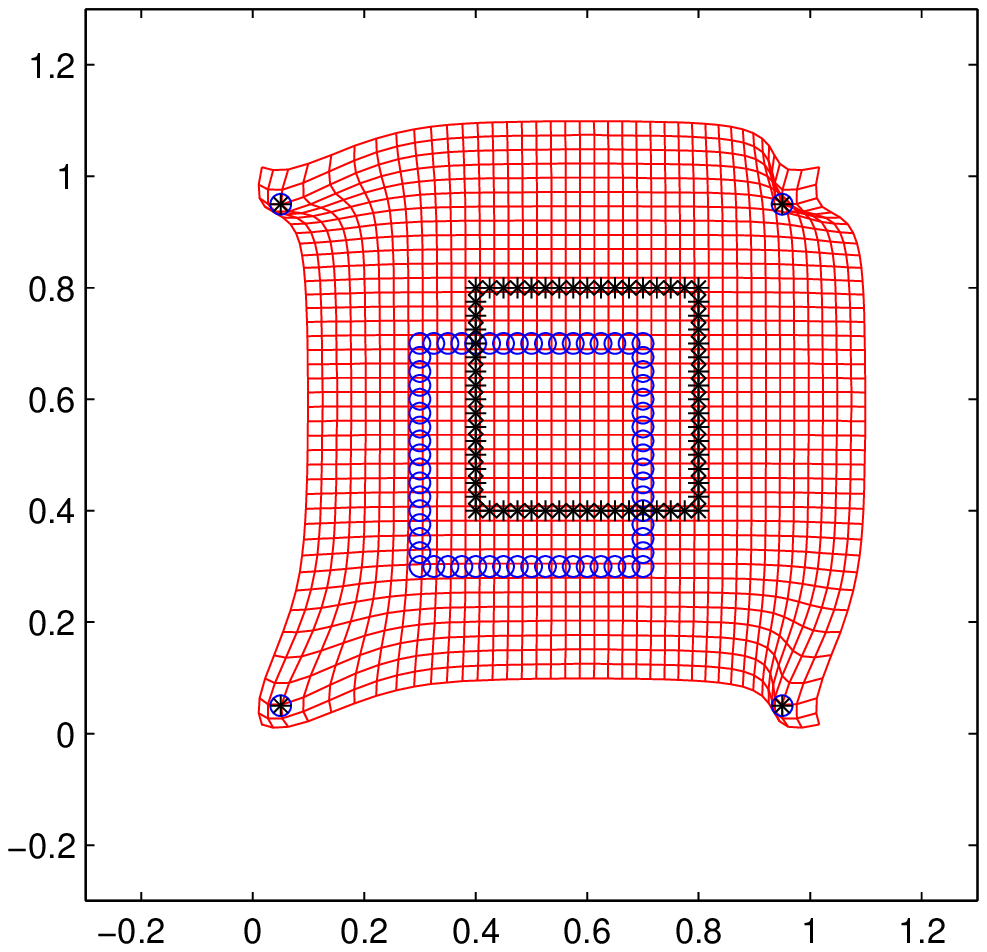}
\centerline{(d) Shep-TPS}
\end{minipage}\\
\begin{minipage}{60mm}
\includegraphics[width=6.cm]{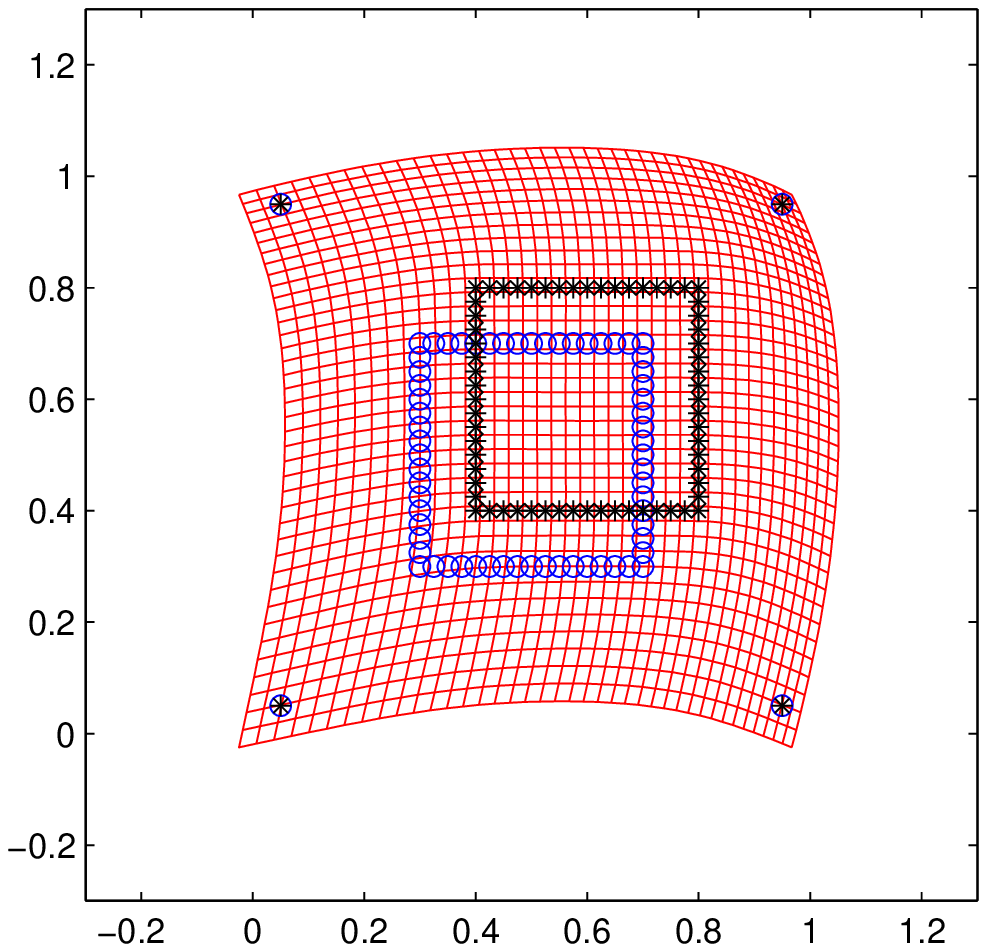}
\centerline{(e) W2-2D, $c=0.1$}
\end{minipage}
\begin{minipage}{60mm}
\includegraphics[width=6.cm]{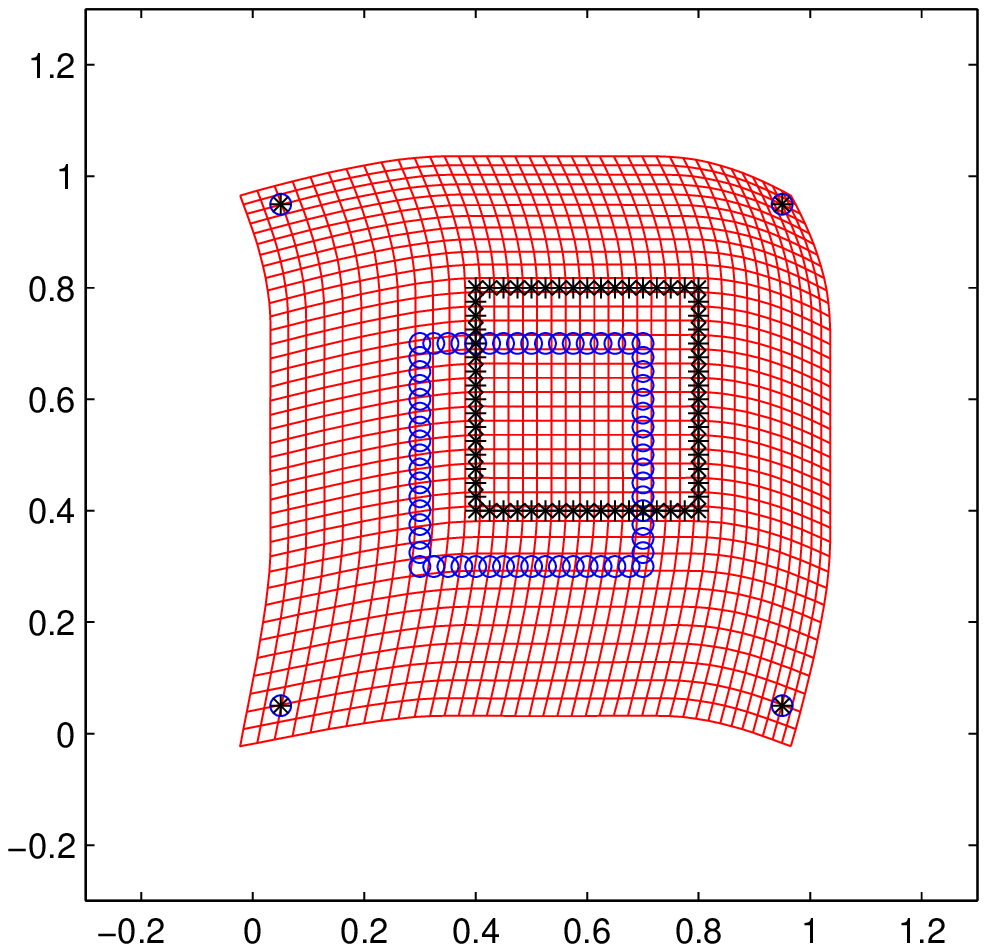}
\centerline{(f) W2-1D$\times$1D, $c=0.2$}
\end{minipage}
\begin{minipage}{60mm}
\includegraphics[width=6.cm]{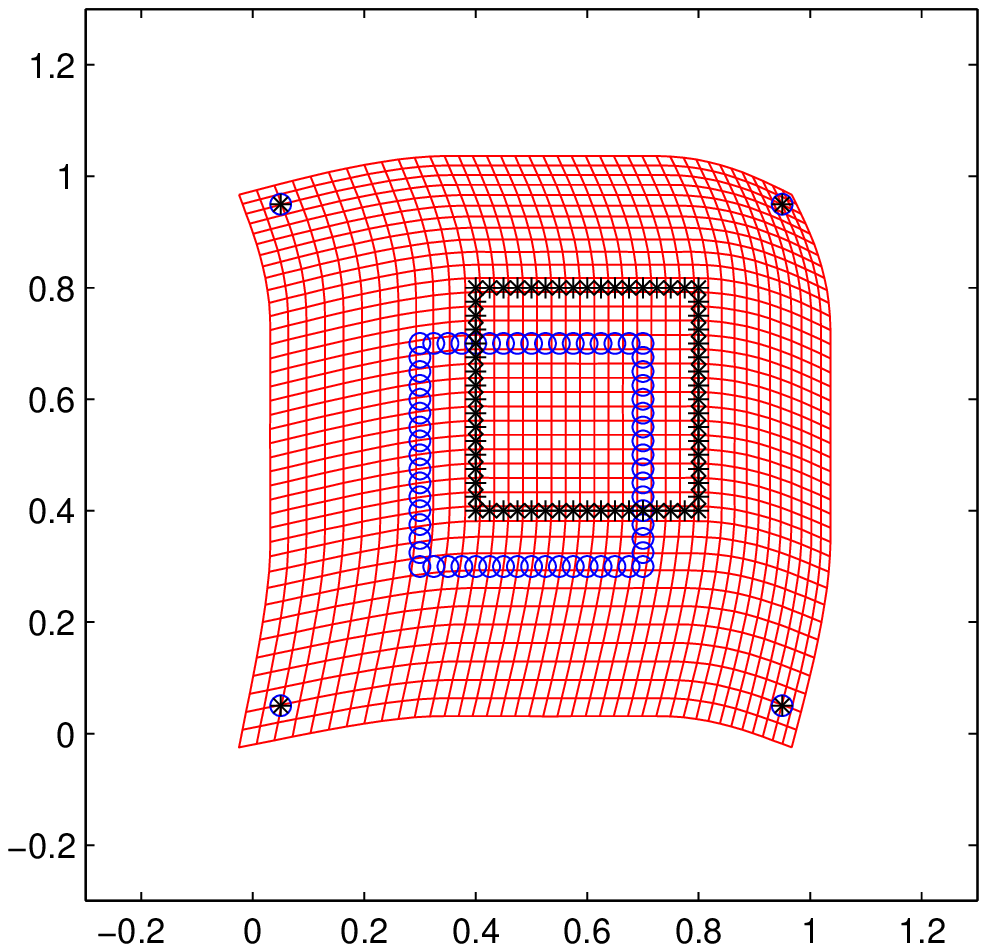}
\centerline{(g) L4, $\alpha=0.6$}
\end{minipage}
\begin{minipage}{60mm}
\includegraphics[width=6.cm]{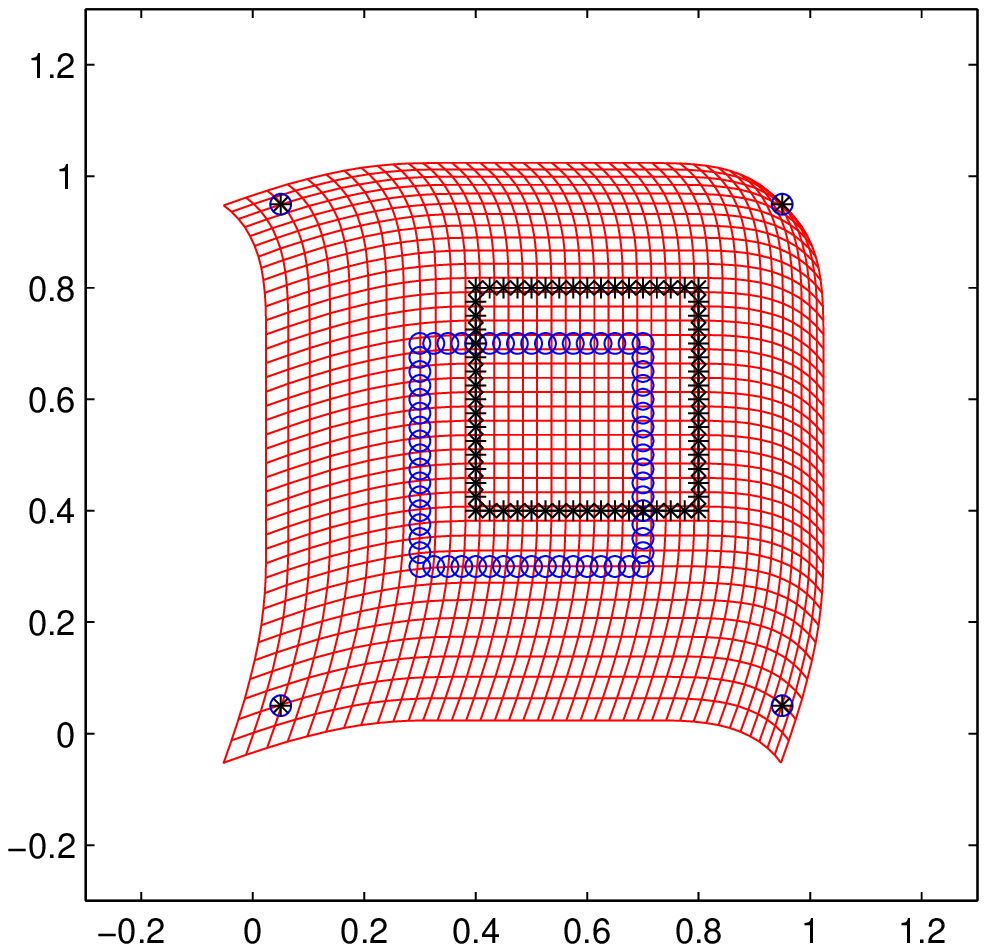}
\centerline{(h) L6, $\alpha=0.2$}
\end{minipage}\\
\end{center}
\caption{Case 3: registration results for the shift of a square using optimal values of $\alpha$ and $c$.}
\label{C3_TEST_opt}
\end{figure}

\begin{figure}[ht!]
\begin{center}
\begin{minipage}{60mm}
\includegraphics[width=6.cm]{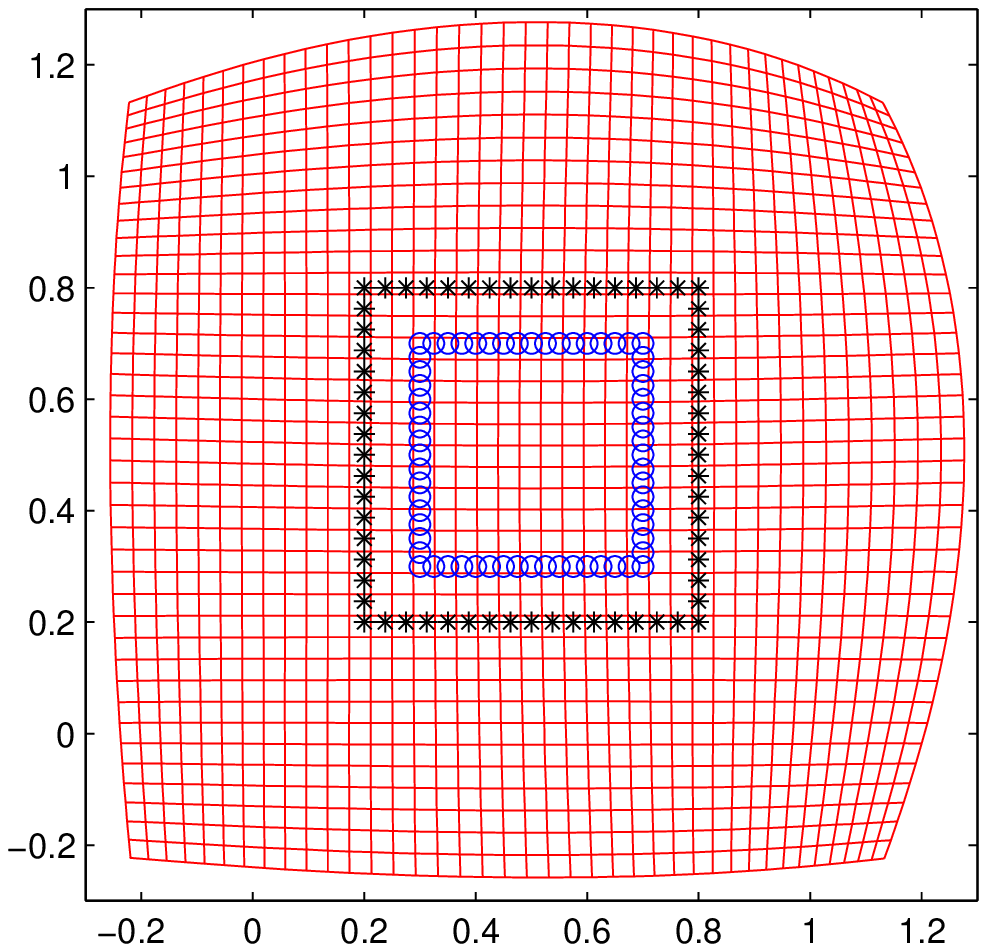}
\centerline{(a) G, $\alpha=2.0$}
\end{minipage}
\begin{minipage}{60mm}
\includegraphics[width=6.cm]{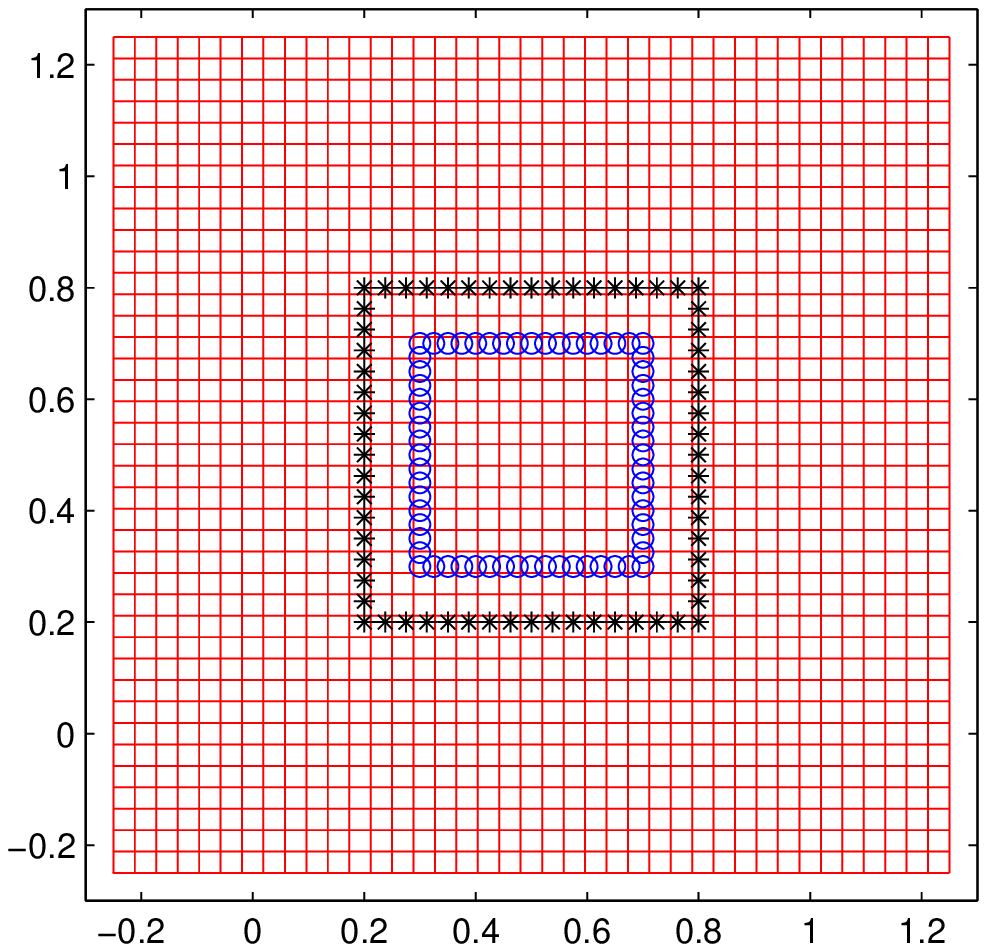}
\centerline{(b) TPS}
\end{minipage}\\
\begin{minipage}{60mm}
\includegraphics[width=6.cm]{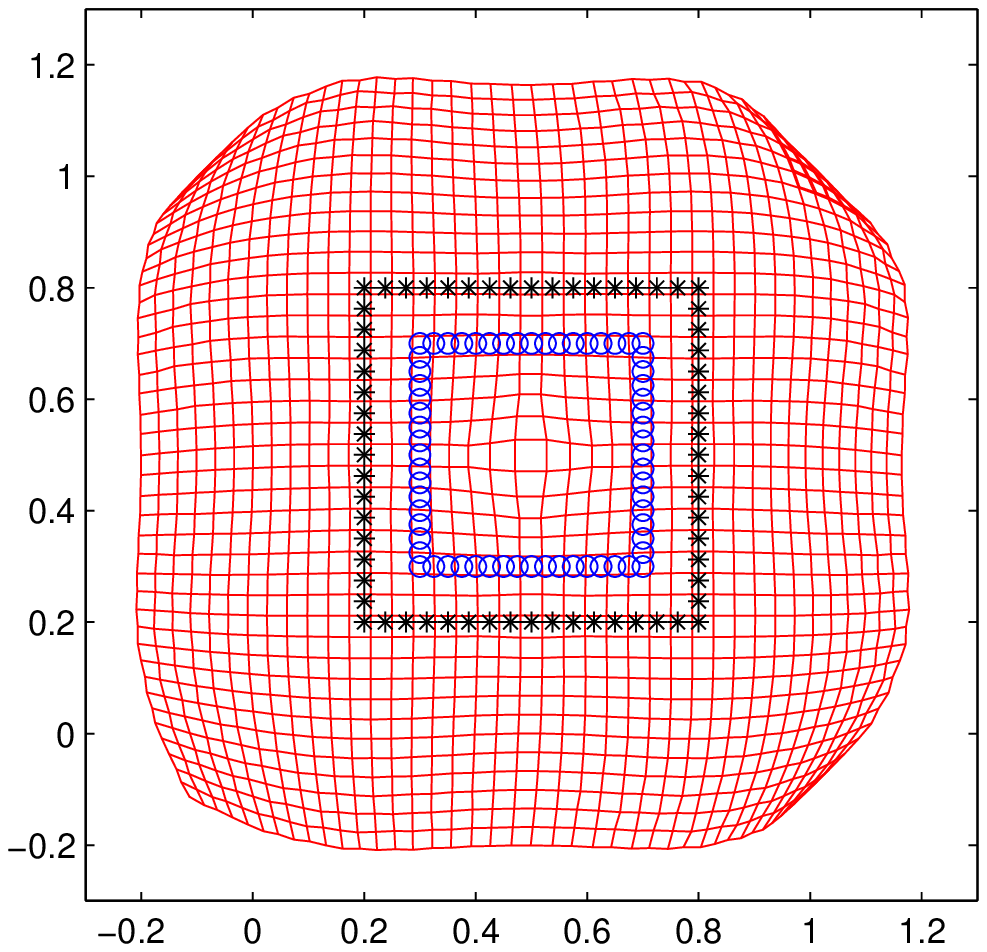}
\centerline{(c) Shep-G, $\alpha=2.0$}
\end{minipage}
\begin{minipage}{60mm}
\includegraphics[width=6.cm]{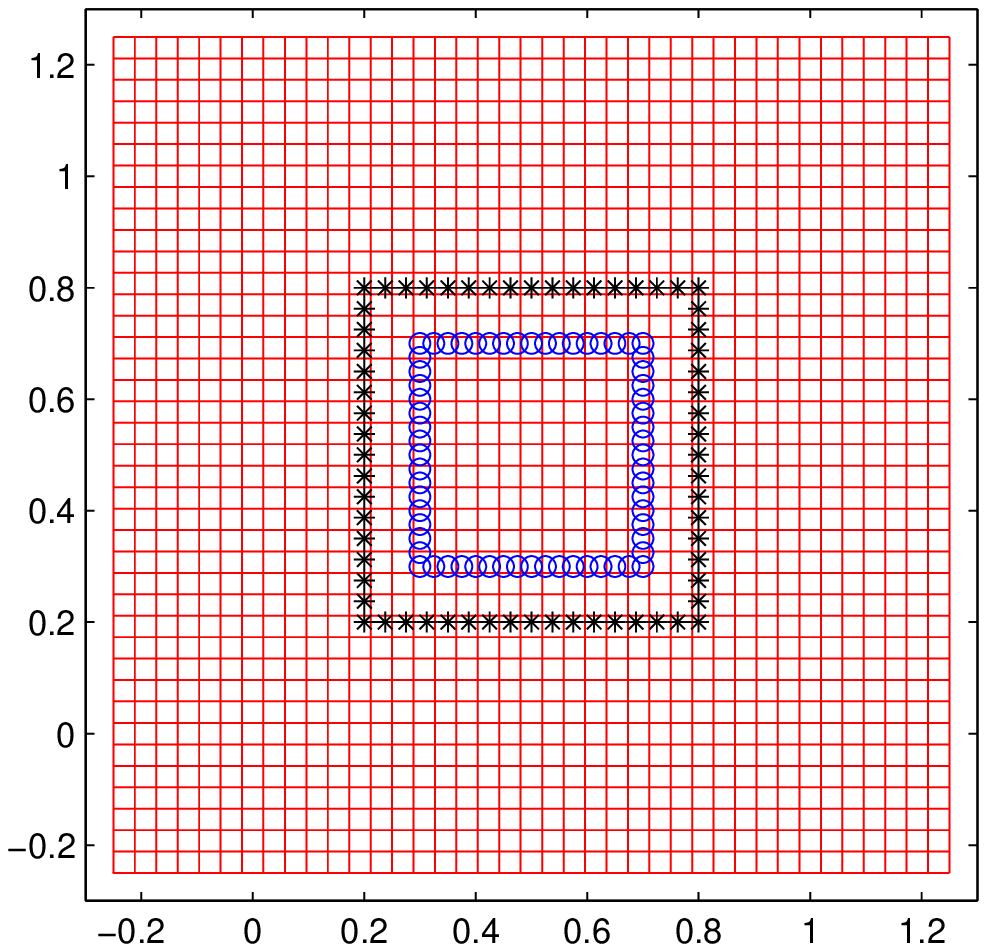}
\centerline{(d) Shep-TPS}
\end{minipage}\\
\begin{minipage}{60mm}
\includegraphics[width=6.cm]{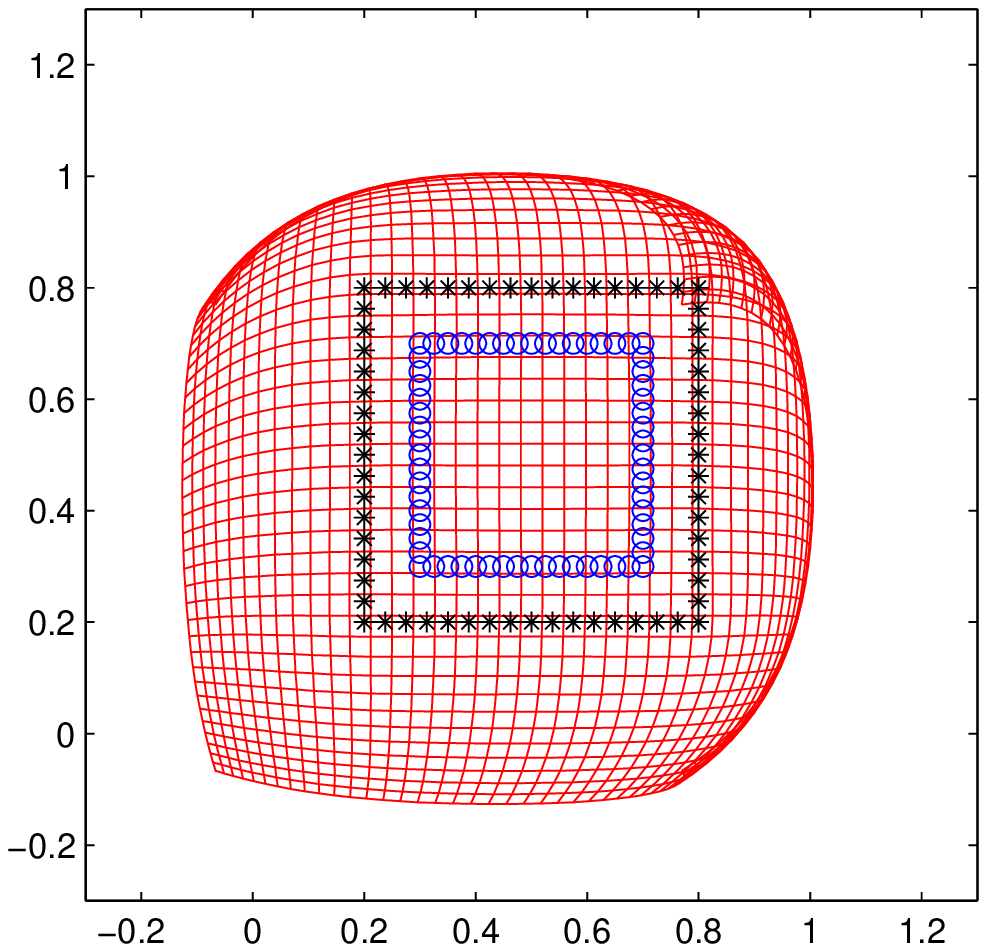}
\centerline{(e) W2-2D, $c=0.5$}
\end{minipage}
\begin{minipage}{60mm}
\includegraphics[width=6.cm]{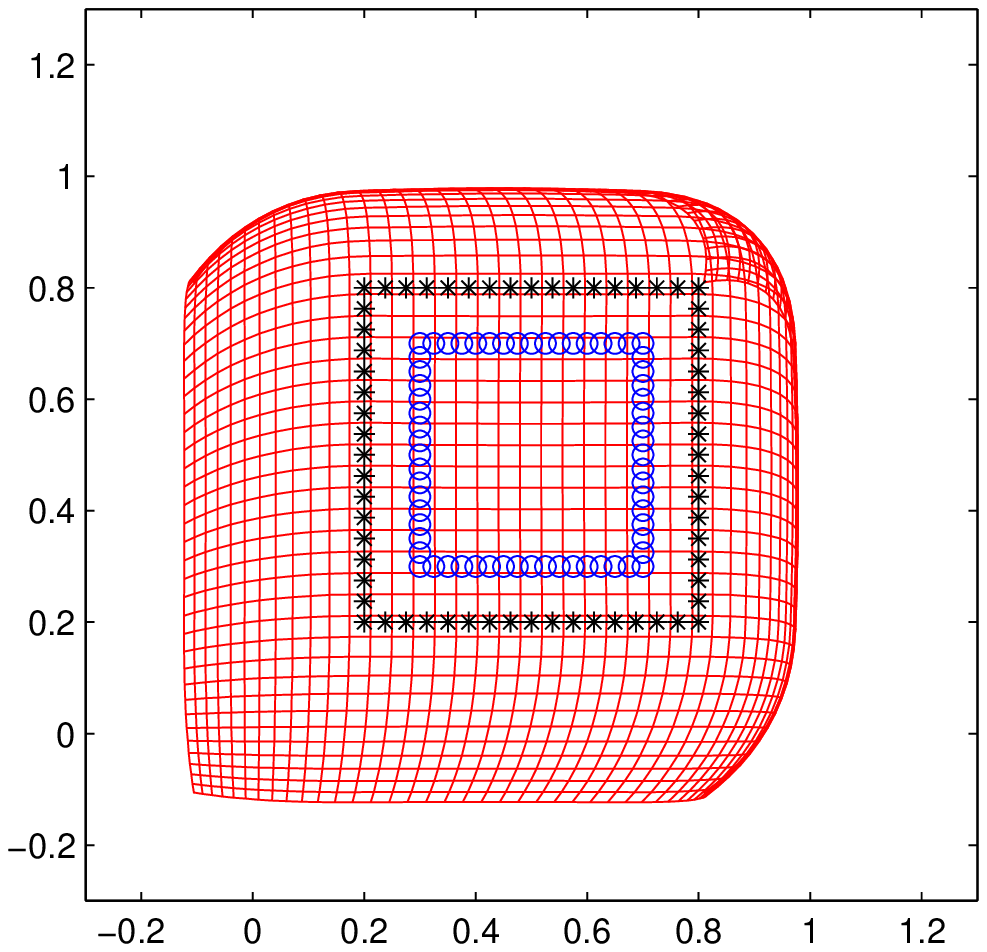}
\centerline{(f) W2-1D$\times$1D, $c=0.6$}
\end{minipage}
\begin{minipage}{60mm}
\includegraphics[width=6.cm]{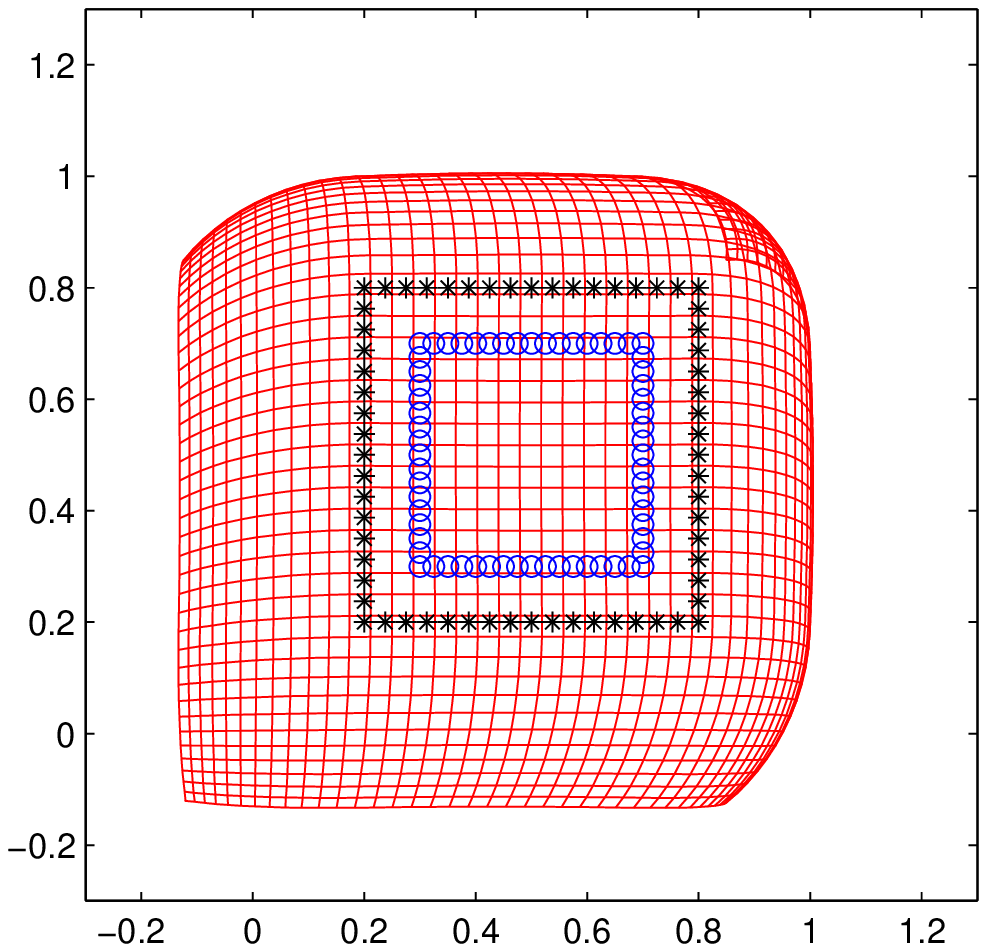}
\centerline{(g) L4, $\alpha=2.0$}
\end{minipage}
\begin{minipage}{60mm}
\includegraphics[width=6.cm]{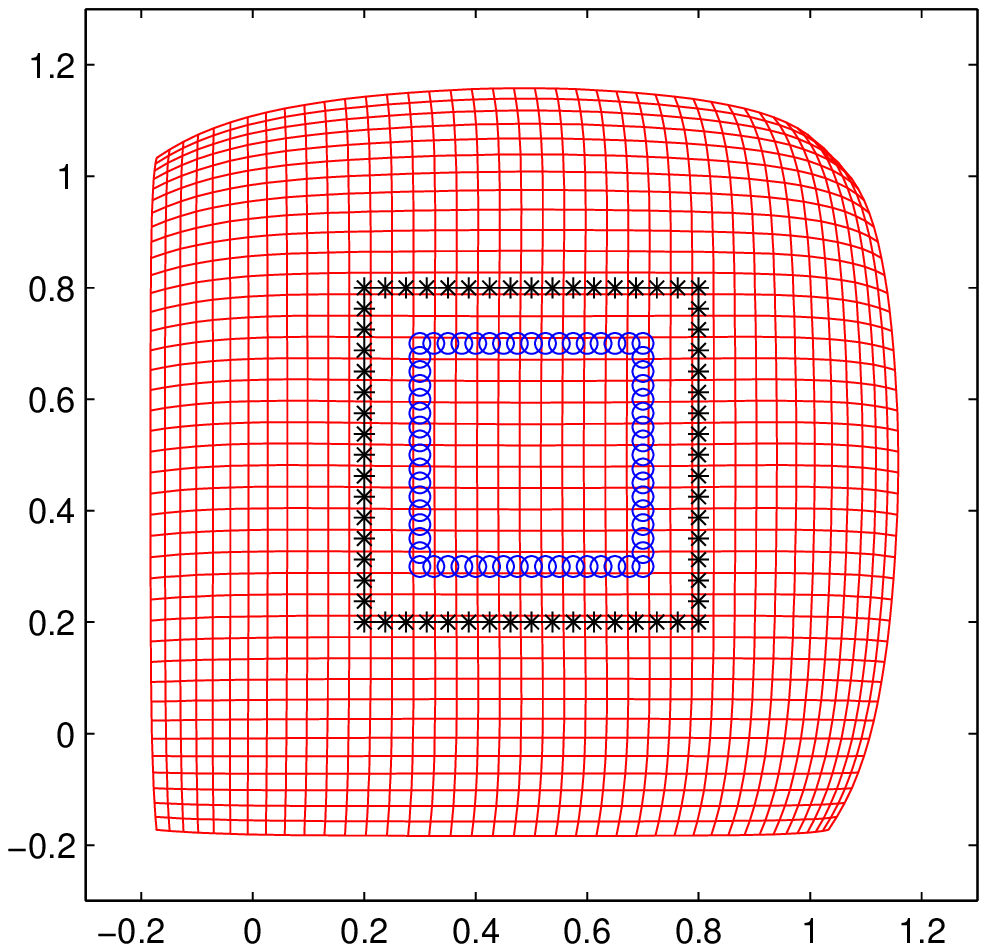}
\centerline{(h) L6, $\alpha=2.0$}
\end{minipage}\\
\end{center}
\caption{Case 4: registration results for the scaling of a square using optimal values of $\alpha$ and $c$.}
\label{C4_TEST_opt}
\end{figure}

Finally, changing the values of shape parameters, sometimes it is possible to obtain better registration results from the point of view of smoothness. In particular, we focus our attention on compactly supported functions, which seem to be more effective for image registration. In Figure 11 we show some transformed grids using Lobachevsky L4 and Wendland W2-2D where the values of parameters are not optimal choices with regard of errors ($\alpha=1.0$ and $c=0.2$, respectively), but such that the registration results are smoother. RMSEs are given by $1.3931\text{E}-1$ and $1.2132\text{E}-1$, respectively, and therefore are greater than the errors found using optimal values for the shape parameters.

Summing all results and considerations on these test cases, we can say that Wendland's functions work particularly well when $c\in [0.1,0.5]$. In fact, in this case, the behavior of RMSEs varying $c$ shows that accuracy is preserved (see Figures 3--6). For the Gaussian the RMSE graphs are flat and all parameter values give equivalent results, but the accuracy is slight better for $\alpha\in[1,2]$. We observe the same behavior for the graphs of Lobachevsky L4 and L6. Precisely, in test Cases 1 and 3 it is better to choose the parameter value in the interval $[0.1,1]$, whereas in test Cases 2 and 4 the choice $\alpha\in[1,2]$ seems better. However, all the considered values give a good approximation and the errors vary very little.

\begin{figure}[ht!]
\begin{center}
\begin{minipage}{60mm}
\includegraphics[width=6.cm]{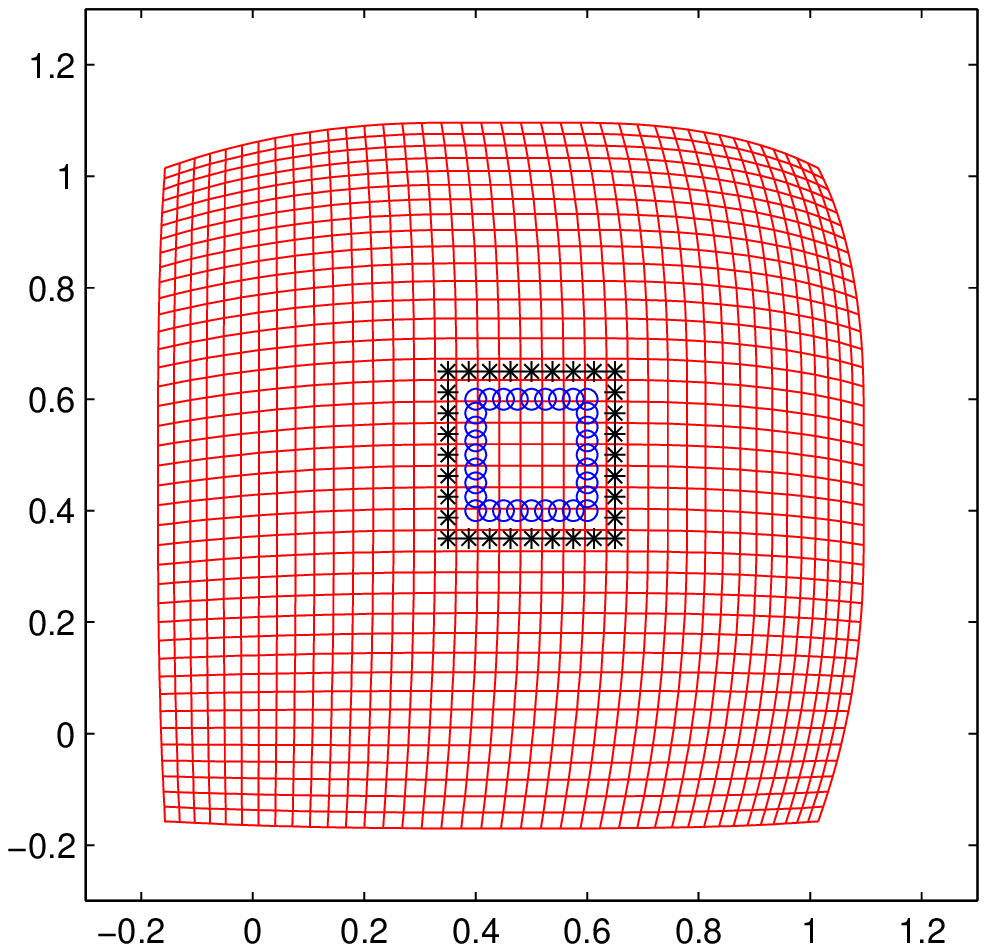}
\centerline{L4}
\end{minipage}
\begin{minipage}{60mm}
\includegraphics[width=6.cm]{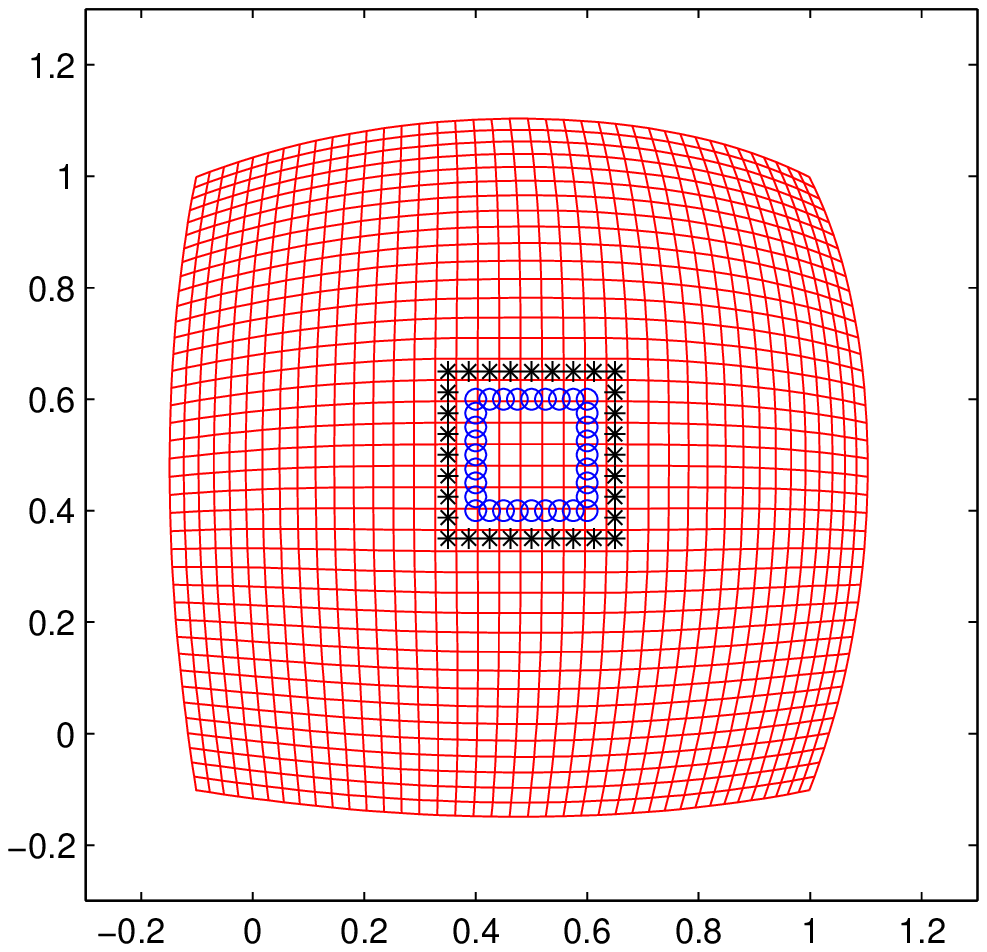}
\centerline{W2-2D}
\end{minipage}
\end{center}
\caption{Case 2, registration results for the scaling of a square, choice of non-optimal parameters.}
\label{C2_TEST}
\end{figure}

\subsection{Test example 2: circle expansion and contraction}
\begin{figure}[ht!]
\begin{center}
  \includegraphics[width=6.cm]{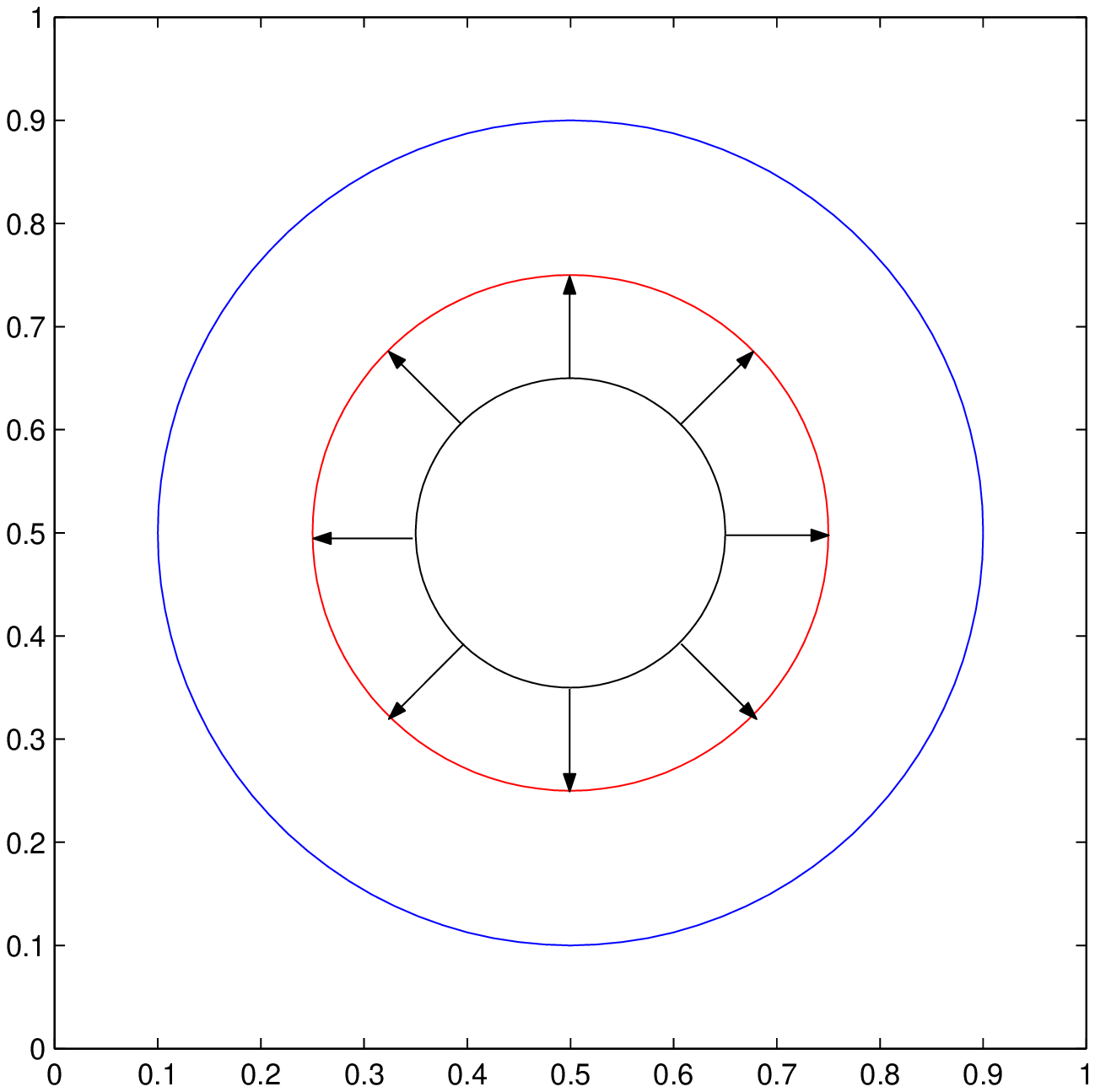}\hfil
  \includegraphics[width=6.cm]{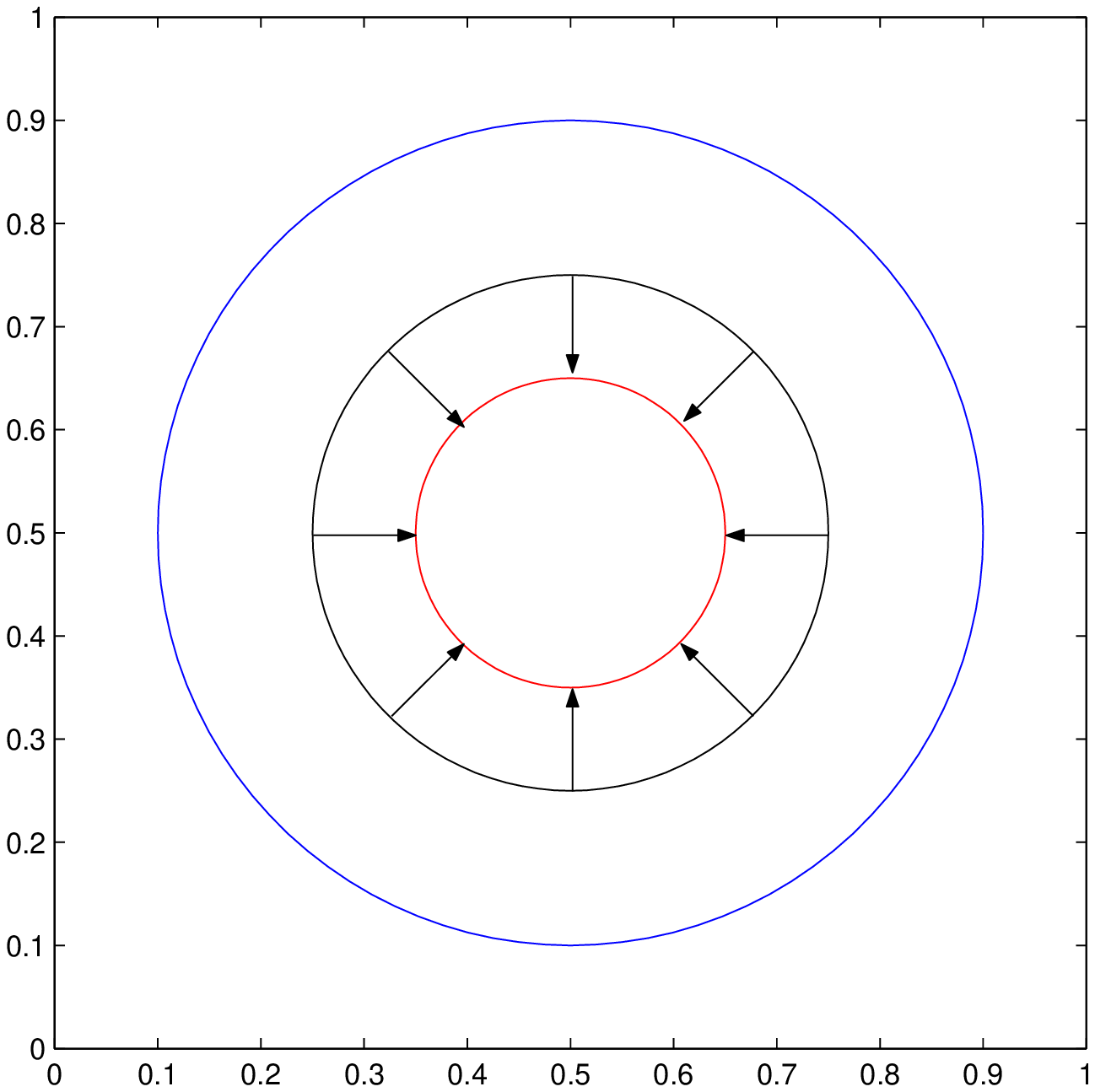}
\end{center}  
  \caption{Models of an expanding tumor (left) and a tumor resection (right).}
  \label{fig:tumor}
\end{figure}

In this subsection we consider two opposite radial transformations, that is, the expansion and the contraction of a circle (see Figure \ref{fig:tumor}). They may offer  very schematic models for the growing and the resection of a tumor in surrounding elastic brain tissue. In these models (given in \cite{kohlrausch05}) the outer circle corresponds to the skull bone, which is assumed to be rigid. The inner circle represents the boundary of the tumor, whereas the space between the inner and the outer circle is assumed to be filled with elastic material, which corresponds to brain tissue.

\begin{figure}[ht!]
\begin{center}
\includegraphics[width=6.cm]{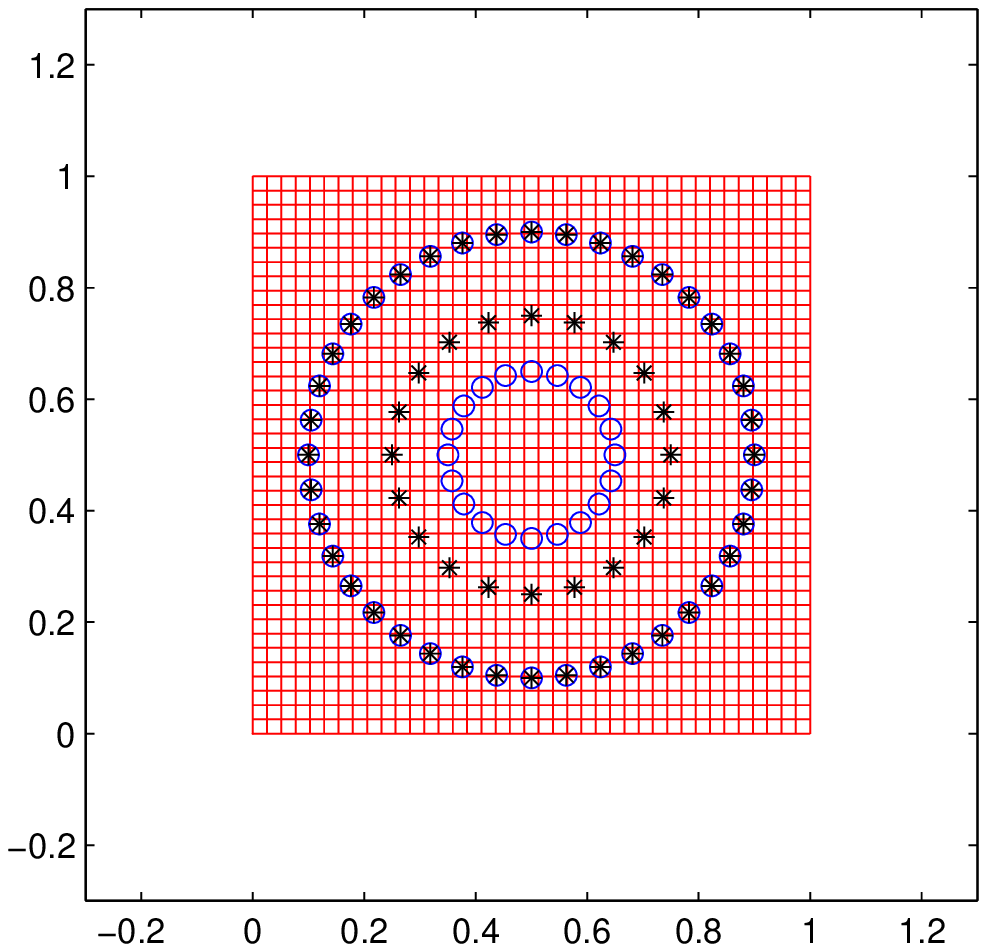}
\includegraphics[width=6.cm]{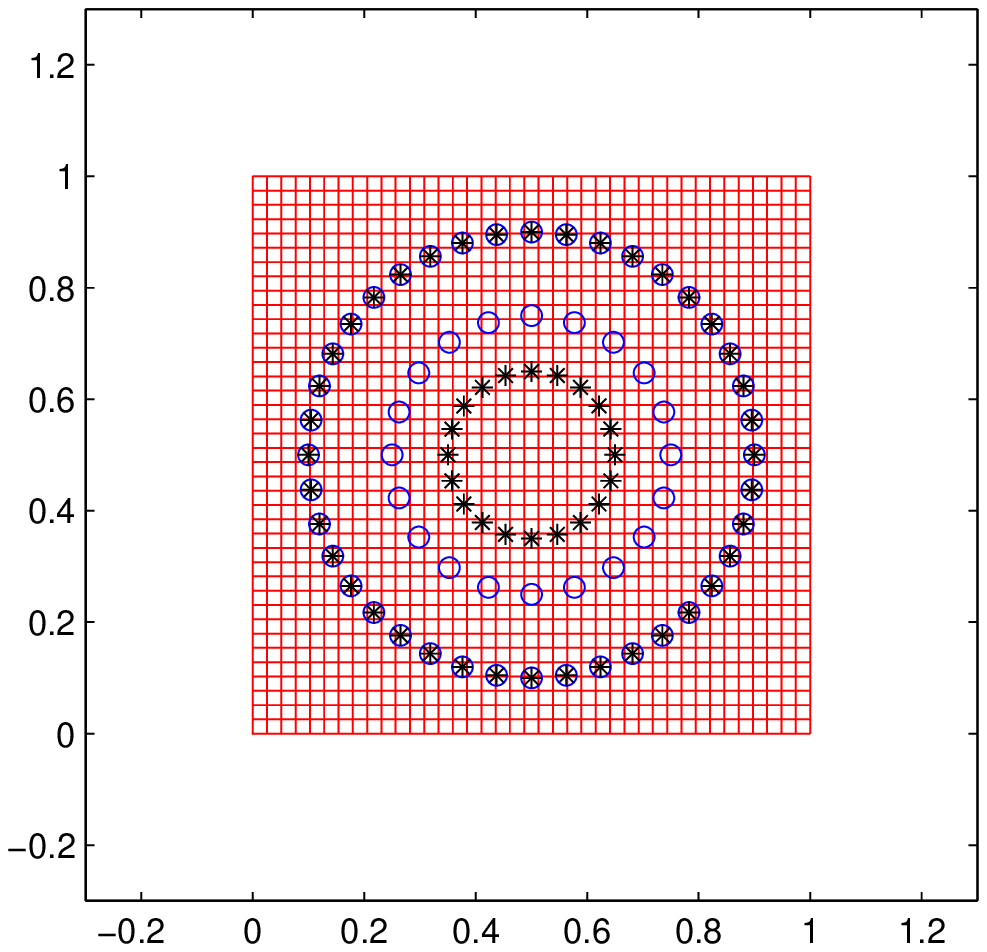}
\end{center}  
  \caption{Circle expansion (left) and circle contraction (right): source and target landmarks.}
\label{TUMOR_EXP-RES}
\end{figure}

The grids are transformed using 20 equidistant landmarks placed on the inner circle and, to prevent an overall shift, also 40 quasi--landmarks, i.e. landmarks at invariant positions, at the outer circle in the source and target images. These point-landmarks, shown in Figure \ref{TUMOR_EXP-RES} for the circle expansion (left) and the circle contraction (right), are marked by a circle ($\circ$) and a star ($\star$), respectively.

In these experiments we compare at first the registration results obtained by using the modified Shepard's method, with TPSs as local interpolants, and the global TPS method. 

Hence, taking optimal values for $N_L$ and $N_W$, in Figure 14 we show the registration results obtained by applying TPS (a) and Shep-TPS (b). Their comparison mainly points out the goodness and the effectiveness of Shepard-type method for this situation (see also \cite{Cavoretto-DeRossi08a,Cavoretto-DeRossi08b}). Specifically, in the circle expansion we consider $N_L=16$ and $N_W=60$ for the modified Shepard's formula, obtaining a RMSE = $4.5853{\rm E}-2$, while the global approach produces a RMSE = $6.0964{\rm E}-2$. On the other hand, in the circle contraction the local method yields a RMSE = $3.5643{\rm E}-2$ for $N_L=5$ and $N_W=60$, whereas the global one gives a RMSE = $7.7354{\rm E}-2$.

The following Figures 14 (c)--(d) show registration results using Wendland W2-2D ($c=0.1$ and $c=0.4$ for expansion and contraction, respectively)  and W2-1D$\times$1D ($c=0.7$ and $c=0.8$, respectively). In particular, in the circle expansion we obtain a RMSE = $9.1226{\rm E}-2$ for W2-2D and a RMSE = $1.5294{\rm E}-1$ for W2-1D$\times$1D, whereas in the circle contraction we have a RMSE = $6.4178{\rm E}-2$ for W2-2D and a RMSE = $7.8952{\rm E}-2$ for W2-1D$\times$1D. We note that deformations of transformed images using Wendland's functions are very large. For the other methods we observe even larger deformations and therefore we present no results.

\begin{figure}[ht!]
\begin{center}

\begin{minipage}{60mm}
  \includegraphics[width=6.cm]{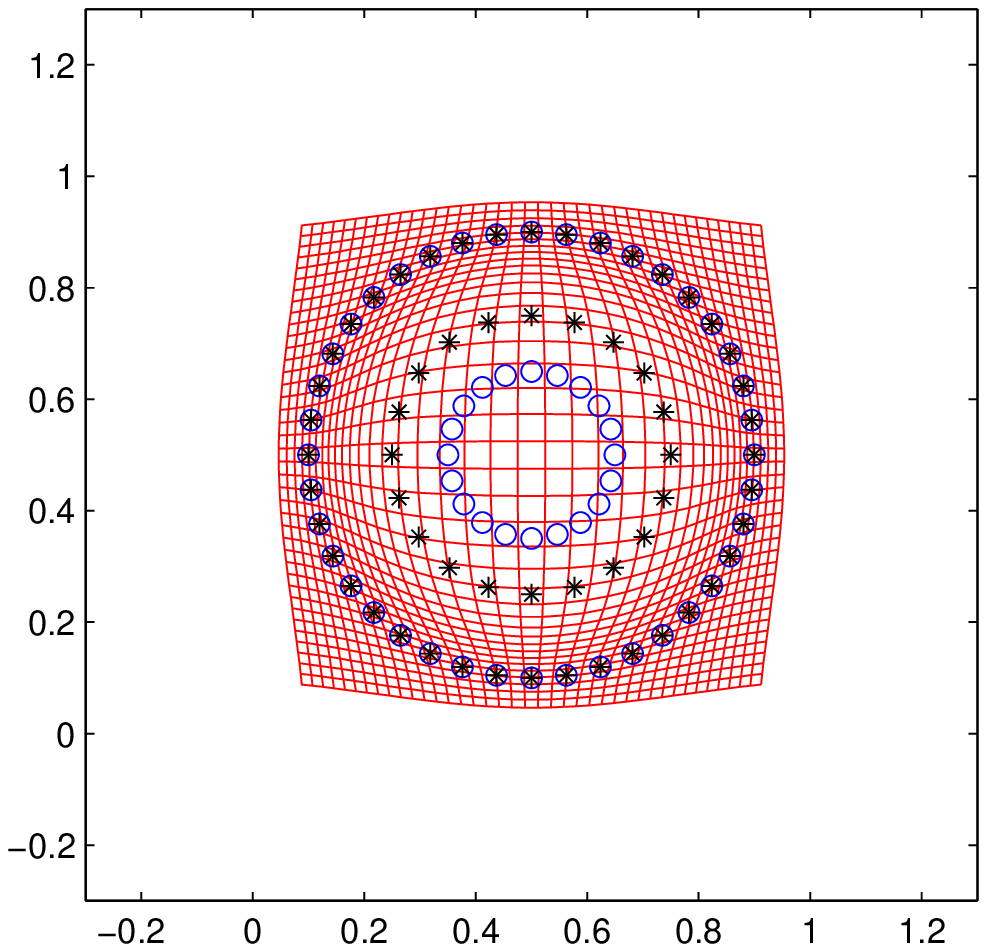}
  \centerline{(a) TPS}
  \end{minipage}
  \begin{minipage}{60mm}
   \includegraphics[width=6.cm]{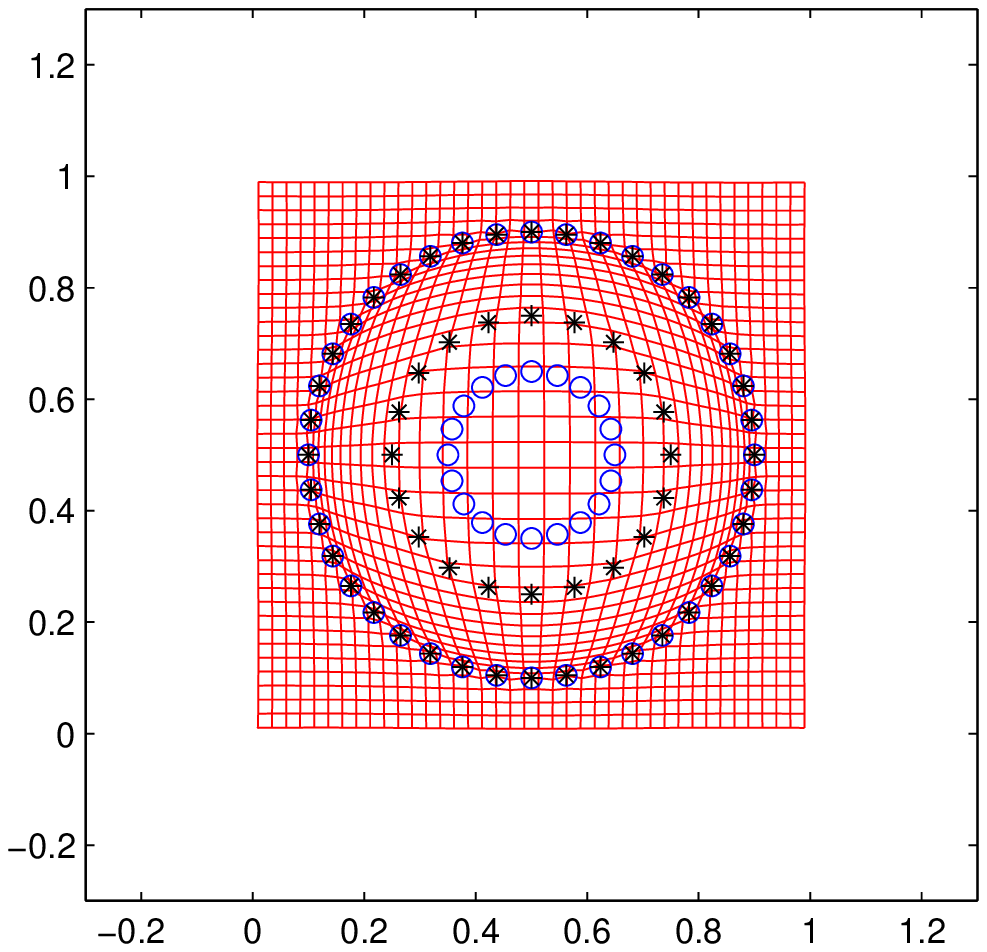}
   \centerline{(b) Shep-TPS}
  \end{minipage}\\
   
\begin{minipage}{60mm}
  \includegraphics[width=6.cm]{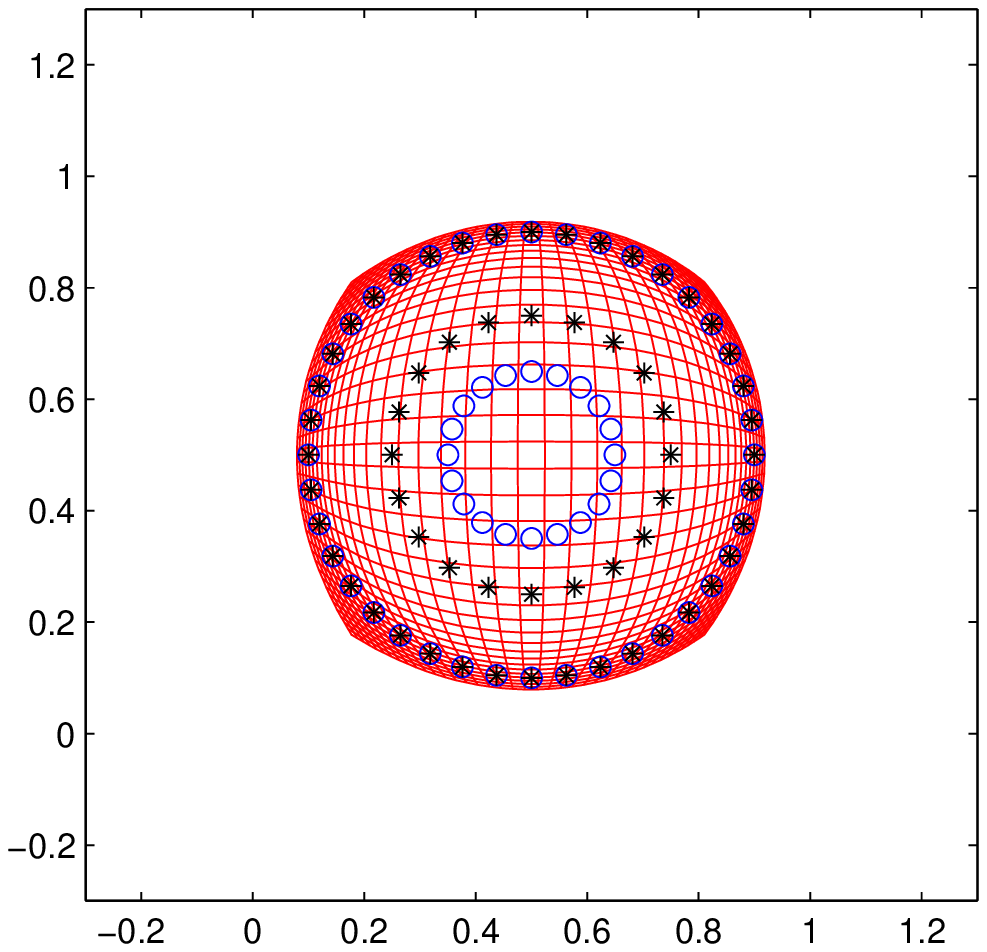}
  \centerline{(c) W2-2D}
  \end{minipage}
  \begin{minipage}{60mm}
   \includegraphics[width=6.cm]{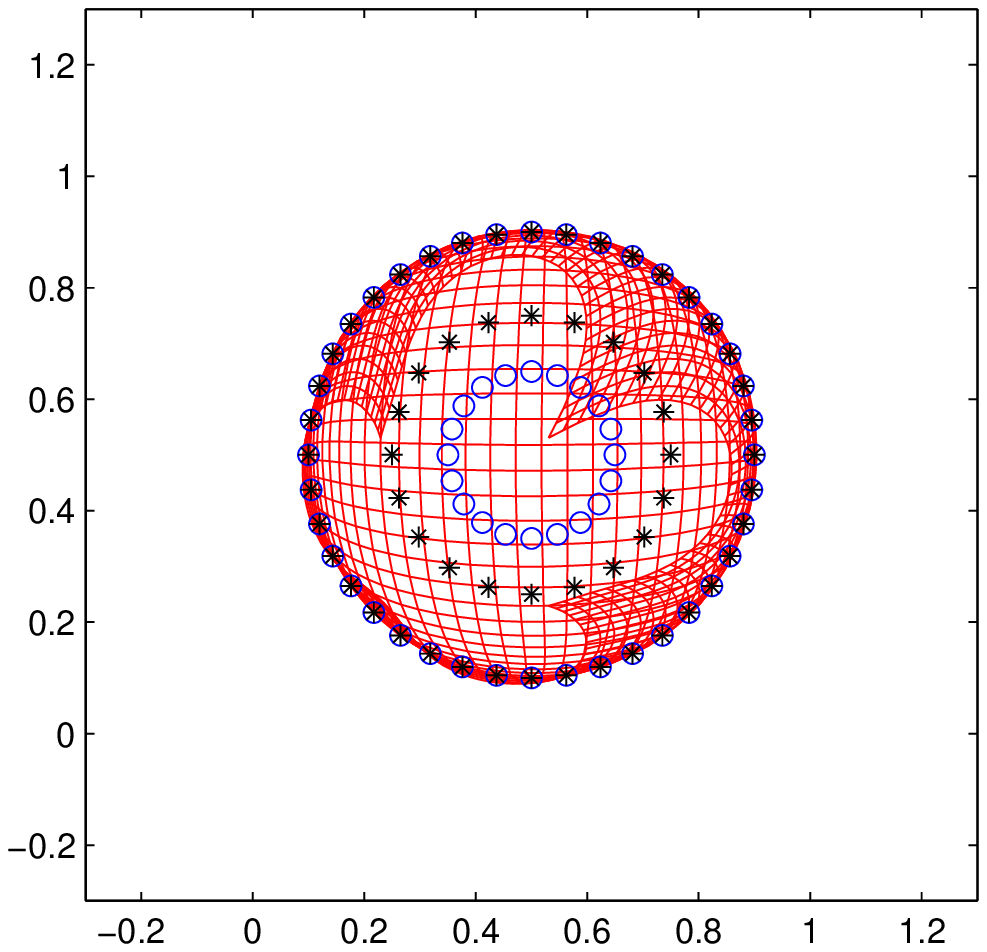}
   \centerline{(d) W2-1D$\times$1D}
  \end{minipage}\\
  
  \begin{minipage}{60mm}
  \includegraphics[width=6.cm]{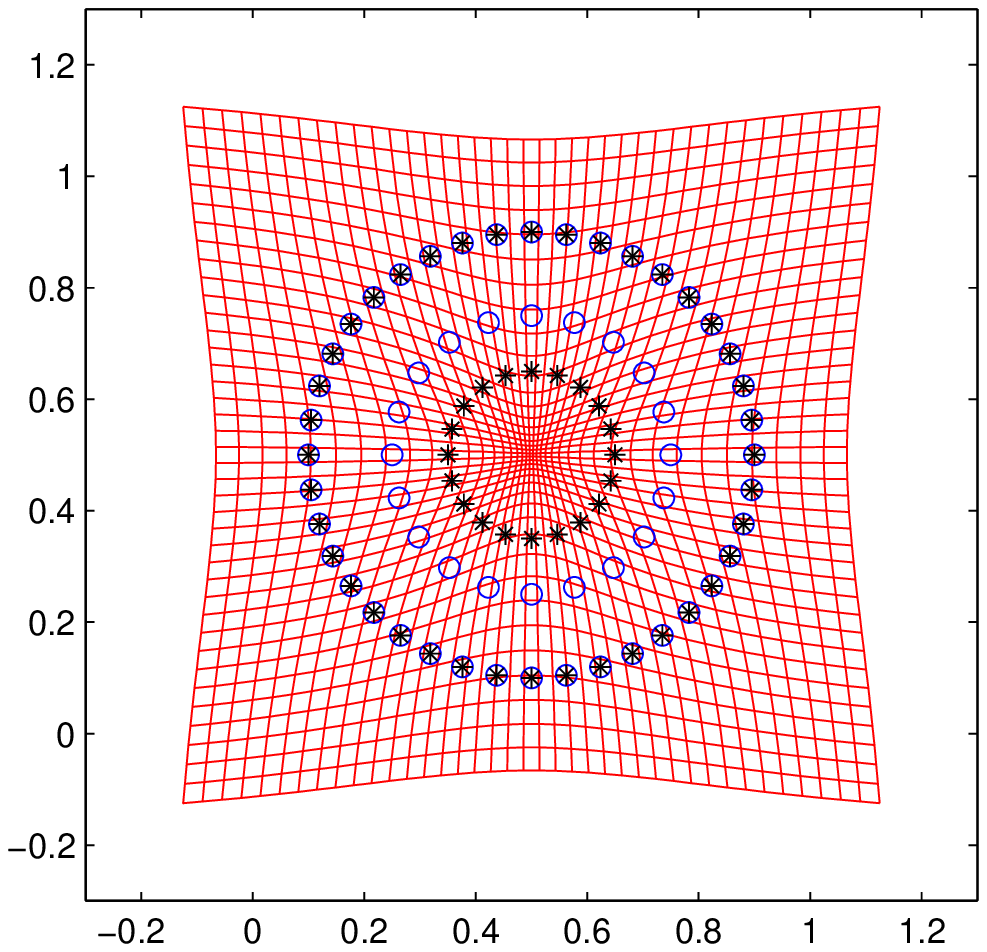}
  \centerline{(e) TPS}
  \end{minipage}
  \begin{minipage}{60mm}
   \includegraphics[width=6.cm]{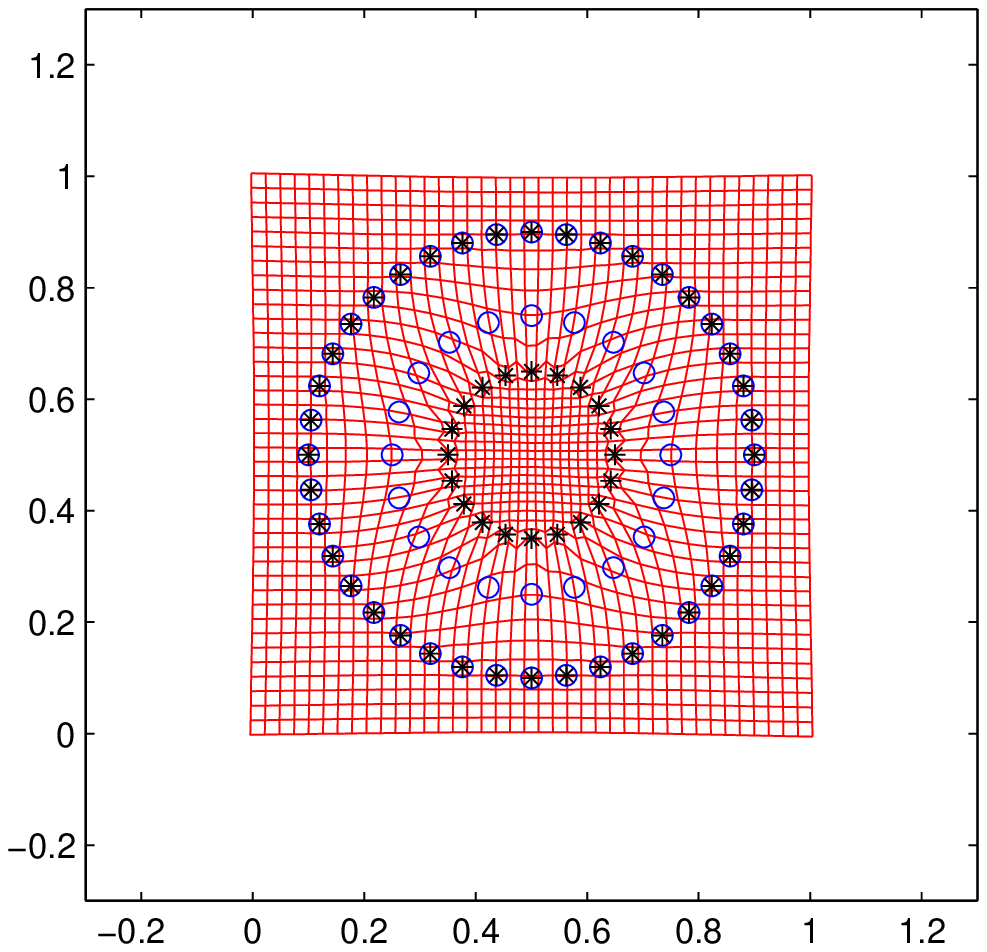}
   \centerline{(f) Shep-TPS}
  \end{minipage}\\

\begin{minipage}{60mm}
  \includegraphics[width=6.cm]{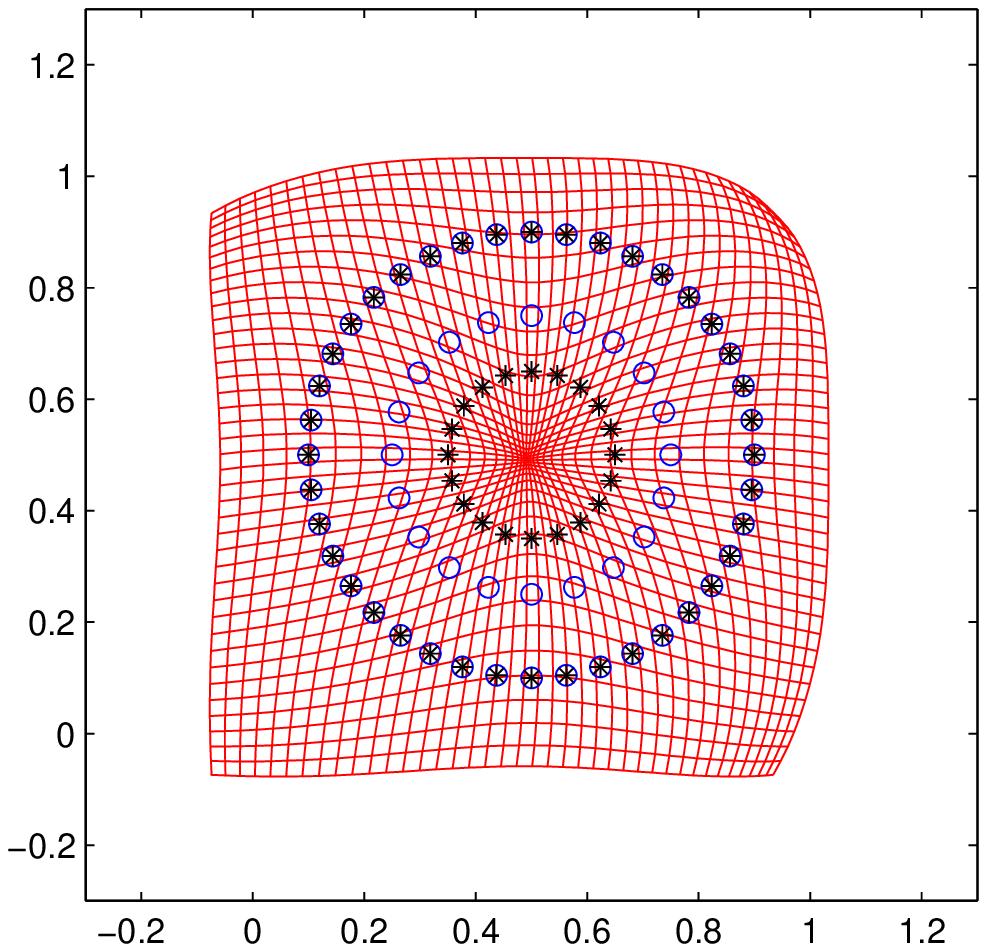}
  \centerline{(g) W2-2D}
  \end{minipage}
  \begin{minipage}{60mm}
   \includegraphics[width=6.cm]{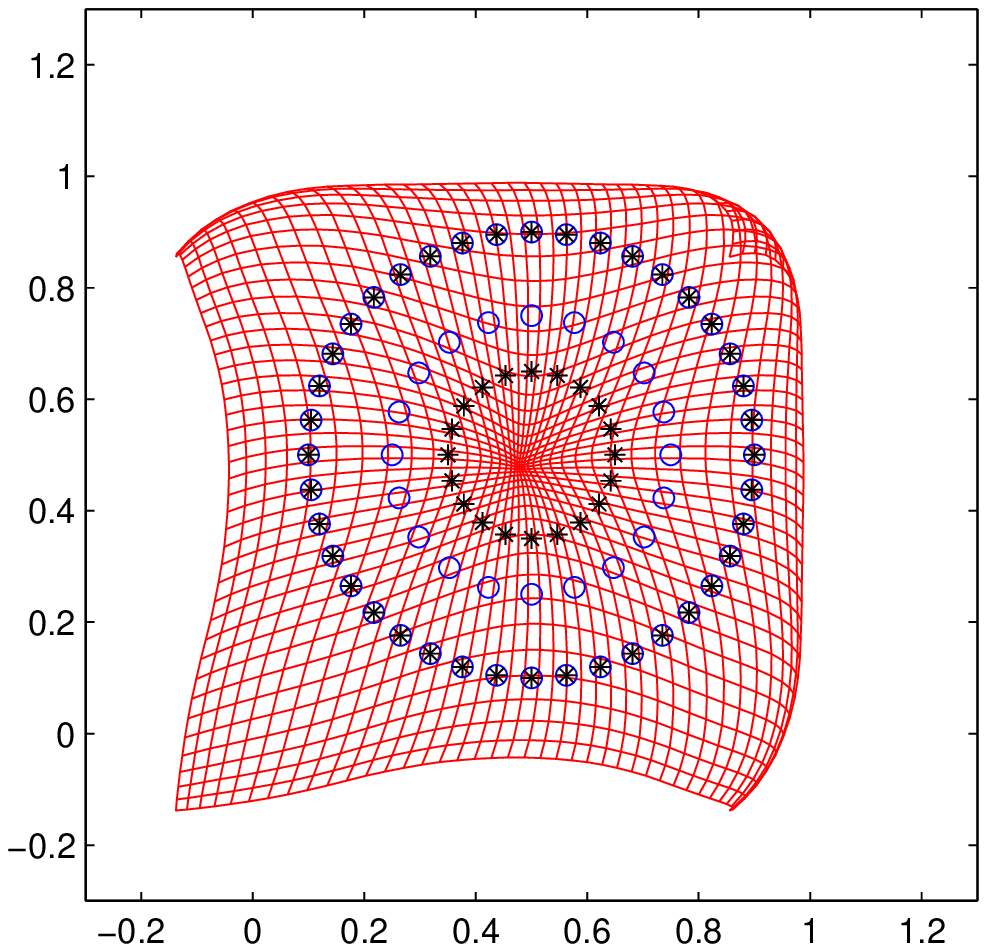}
   \centerline{(h) W2-1D$\times$1D}
  \end{minipage}\\

\end{center}  
  \caption{Circle expansion: registration results obtained by using TPS (a) and Shep-TPS (b), W2-2D (c) and W2-1D$\times$1D (d); circle contraction: registration results obtained by using TPS (e) and Shep-TPS (f), W2-2D (g) and W2-1D$\times$1D (h).}
  \label{fig:14}
\end{figure}

\clearpage

\subsection{Real-life case: an application to medical images}

In this subsection we shortly present some experimental results obtained by applying Gaussians, thin plate splines, Lobachevsky splines, bivariate Wendland's functions and products of univariate Wendland's functions to real image data. More precisely, we consider two X-ray images of the cervical of an anonymous patient taken at different times. The aim of this application is to present a real-life example as, to say, that presented by Modersitzki in \cite{Modersitzki09}, an exaustive discussion based on medical and computer science details being left out. In Figure 15 we show the two images along with landmarks and quasi-landmarks, setting on the left the source image and on the right the target one. The size of both images is $512 \times 512$ pixels. In particular, within each of the two images, there have been selected 6 landmarks, whose coordinates are reported in Table \ref{MR_land}; moreover, to fix transformation and to prevent an overall shift, 12 quasi-landmarks have been added on the boundaries of the source and target images.

\begin{Large}
\begin{table}[ht!]
\begin{center}
\begin{tabular}{|c|c|c|c|c|}
\hline
\rule[0mm]{0mm}{3ex}
 $i$-th landmark &  $x_{i1}$  & $x_{i2}$  &  $t_{i1}$ & $t_{i2}$   \\
\hline
\hline
\rule[0mm]{0mm}{3ex}
$1$ & $0.3135$ & $0.8232$ & $0.3467$ & $0.8525$  \\
\rule[0mm]{0mm}{3ex}
$2$ & $0.3330$ & $0.7080$ & $0.3584$ & $0.7334$  \\
\rule[0mm]{0mm}{3ex}
$3$ & $0.3643$ & $0.5967$ & $0.3877$ & $0.6162$  \\
\rule[0mm]{0mm}{3ex}
$4$ & $0.4131$ & $0.5068$ & $0.4229$ & $0.5068$  \\
\rule[0mm]{0mm}{3ex}
$5$ & $0.4580$ & $0.4053$ & $0.4600$ & $0.3916$  \\
\rule[0mm]{0mm}{3ex}
$6$ & $0.5146$ & $0.3057$ & $0.5205$ & $0.2783$  \\
\hline 
\end{tabular}
\end{center}
\caption{Real-life case: coordinates of source and target landmarks.}
\label{MR_land}
\end{table}
\end{Large}

Each result in Figure 16 (a)--(f) represents a transformed image, obtained using Gaussian, thin plate spline, Wendland's W2-2D, W4-2D, W2-1D$\times$1D, Lobachevsky spline L4 transformations, respectively. For the Gaussian and Lobachevsky L4 transformations we have used the parameter value $\alpha=1.6$, and for Wendland's transformations the value $c=0.1$. We observe that Gaussian transformation strongly deforms the image, while using thin plate and Lobachevsky splines the deformations are less significant. Better registration results are obtained with the Wendland's transformations, especially when we use W2-2D. Moreover, we point out that in general classical Wendland's transformations (both W2-2D and W4-2D) give less deformations than product Wendland's ones, and registration results are better if we use \lq\lq small\rq\rq\ values for the parameter $c$. Finally, in order to provide a more precise analysis among the different methodologies, in Table \ref{MR_err} we also report an error quantification, expressed in terms of the RMSEs and suitably computed on a significant subset of evaluation points automatically chosen. This study confirms once more, from a numerical standpoint, what we have previously observed.

In this application we take into account the analyses performed in Subsection 7.1, especially for the choice of the parameters. Here for the Gaussian and Lobachevsky L4 transformations we have used the parameter value $\alpha=1.6$, and for Wendland's transformations the value $c=0.1$, but, as pointed out in Subsection 7.1, the parameter values for the Wendland and Lobachevsky functions can be chosen in a wide range. This consideration is supported by results obtained with different parameter values. However, it is not astonishing that the optimal parameter values obtained in this example do not coincide with those given in previous test Examples 1 and 2. In fact, small differences of the optimal parameter values in the various situations are to be expected, because the number and collocation of landmarks are considerably changed.

\begin{Large}
\begin{table}[ht!]
\begin{center}
\begin{tabular}{|c|c|c|c|c|c|}
\hline
\rule[0mm]{0mm}{3ex}
 G &  TPS  & W2-2D  &  W4-2D & W2-1D$\times$1D & L4  \\
\hline
\hline
\rule[0mm]{0mm}{3ex}
$1.0314\text{E}-1$ & $1.9685\text{E}-2$ & $1.9526\text{E}-2$ & $2.5981\text{E}-2$ & $2.7157\text{E}-2$ & $1.9843\text{E}-2$ \\
\hline 
\end{tabular}
\end{center}
\caption{Real-life case: RMSEs obtained by using $\alpha=1.6$ and $c=0.1$.}
\label{MR_err}
\end{table}
\end{Large}

\begin{figure}[ht!]
\begin{center}
\includegraphics[width=8.cm]{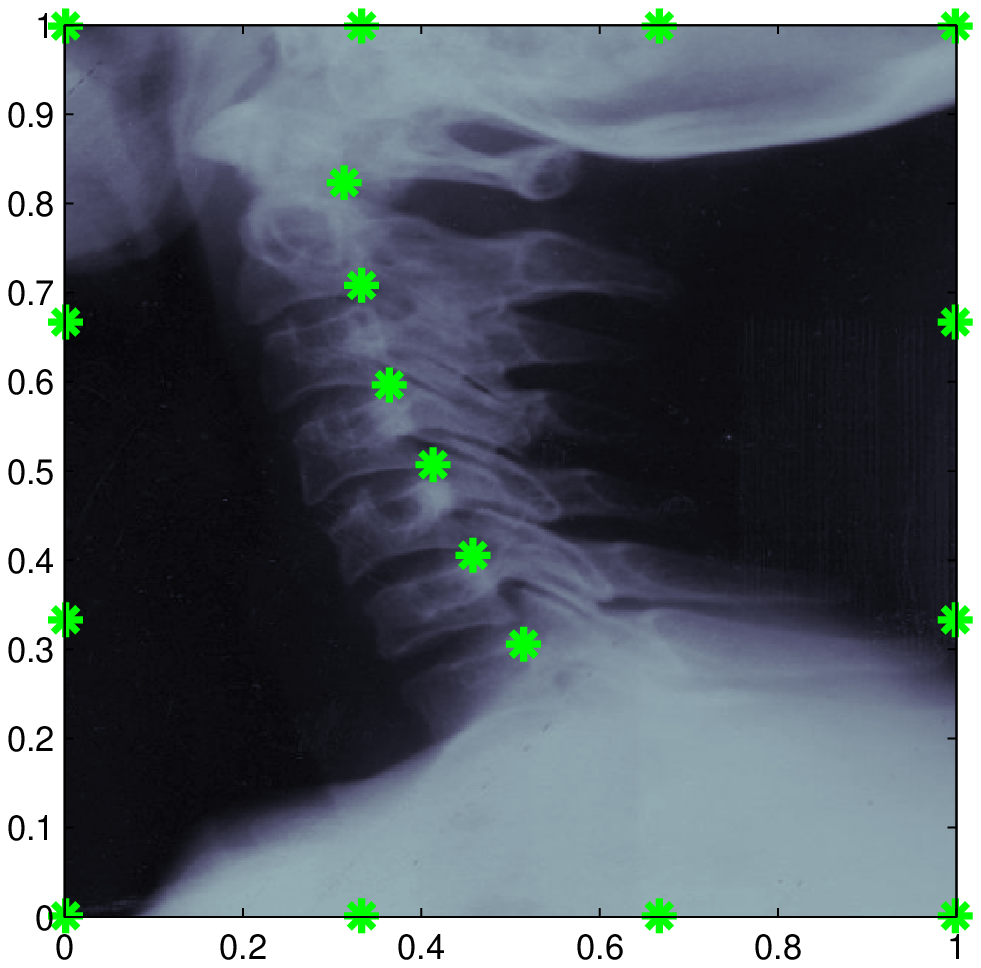}
\includegraphics[width=8.cm]{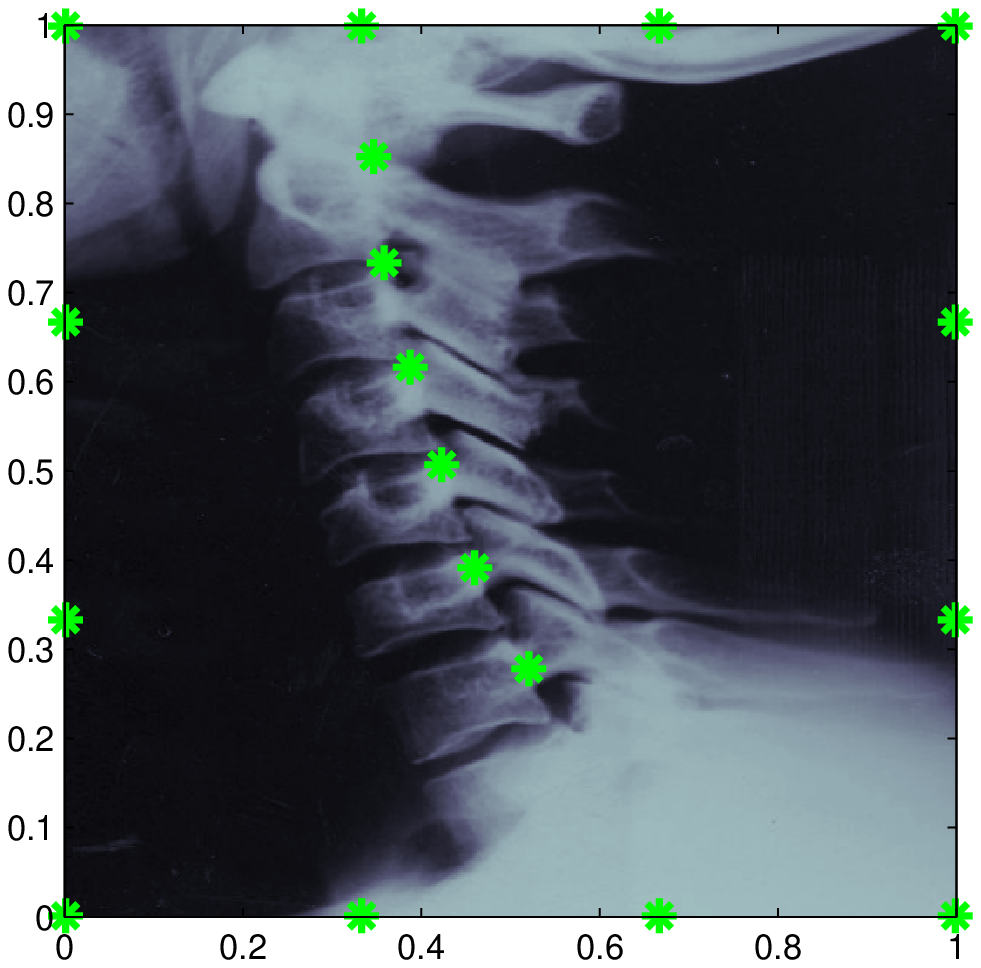}
\end{center}  
  \caption{Source and target cervical images with landmarks and quasi-landmarks (left to right).}
\label{CERV-MI}
\end{figure}

\begin{figure}[ht!]
\begin{center}
\begin{minipage}{80mm}
\includegraphics[width=8cm]{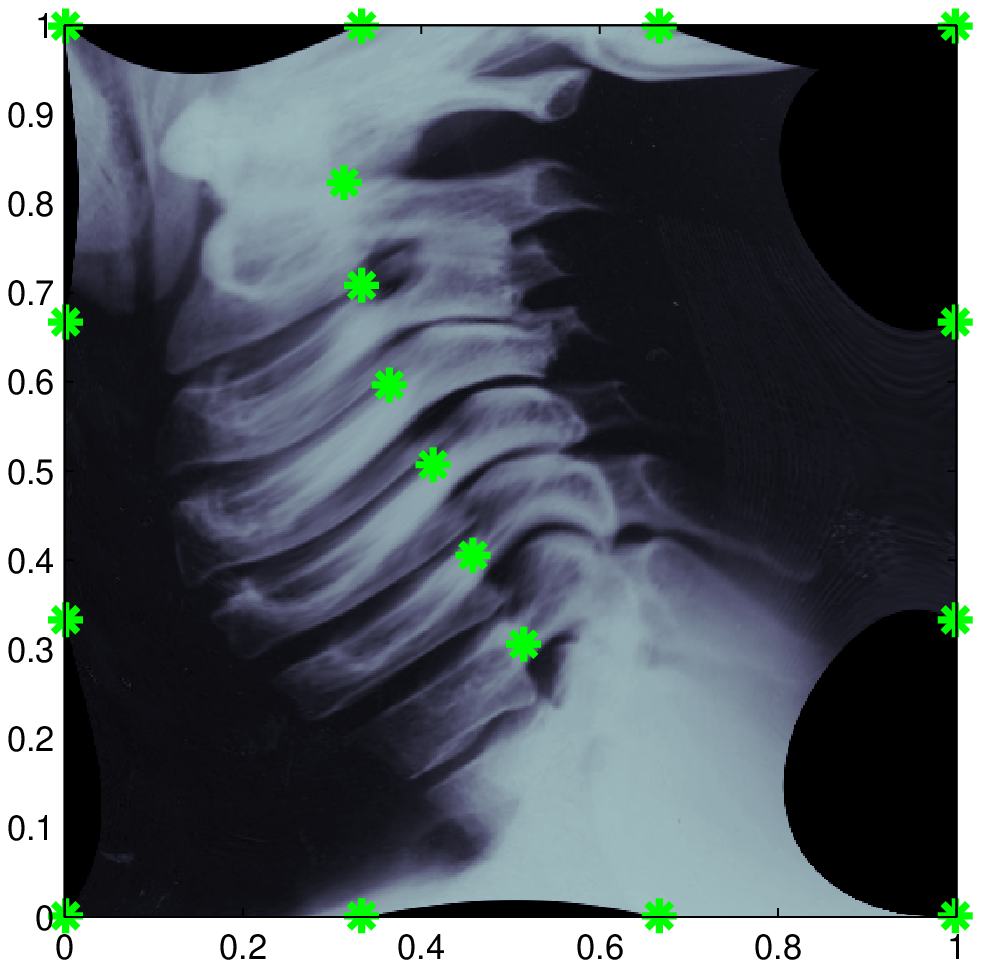}
\centerline{(a) G, $\alpha=1.6$}
\end{minipage}
\begin{minipage}{80mm}
\includegraphics[width=8cm]{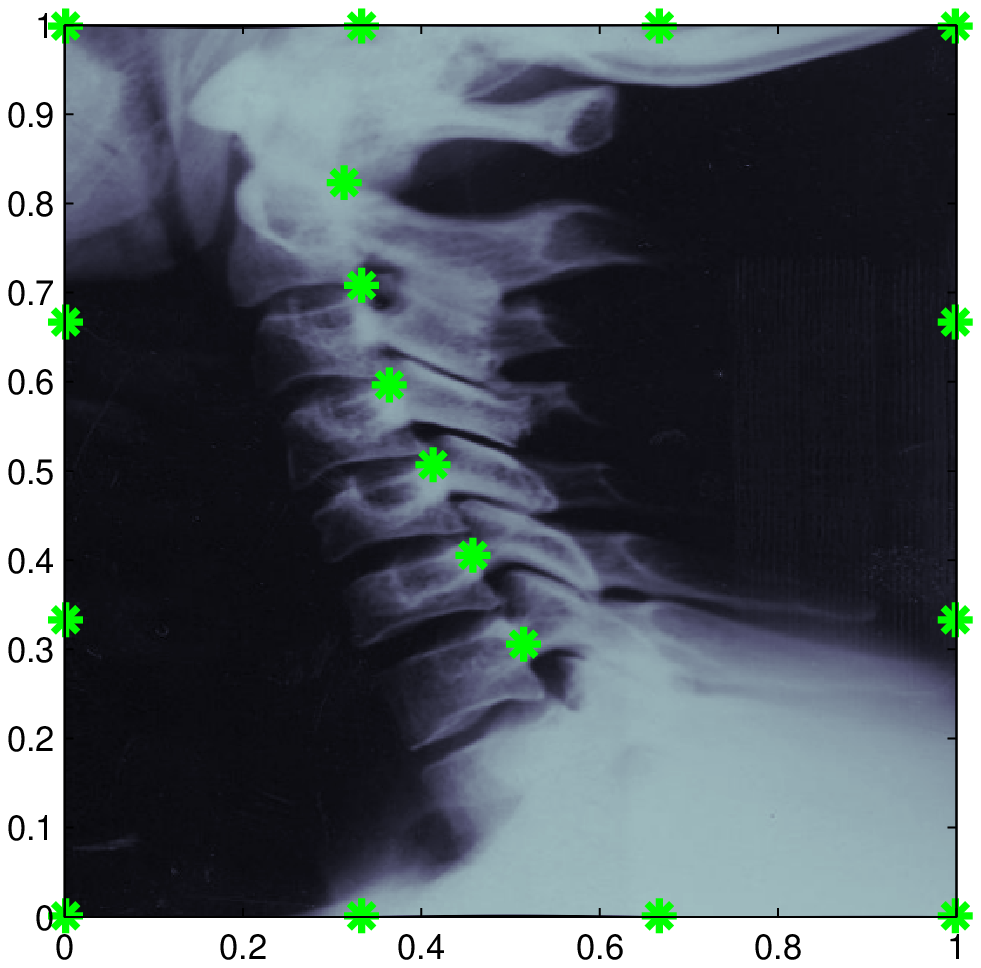}
\centerline{(b) TPS}
\end{minipage}\\
\begin{minipage}{80mm}
\includegraphics[width=8cm]{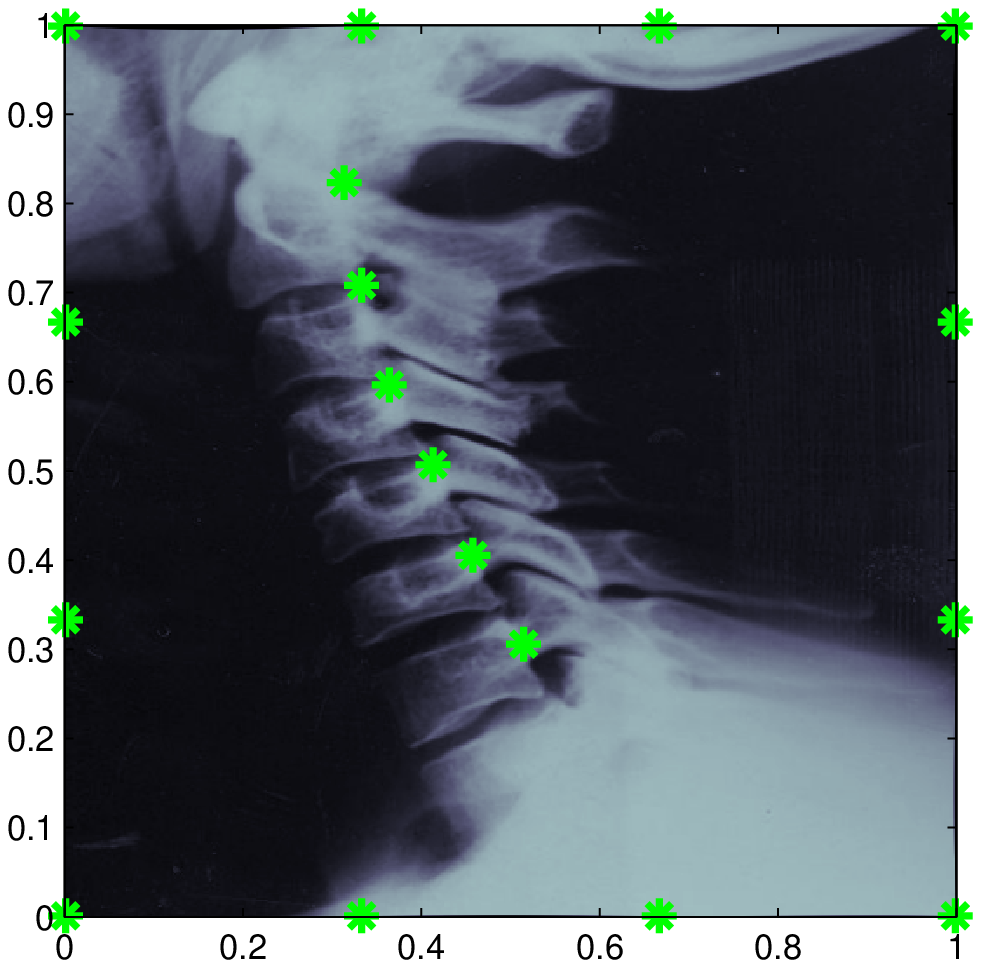}
\centerline{(c) W2-2D, $c=0.1$}
\end{minipage}
\begin{minipage}{80mm}
\includegraphics[width=8cm]{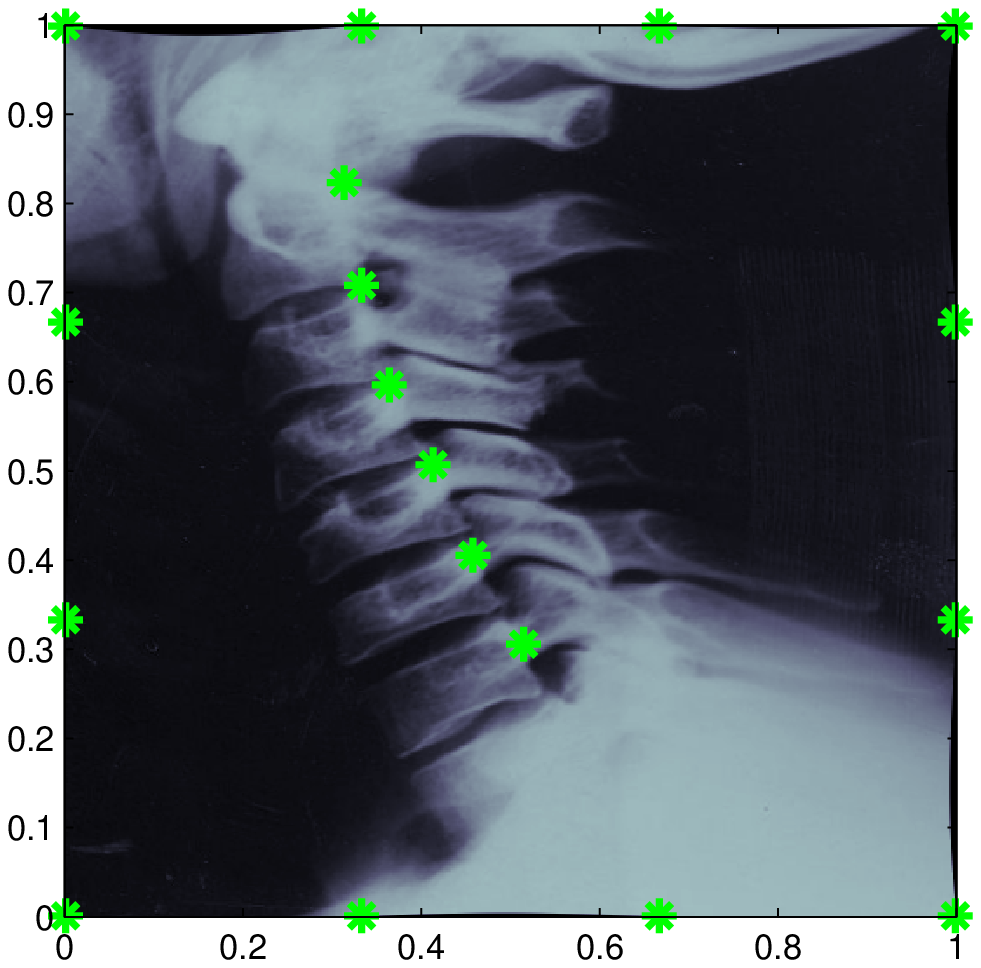}
\centerline{(d) W4-2D, $c=0.1$}
\end{minipage}\\
\begin{minipage}{80mm}
\includegraphics[width=8cm]{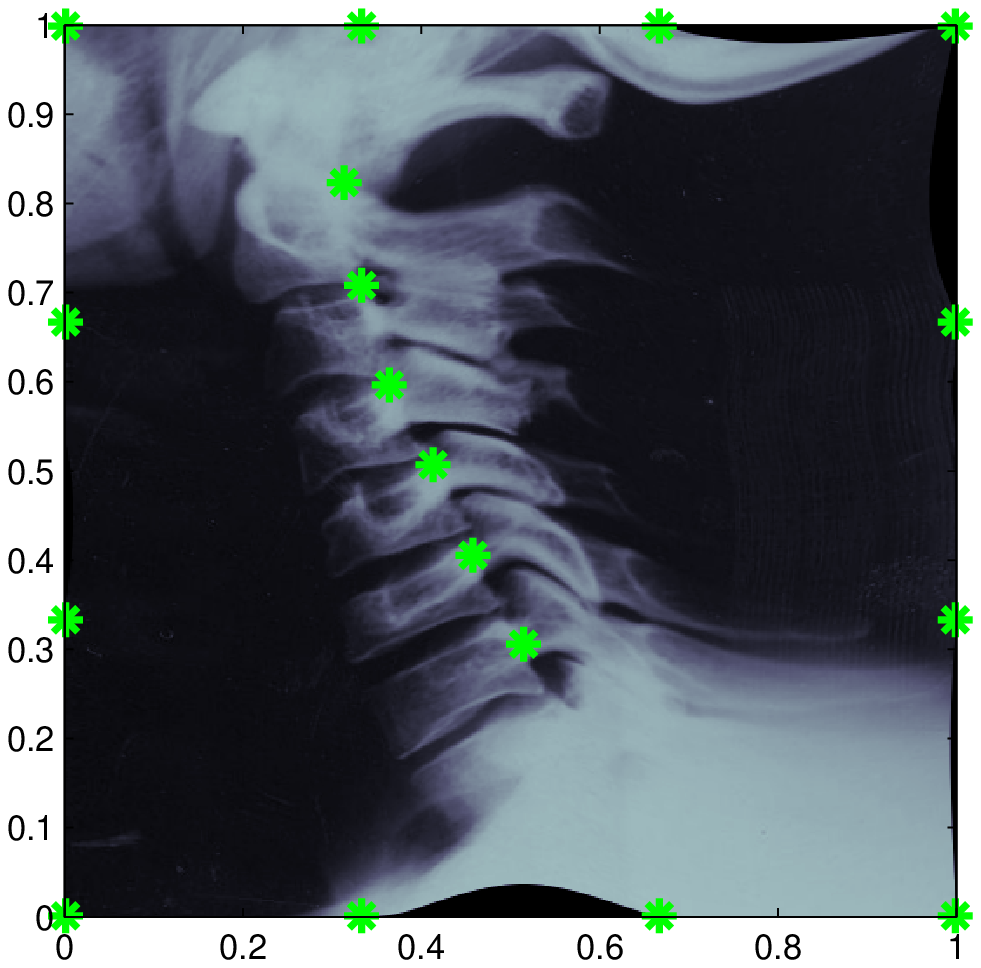}
\centerline{(e) W2-1D$\times$1D, $c=0.1$}
\end{minipage}
\begin{minipage}{80mm}
\includegraphics[width=8cm]{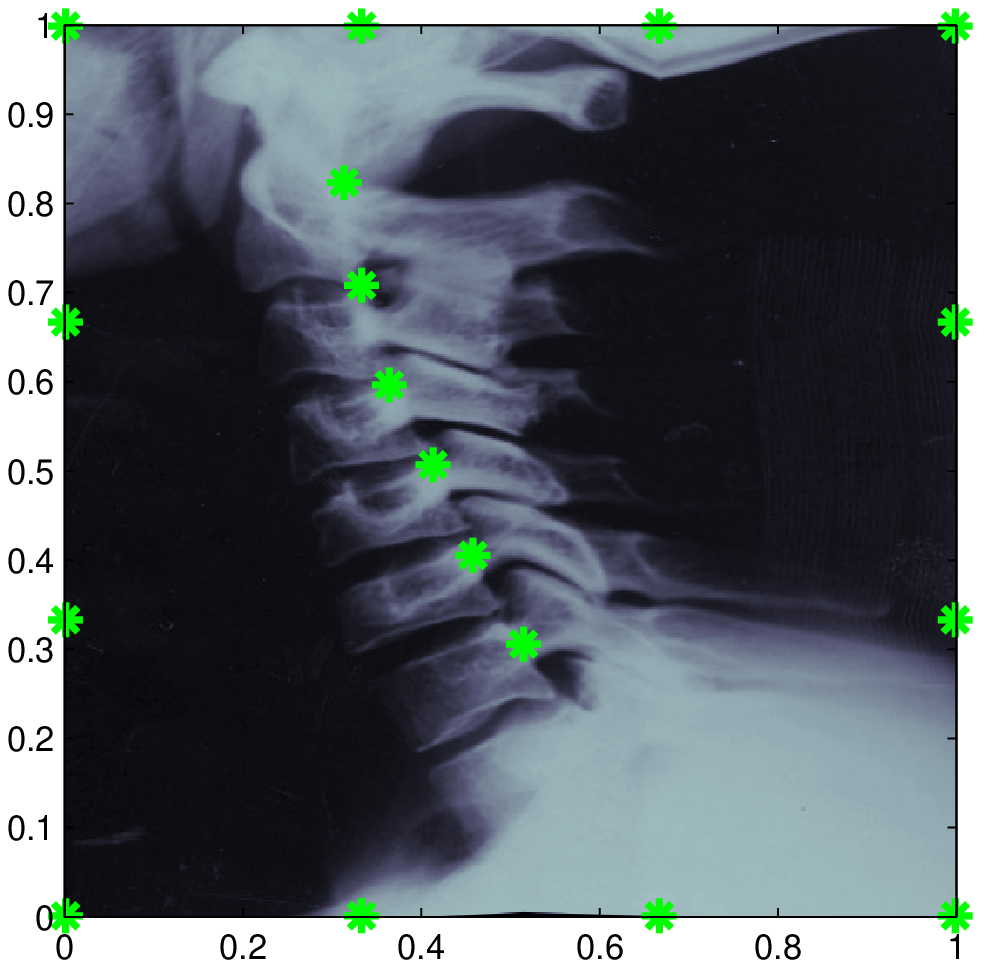}
\centerline{(f) L4, $\alpha=1.6$}
\end{minipage}\\
\end{center}
\caption{Real-life case: registration results obtained by Gaussian (a), thin plate spline (b), Wendland's function (c)-(d)-(e) and Lobachevsky spline (f).}
\label{Ris_cervical}
\end{figure}


\section{Conclusions}

We compared some well established interpolation methods, such as compactly supported radial basis functions, specifically Wendland's functions, and some techniques recently proposed in the context of image registration, namely the modified Shepard's method, in which the nodal functions are obtained by using radial basis functions, and compactly supported Lobachevsky splines. To this aim, we briefly recalled all the above local interpolation schemes and made comments on their performances in landmark-based image registration, taking into account the wide literature on the topic. For the first time, as far as we know, we proposed the use of products of univariate compactly supported radial basis functions, which perform well when they are compared with other schemes even if, in general, multivariate radial Wendland's functions provided better registration results.

In particular, we were interested to point out which scheme may be preferable in a specific situation. Moreover, since many techniques require the use of a shape parameter which might greatly influence registration results, we pointed out that sometimes may be necessary a compromise between accuracy and smoothness, whereas in other situations a significative improving of accuracy, with the use of optimal parameters, is recommendable.

With regard to all the analyzed test cases and, in particular, considering RMSEs and smoothness of the transformed images, we can claim that among the ten considered schemes the following four are preferable: TPS, Shep-TPS, W2-2D, L4.  The TPS method is global and does not need a parameter. This can be advantageous, but as a consequence the scheme is also less adaptive to the different situations that arise in real cases. On the contrary, referring to Shep-TPS, the problem mainly consists in finding optimal choices for the three parameters. Regarding to Wendland's functions W2-2D and Lobachevsky L4, it is remarkable that their performances are good in all the considered different situations and, moreover, the parameter selection can be made by the user in a large interval.

Here, we used interpolating transformations which accomplish an exact match of corresponding landmarks. This implicitly means that the landmark positions are exactly known. However, if we have to deal with landmark localization errors, then it would be advantageous to weaken the interpolation conditions by introducing an approximation scheme. Further investigations in this direction are still required and ongoing. 


\section*{Acknowledgements} 
The authors thank the anonymous referees for their detailed and valuable comments which helped to greatly improve the quality
of the paper.
Moreover, the authors gratefully acknowledge the support of Department of Mathematics \lq\lq G. Peano\rq\rq~, University of Turin, project \lq\lq Modeling and approximation of complex systems (2010)\rq\rq. 
Finally, the second author is grateful to the \lq\lq Istituto Nazionale di Alta Matematica\rq\rq\ (INdAM) for its financial support by a research grant.


\end{document}